\PassOptionsToPackage{dvipsnames}{xcolor} 
\documentclass[11pt,letter]{article}
\usepackage{amssymb,amsmath,amsthm,multirow,color,xcolor,cancel,bbm,tikz,mathrsfs,hyperref,authblk,calc}
\usepackage[textheight=24cm,textwidth=16cm]{geometry}
\tikzstyle{main node}=[draw,circle,inner sep=1,outer sep=2,thick,minimum size=12pt]

\hypersetup{
    colorlinks=true,
    linkcolor=blue,
    citecolor=cyan,
    filecolor=magenta,      
    urlcolor=green,
    pdftitle={Overleaf Example},
    pdfpagemode=FullScreen,
    }

\newtheorem{proposition}{Proposition}
\newtheorem{lemma}{Lemma}
\newtheorem{theorem}{Theorem}

\newtheorem{remark}{Remark}
\newtheorem{example}{Example}

\newcommand{\BF}[1]{{{\boldmath{\bf #1}\unboldmath}}}
\newcommand{\EM}[1]{{\em\textcolor{Maroon}{#1}}}
\newcommand{\B}{\{0,1\}}
\newcommand{\ONE}{\mathbf{1}}
\newcommand{\ZERO}{{\mathbf{0}}}

\newcommand{\q}{{\boldsymbol{q}}}
\newcommand{\W}{{\mathbf{w}}}
\newcommand{\LAMBDA}{\boldsymbol{\lambda}}
\newcommand{\MU}{\boldsymbol{\mu}}
\newcommand{\ELL}{\boldsymbol{\ell}}
\newcommand{\C}{\boldsymbol{c}}
\newcommand{\fix}{\mathrm{fix}}
\newcommand{\synch}{\mathrm{synch}}
\newcommand{\fp}{\mathrm{fp}}
\newcommand{\att}{\mathrm{att}}

\title{Synchronizing Boolean networks asynchronously}

\date{March 10, 2022; revised March 13, 2023}

\author{
Julio Aracena\footnote{\scriptsize CI2MA and Departamento de Ingenier\'ia Matem\'atica, Universidad~de~Concepci\'on, Chile. ({\tt jaracena@ing-mat.udec.cl})},
Adrien Richard\footnote{\scriptsize Universit\'e Côte d’Azur, CNRS, I3S, Sophia Antipolis, France. ({\tt  adrien.richard@cnrs.fr})}~\footnote{\scriptsize Universidad de Chile, CMM, Santiago, Chile.},
Lilian Salinas\footnote{\scriptsize Department of Computer Sciences and CI2MA, University of Concepci\'on, Chile.  Email: ({\tt lilisalinas@udec.cl})}
}

\begin{document}
%%%%%%%%%%%%%%%%%%%%%%%%%%%%%%%%%%%%%%%%%%%%%%%%%%%%%%%%%%%%%%%%%%%%%%%%%%%%%%%%%%%%%%%%%%%%%%%%%%%%%%%
%%%%%%%%%%%%%%%%%%%%%%%%%%%%%%%%%%%%%%%%%%%%%%%%%%%%%%%%%%%%%%%%%%%%%%%%%%%%%%%%%%%%%%%%%%%%%%%%%%%%%%%

\maketitle

\begin{abstract}
The {\em asynchronous automaton} associated with a Boolean network $f:\B^n\to\B^n$, considered in many applications, is the finite deterministic automaton where the set of states is $\B^n$, the alphabet is $[n]$, and the action of letter $i$ on a state $x$ consists in either switching the $i$th component if $f_i(x)\neq x_i$ or doing nothing otherwise. These actions are extended to words in the natural way. A word is then {\em synchronizing} if the result of its action is the same for every state. In this paper, we ask for the existence of synchronizing words, and their minimal length, for a basic class of Boolean networks called and-or-nets: given an arc-signed digraph $G$ on $[n]$, we say that $f$ is an {\em and-or-net} on $G$ if, for every $i\in [n]$, there is $a$ such that, for all state $x$, $f_i(x)=a$ if and only if $x_j=a$ ($x_j\neq a$) for every positive (negative) arc from $j$ to $i$; so if $a=1$ ($a=0$) then $f_i$ is a conjunction (disjunction) of positive or negative literals. Our main result is that if $G$ is strongly connected and has no positive cycles, then either every and-or-net on $G$ has a synchronizing word of length at most $10(\sqrt{5}+1)^n$, much smaller than the bound $(2^n-1)^2$ given by the well known \v{C}ern\'y's conjecture, or $G$ is a cycle and no and-or-net on $G$ has a synchronizing word. This contrasts with the following complexity result:  it is coNP-hard to decide if every and-or-net on $G$ has a synchronizing word, even if $G$ is strongly connected or has no positive cycles. 
\end{abstract}

%%%%%%%%%%%%%%%%%%%%%%%%%%%%%%%%%%%%%%%%%%%%%%%%%%%%%%%%%%%%%%%%%%%%%%%%%%%%%%%%%%%%%%%%%%%%%%%%%%%%%%%
\section{Introduction}
%%%%%%%%%%%%%%%%%%%%%%%%%%%%%%%%%%%%%%%%%%%%%%%%%%%%%%%%%%%%%%%%%%%%%%%%%%%%%%%%%%%%%%%%%%%%%%%%%%%%%%%

A \EM{Boolean network} (BN) is a finite dynamical system usually defined by a function  
\[
f:\B^n\to\B^n,\qquad x=(x_1,\dots,x_n)\mapsto f(x)=(f_1(x),\dots,f_n(x)).
\]

BNs have many applications. In particular, since the seminal papers of McCulloch and Pitts \cite{MP43}, Hopfield \cite{H82}, Kauffman \cite{K69,K93} and Thomas \cite{T73,TA90}, they are omnipresent in the modeling of neural and gene networks (see \cite{B08,N15} for reviews). They are also essential tools in computer science, see \cite{ANLY00,GRF16,BGT14, CFG14a, GR15b} for instance. 

\medskip
The ``network'' terminology comes from the fact that the \EM{interaction digraph} of $f$ is often considered as the main parameter of $f$: the vertex set is $[n]=\{1,\dots,n\}$ and there is an arc from $j$ to $i$ if $f_i$ depends on input $j$. The \EM{signed interaction digraph} provides useful additional information about interactions, and is commonly consider in the context of gene networks: the vertex set is $[n]$ and there is a positive (negative) arc from $j$ to $i$ if there is $x,y\in \B^n$ that only differ in $x_j<y_j$ such that $f_i(y)-f_i(x)$ is positive (negative). Note that the presence of both a positive and a negative arc from one vertex to another is allowed. If $G$ is the signed interaction digraph $f$ then we say that $f$ is a BN \EM{on} $G$. 

\medskip
From a dynamical point of view, the successive iterations of $f$ describe the so called synchronous dynamics: if $x^t$ is the state of the system at time $t$, then $x^{t+1}=f(x^t)$ is the state of the system at the next time. Hence, all components are updated in parallel at each time step. However, when BNs are used as models of natural systems, such as gene networks, synchronicity can be an issue. This led researchers to consider the (fully) asynchronous dynamics, where one component is updated at each time step (see e.g. \cite{T91,TA90,TK01,A-J16}). This asynchronous dynamics can be described by the paths of a deterministic finite automaton called \EM{asynchronous automaton} of $f$: the set of states is $\B^n$, the alphabet is $[n]$ and   the action of letter $i$ on a state $x$, denoted $f^i(x)$, is the state obtained from $x$ by updating the component $i$ only, that is, $f^i(x)=(x_1,\dots,f_i(x),\dots,x_n)$. These actions are extended to any word $w=i_1,\dots,i_\ell$ over the alphabet $[n]$ in the natural way, by setting $f^w=f^{i_\ell}\circ f^{i_{\ell-1}}\circ\dots\circ f^{i_1}$.

\medskip
In this paper, we study synchronizing properties of this asynchronous automaton, as proposed in \cite{AGRS20}. A word $w$ over $[n]$ \EM{synchronizes} $f$ if it synchronizes its asynchronous automaton, that is, if $f^w$ is a constant function. If $w$ synchronizes $f$ then $w$ is a \EM{synchronizing word} for $f$, and if $f$ admits a synchronizing word then $f$ is \EM{synchronizing}. The central open problem concerning synchronization is the famous \v{C}ern\'y's conjecture \cite{C64,CPR71}, which says that any synchronizing deterministic automaton with $q$ states has a synchronizing word of length at most $(q-1)^2$. This conjecture has been proved for several classes of automata, but the best general bounds are only cubic in $q$, see \cite{V08} for a review.

\medskip
A common research direction concerning BNs tries to deduce from a signed digraph $G$ on $[n]$ the dynamical properties of the BNs $f$ on $G$. An influential result is this direction is the following \cite{A08}: {\em (i)} if $G$ has no negative cycles, then $f$ has at least one fixed point, and {\em (ii)} if $G$ has no positive cycles, then $f$ has at most one fixed point. It is then natural to follow this line of research, and ask what can be said on synchronizing properties under these hypothesis. 

\medskip
First, suppose that $G$ has no negative cycles. If $f$ has at least two fixed points then it is not synchronizing. Otherwise, by {\em (i)}, $f$ has a unique fixed point  and it is not difficult to show that $f$ has a synchronizing word of length $n$ (see Proposition~\ref{pro:positive_G} in Appendix \ref{app:one_fixed_point}). This completely solves the case where $G$ has no negative cycles. 

\medskip
Second, suppose that $G$ has no positive cycles. By {\em (ii)}, $f$ has at most one fixed point. If $f$ has indeed a fixed point, then one can show, as previously, that $f$ has a synchronizing word of length $n$ (see Proposition~\ref{pro:negative_G} in Appendix \ref{app:one_fixed_point}). However, what happens if $f$ has no fixed points? Here some difficulties come and our main result provides a partial answer. 

\medskip
Suppose that $G$ has no positive cycles and is, in addition, strongly connected and non-trivial (that is, contains at least one arc). Then $f$ has no fixed points \cite{A08}, which is the interesting case mentioned above, thus these additional assumptions are natural. It is still difficult to understand synchronizing properties, but our main result, Theorem~\ref{thm:main} below, gives a clear picture when $f$ belongs to a well studied class of BNs, called and-or-nets (see e.g. \cite{GH00,ADG04b,CLP2005,MTA10,VL12,RR13,VAL15,ARS17b} and the references therein for studies about this class of BNs). 

\medskip
We say that $f$ is an \EM{and-or-net} if, for every $i\in [n]$, the component $f_i$ of $f$ is a conjunction or a disjunction of positive or negative literals, that is, there is $a\in \B$ such that, for every state $x$, we have $f_i(x)=a$ if and only if $x_j=a$ for all positive arcs of $G$ from $j$ to $i$ and $x_j\neq a$ for all negative arcs of $G$ from $j$ to $i$ (so if $a=1$ then $f_i$ is a conjunction, and if $a=0$ then $f_i$ is a disjunction). We can now state our main result.

%\cite{GH00,ADG04b,CLP2005,TA09,TA09b,JLV10,MTA10,VL12,GN12,AKMT12,VBHKWL12,RR13,ARS14,VAL15,ARS17b,TS20}

\begin{theorem}\label{thm:main}
Let $G$ be a strongly connected signed digraph on $[n]$ without positive cycles. Either
\begin{itemize}
\item $G$ is a cycle and no BN on $G$ is synchronizing, or
\item every and-or-net on $G$ has a synchronizing word of length at most $10(\sqrt{5}+1)^n$.
\end{itemize}
\end{theorem}

Since $\sqrt{5}+1<4$, the bound $10(\sqrt{5}+1)^n$ is sub-quadratic according to the number $2^n$ of states, and thus much smaller than the bound $(2^n-1)^2$ predicted by \v{C}ern\'y's conjecture. 

\medskip
The difficult part is the second case. For that, the main tool is the following lemma, which should be of independent interest since it holds for every BNs, and not only for and-or-nets. It says that if $G$ is non-trivial, strongly connected and without positive cycles, then no components are definitively fixed in the asynchronous dynamics. 

\begin{lemma}\label{lem:flipping_lemma}
Let $G$ be a non-trivial strongly connected signed digraph on $[n]$ without positive~cycles, and let $f$ be a BN on $G$. For all $i\in [n]$ and $x\in\B^n$, we have $f^w(x)_i\neq x_i$ for some word~$w$.
\end{lemma}

Let $G$ be as in Theorem~\ref{thm:main} and non-trivial. To obtain the bound $10(\sqrt{5}+1)^n$, we prove that if $f$ is an and-or-net on $G$, the word $w$ in the previous lemma is of length at most $F(n+2)$, the $(n+2)$th Fibonacci number. From that, we next prove that if $G$ is not a cycle, then for every states $x,y$ there is a word $w$ of length at most $3F(n+4)$ such that $f^w(x)=f^w(y)$. This immediately implies that $f$ has a synchronizing word of length at most $3F(n+4)(2^n-1)$ and an easy computation shows that this is less than $10(\sqrt{5}+1)^n$. However, we will show that at least one and-or-net on $G$ can be synchronized much more quickly, with a word of sub-linear length according to the number of states. 

\begin{theorem}\label{thm:fast}
Let $G$ be a strongly connected signed digraph on $[n]$ without positive cycles, which is not a cycle. At least one and-or-net on $G$ has a synchronizing word of length at most $5n(\sqrt{2})^n$. 
\end{theorem}

Concerning complexity issues, by Theorem~\ref{thm:main}, if $G$ is strongly connected and has no positive cycles, then one can decide in linear time if every and-or-net on $G$ is synchronizing. This contrasts with the following results, which shows that if one of the two hypotheses made on $G$ (strongly connected, no positive cycles) is removed, the decision problem becomes much more harder. 

\begin{theorem}\label{thm:complexity}
Let $G$ be a signed digraph on $[n]$. If $G$ is not strongly connected or has a positive cycle, then it is coNP-hard to decide if every and-or-net on $G$ is synchronizing (even if $G$ has maximum in-degree at most $2$ and a vertex meeting every cycle). 
\end{theorem}

Finally, the following (easy) property shows that the second case of Theorem~\ref{thm:main} cannot be generalized to all the BNs on $G$. We say that $G$ is \EM{simple} if it does not contains both a positive and negative arc from one vertex to another (if $G$ is strongly connected and has no positive cycles then it is necessarily simple).  

\begin{proposition}\label{pro:02}
If $G$ is a simple signed digraph without vertices of in-degree $0$ or $2$, then at least one BN on $G$ is not synchronizing.
\end{proposition}  

This proposition also shows that conjunctions and disjunctions are in some sense necessary, since if $G$ is simple and $i$ is a vertex of in-degree at most $2$, then $f_i$ is either a conjunction or a disjunction for every BN $f$ on $G$.

\begin{figure}[p]
\[
\begin{tikzpicture}[scale=0.85]
\tikzstyle{every node}=[font=\small,rectangle,fill=white,line width=0.3mm,rounded corners=0.2cm,inner sep=3,outer sep=2,anchor=west,draw=Red!80,text=Red!80]
\node[text width=73] 	at (-0.5,0) 	(same) 		{$G$ has no cycles of opposite signs};
\node 					at (-0.5,12) 	(no-) 		{$G$ has no neg. cycles};
\node[text width=57]	at (-0.5,9) 	(strong) 	{$G$ is strongly connected};
\node 					at (-0.5,6) 	(no+) 		{$G$ has no pos. cycles};
\node 					at (-0.5,4)		(nocy)		{$G$ is not a cycle};
\node[anchor=center]	at (14.3,4)		(cy)		{$G$ is a cycle};
\node					at (-0.5,-2)	(acy)		{$G$ is a acyclic};
\tikzstyle{every node}=[font=\small,rectangle,fill=white,line width=0.3mm,rounded corners=0.2cm,inner sep=3,outer sep=2,anchor=west,draw=Blue!70,text=Blue!70]
\node[anchor=center]	at (4.7,-2) 	(fp1) 		{$\fp(f)=1$};
\node 					at (6.7,-2) 	(fix=synch) 		{$\fix(f)=\synch(f)$};
\node 					at (6.7,0) 		(synchn) 	{$\synch(f)\leq n$};
\node 					at (-0.5,2)		(andor)		{$f$ is an and-or-net};
\node 					at (6.7,14) (n3) 		{$\fix(f)=O(n^3)$};
\node 					at (6.7,12) (attfp) 	{$\att(f)=\fp(f)$};
\node 					at (6.7,6) 	(att-1) 	{$\att(f)= 1$};
\node 					at (6.7,8) 	(unstable) 	{$f$ is unstable};
\node 					at (6.7,2)	(bound1)	{$\synch(f)\leq 10(\sqrt{5}+1)^n$};
\node[text width=105] 	at (6.7,4)	(bound2)	{$\synch(h)\leq 5n(\sqrt{2})^n$ for some and-or-net $h$ on $G$};
\node at (10.2,14) 	(fixable) 	{$\fix(f)<\infty$};
\node at (10.2,12) 	(fp+1) 		{$\fp(f)\geq 1$};
\node at (10.2,10) 	(fp+2) 		{$\fp(f)\geq 2$};
\node at (10.2,6) 	(fp-1) 		{$\fp(f)\leq 1$};
\node at (10.2,8) 	(fp0) 		{$\fp(f)=0$};
\node[anchor=center] at (14.3,10)	(att+2) 	{$\att(f)\geq 2$};
\node[anchor=center] at (14.3,2)		(synch) 	{$\synch(f)<\infty$};
\node[anchor=center] at (14.3,6)	(nosynch) 	{$\synch(f)=\infty$};
\node[outer sep=0,draw=black,text=black,circle,anchor=center,anchor=center] 		at (4.7,10) (and1) {\tiny$\land$};
\node[outer sep=0,draw=black,text=black,circle,anchor=center,anchor=center] 		at (4.7,0)	(and2) {\tiny$\land$};
\node[outer sep=0,draw=black,text=black,circle,anchor=center,anchor=center] 		at (4.7,8) 	(and4) {\tiny$\land$};
\node[outer sep=0,draw=black,text=black,circle,anchor=center,anchor=center] 		at (4.7,4) 	(and5) {\tiny$\land$};
\node[outer sep=0,draw=black,text=black,circle,anchor=center,anchor=center] 		at (4.7,2) 	(and6) {\tiny$\land$};
\tikzstyle{every node}=[inner sep=0,outer sep=3,font=\footnotesize,anchor=center]

\path[->,thick]
(no-.east) 		edge[black] node[pos=0.8,above] {\cite{R10}} (attfp)
(no-.east) 		edge[black] node[pos=0.8,above left] {\cite{AGRS20}} (n3.west)
(no-.east) 		edge[black,-] (and1)
(no-.east) 		edge[bend right=20] node[pos=0.85,below] {\cite{A08}} (fp+1)
(no+.east) 		edge[ultra thick,-] (and4)
(no+.east) 		edge[bend left=20] node[pos=0.85,above] {\cite{A08,aracena2021finding}} (fp-1)
(no+.east) 		edge node[pos=0.8,below] {\cite{RC07}} (att-1)
(strong.east) 	edge[black,-] (and1)
(strong.east) 	edge[ultra thick,-] (and4)
(fixable) 		edge[Blue!70,<->] (attfp)
(bound1) 		edge[Blue!70] (synch)
(synch.160)		edge[Blue!70,bend right=25] (att-1)
(unstable) 		edge[Blue!70] (fp0)
(acy)	 		edge node[pos=0.5,above] {\cite{R80}} (fp1)
(fp0) 			edge[Blue!70] (fp-1)
(fp+2) 			edge[Blue!70] (att+2)
(fp+2) 			edge[Blue!70] (fp+1)
(fp1) 			edge[ultra thick,black,-] (and2)
(fp1) 			edge[Blue!70,->] (fix=synch)
(n3) 			edge[Blue!70] (fixable)
(attfp) 		edge[Blue!70] (fp+1)
(att+2) 		edge[Blue!70] (nosynch)
(att-1) 		edge[Blue!70] (fp-1)
(same) 			edge[ultra thick,black,-] (and2)
(nocy) 			edge[ultra thick,-] (and5)
(cy) 			edge[ultra thick] node[pos=0.5,right,outer sep=3.8] {{\bf Lem. \ref{lem:initial_cycle}}} (nosynch)
(andor) 		edge[ultra thick,-] (and6)
(synchn.0) 		edge[Blue!70] (synch.190)
(and1) 			edge [black] node[pos=0.79,above] {\cite{A08}} (fp+2)
(and2) 			edge[ultra thick] node[pos=0.5,above,outer sep=3.8] {{\bf Pro. \ref{pro:positive_G}}} node[pos=0.5,below,outer sep=3.8] {{\bf Pro. \ref{pro:negative_G}}} (synchn)
(and4) 			edge[ultra thick] node[pos=0.5,below,outer sep=3.8] {{\bf Lem. \ref{lem:flipping_lemma}}} (unstable)
(and4) 			edge[bend left=25] node[pos=0.785,above] {\cite{A08}} (fp0)
(and4) 			edge[ultra thick,->] (and5)
(and5) 			edge[ultra thick] node[pos=0.5,below,outer sep=3.8] {{\bf Thm. \ref{thm:fast}}} (bound2)
(and5) 			edge[ultra thick,->] (and6)
(and6) 			edge[ultra thick] node[pos=0.5,below,outer sep=3.8] {{\bf Thm. \ref{thm:main}}} (bound1)
;
\draw[->,thick,Red!80] (no+.180) -- (-0.8,6) -- (-0.8,0) -- (same.180);
\draw[->,thick,Red!80] (no-.180) -- (-1.2,12) -- (-1.2,-0.24) -- (same.188);
\end{tikzpicture}
\]
\caption{\label{fig:summary}
Summary of the main results and related results from the literature. $G$ is a signed digraph on $[n]$, and $f$ is a BN on $G$. $\fp(f)$ is the number of fixed points of $f$. $\att(f)$ is the number of asynchronous attractors of $f$, that is, the number of inclusion-minimal set $A\subseteq\B^n$ such that $f^i(A)\subseteq A$ for every $i\in [n]$. We have $\att(f)\geq 1$ and $\att(f)\geq \fp(f)$ since $x$ is a fixed point of $f$ if and only if $\{x\}$ is an asynchronous attractor of $f$. $f$ is unstable if, for all $i\in [n]$ and $x\in\B^n$, we have $f^w(x)_i\neq x_i$ for some word~$w$ (conclusion of Lemma~\ref{lem:flipping_lemma}). $\fix(f)$ is the minimum length of a word $w$ over $[n]$ such that $f^w(x)$ is a fixed point of $f$ for any $x\in\B^n$; if no such word exists then $\fix(f)=\infty$. $\synch(f)$ is the minimum length of a synchronizing word for $f$; if no such word exists then $\synch(f)=\infty$. %One easily check that if $\fp(f)=1$ then $\synch(f)=\fix(f)$. 
}
\end{figure}
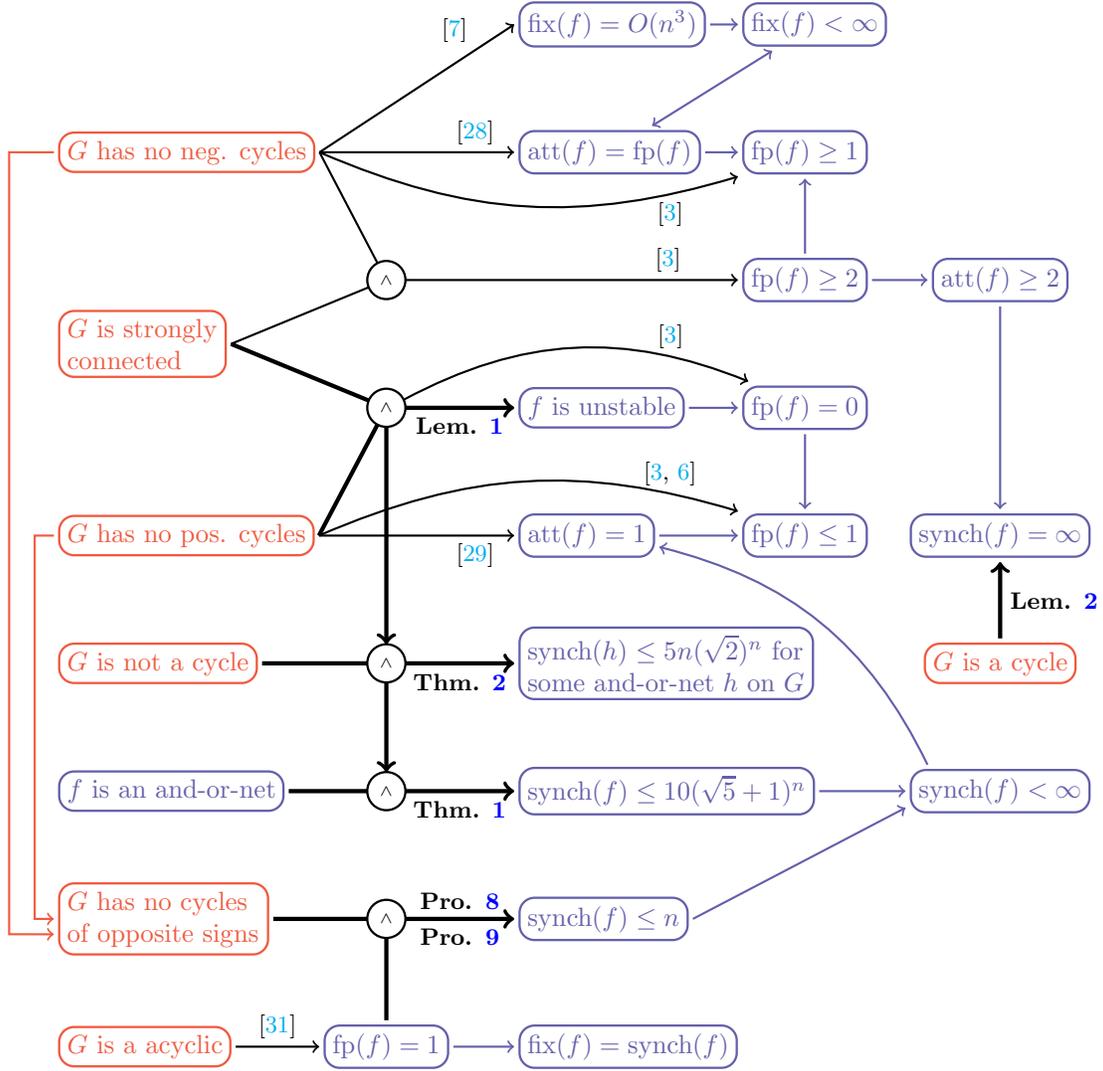

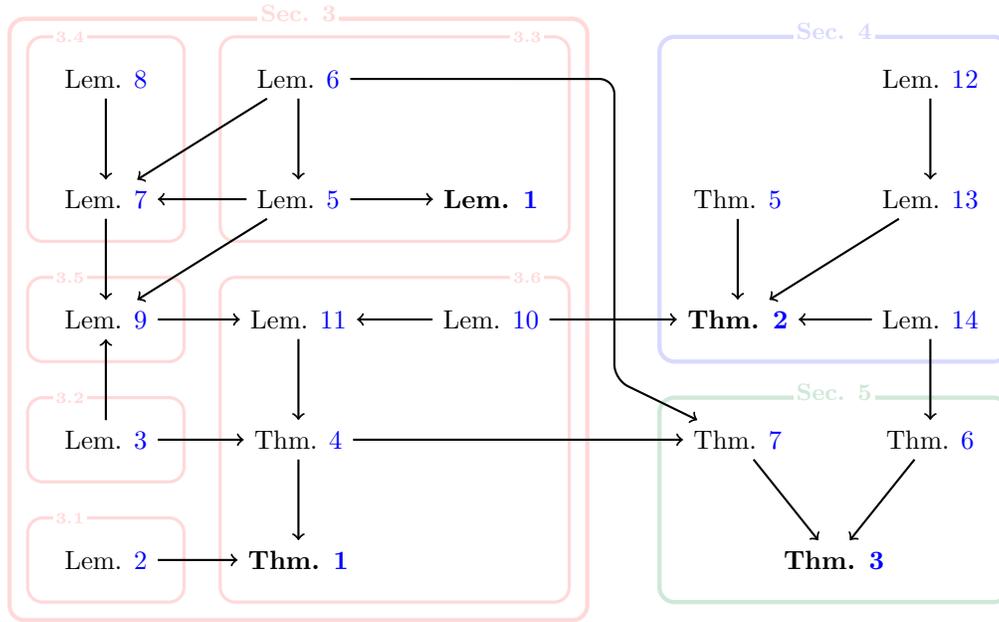
\begin{figure}
\[
\begin{tikzpicture}[scale=0.8]
\tikzstyle{every node}=[font=\small,fill=white,inner sep=1,outer sep=3,anchor=center]
\def\h{0.8}
\def\s{0.9}
\node at (0*\h,8) 	(L8) {Lem. \ref{lem:coherent_local_configuration}};
\node at (0*\h,6) 	(L7) {Lem. \ref{lem:fixing_a_vertex}};
\node at (0*\h,4) 	(L9) {Lem. \ref{lem:2-synchronization}};
\node at (0*\h,2) 	(L3) {Lem. \ref{lem:invariant_homogeneity}};
\node at (0*\h,0) 	(L2) {Lem. \ref{lem:initial_cycle}};
\node at (4*\h,8) 	(L6) {Lem. \ref{lem:diffusion}};
\node at (4*\h,6) 	(L5) {Lem. \ref{lem:no_fixed_vertex_2}};
\node at (4*\h,4) 	(L11) {Lem. \ref{lem:synchro_and_net}};
\node at (4*\h,2) 	(T4) {Thm. \ref{thm:main2}};
\node at (4*\h,0) 	(T1) {{\bf Thm. \ref{thm:main}}};
\node at (8*\h,6) 	(L1) {{\bf Lem. \ref{lem:flipping_lemma}}};
\node at (8*\h,4) 	(L10) {Lem. \ref{lem:2-synch}};
\node at (12*\h+\s,6) 	(T5) {Thm. \ref{thm:chromatic}};
\node at (12*\h+\s,4) 	(T2) {{\bf Thm. \ref{thm:fast}}};
\node at (12*\h+\s,2) 	(T7) {Thm. \ref{thm:every_synch}};
\node at (16*\h+\s,8) 	(L12) {Lem. \ref{lem:fast-2}};
\node at (16*\h+\s,6) 	(L13) {Lem. \ref{lem:fast-3}};
\node at (16*\h+\s,4) 	(L14) {Lem. \ref{lem:acyclic_image}};
\node at (16*\h+\s,2) 	(T6) {Thm. \ref{thm:fixed_point}};
\node at (14*\h+\s,0) 	(T3) {{\bf Thm. \ref{thm:complexity}}};
\def\sv{0.7}
\def\sh{1.3}
\def\svb{1.0}
\def\shb{1.6}
\draw[red!15,ultra thick,rounded corners=0.2cm] ({0*\h-\shb},{0-\svb}) rectangle 	({8*\h+\shb},{8+\svb});
\draw[red!15,very thick,rounded corners=0.2cm] ({0*\h-\sh},{0-\sv}) rectangle 	({0*\h+\sh},{0+\sv});
\draw[red!15,very thick,rounded corners=0.2cm] ({0*\h-\sh},{2-\sv}) rectangle 	({0*\h+\sh},{2+\sv});
\draw[red!15,very thick,rounded corners=0.2cm] ({0*\h-\sh},{4-\sv}) rectangle 	({0*\h+\sh},{4+\sv});
\draw[red!15,very thick,rounded corners=0.2cm] ({0*\h-\sh},{6-\sv}) rectangle 	({0*\h+\sh},{8+\sv});
\draw[red!15,very thick,rounded corners=0.2cm] ({4*\h-\sh},{0-\sv}) rectangle 	({8*\h+\sh},{4+\sv});
\draw[red!15,very thick,rounded corners=0.2cm] ({4*\h-\sh},{6-\sv}) rectangle 	({8*\h+\sh},{8+\sv});
\draw[blue!15,ultra thick,rounded corners=0.2cm] ({12*\h-\sh+\s},{4-\sv}) rectangle ({16*\h+\sh+\s},{8+\sv});
\draw[Green!18,ultra thick,rounded corners=0.2cm] ({12*\h-\sh+\s},{0-\sv}) rectangle ({16*\h+\sh+\s},{2+\sv});
\draw[thick,rounded corners=0.2cm,->] (L6) -- (10*\h+\s*0.5,{8-0.0}) -- (10*\h+\s*0.5,{2+1.0}) -- (T7);
\node[text=red!15,font=\footnotesize] at (4*\h,8+\svb+0.1) 	(S3) {{\bf Sec. 3}};
\node[text=blue!15,font=\footnotesize] at (14*\h+\s,8+\sv+0.1) 	(S4) {{\bf Sec. 4}};
\node[text=Green!18,font=\footnotesize] at (14*\h+\s,2+\sv+0.1) 	(S5) {{\bf Sec. 5}};
\node[text=red!15,font=\tiny] at (0*\h-0.6,8+\sv) 	(S34) {{\bf 3.4}};
\node[text=red!15,font=\tiny] at (0*\h-0.6,4+\sv) 	(S35) {{\bf 3.5}};
\node[text=red!15,font=\tiny] at (0*\h-0.6,2+\sv) 	(S32) {{\bf 3.2}};
\node[text=red!15,font=\tiny] at (0*\h-0.6,0+\sv) 	(S31) {{\bf 3.1}};
\node[text=red!15,font=\tiny] at (8*\h+0.6,8+\sv) 	(S33) {{\bf 3.3}};
\node[text=red!15,font=\tiny] at (8*\h+0.6,4+\sv) 	(S36) {{\bf 3.6}};
\draw[thick,->]
(L6) edge (L5)
(L6) edge (L7)
(L8) edge (L7)
(L5) edge (L7)
(L5) edge (L1)
(L7) edge (L9)
(L3) edge (L9)
(L5) edge (L9)
(L9) edge (L11)
(L10) edge (L11)
(L11) edge (T4)
(L3) edge (T4)
(T4) edge (T1)
(L2) edge (T1)
(T4) edge (T7)
(T7) edge (T3)
(T6) edge (T3)
(L10) edge (T2) 
(L14) edge (T6)
(L14) edge (T2)
(T5) edge (T2)
(L13) edge (T2)
(L12) edge (L13)
;
\end{tikzpicture}
\]
\caption{\label{fig:dependency}
Dependency graph of the main and auxiliary results.
}
\end{figure}

\medskip
The paper is organized as follows. In Section \ref{sec:preliminaries}, basic definitions and results are given; Proposition~\ref{pro:02} is proved there. The proofs of Theorems~\ref{thm:main}, \ref{thm:fast} and \ref{thm:complexity} are given in Sections \ref{sec:main}, \ref{sec:fast} and \ref{sec:complexity}, respectively. Some conclusions and perspectives are given in Section~\ref{sec:conclu}. In Appendix \ref{app:one_fixed_point}, the two results mentioned before the statement of Theorem~\ref{thm:main} are proved (see  Propositions \ref{pro:positive_G} and \ref{pro:negative_G} and in Appendix \ref{app:one_fixed_point}). A summary of the main results and related results from the literature is given in Figure \ref{fig:summary}, and a dependency graph of the main and auxiliary results is given in Figure \ref{fig:dependency}

%%%%%%%%%%%%%%%%%%%%%%%%%%%%%%%%%%%%%%%%%%%%%%%%%%%%%%%%%%%%%%%%%%%%%%%%%%%%%%%%%%%%%%%%%%%%%%%%%%%%%%%
\section{Preliminaries}\label{sec:preliminaries}
%%%%%%%%%%%%%%%%%%%%%%%%%%%%%%%%%%%%%%%%%%%%%%%%%%%%%%%%%%%%%%%%%%%%%%%%%%%%%%%%%%%%%%%%%%%%%%%%%%%%%%%

%%%%%%%%%%%%%%%%%%%%%%%%%%%%%%%%%%%%%%%%%%%%%%%%%%%%%%%%%%%%%%%%%%%%%%%%%%%%%%%%%%%%%%%%%%%%%%%%%%%%%%%
\subsection{Digraphs and signed digraphs}
%%%%%%%%%%%%%%%%%%%%%%%%%%%%%%%%%%%%%%%%%%%%%%%%%%%%%%%%%%%%%%%%%%%%%%%%%%%%%%%%%%%%%%%%%%%%%%%%%%%%%%%

A \EM{digraph} is a couple $G=(V,E)$ where $V$ is a set of vertices and $E\subseteq V^2$ is a set of arcs. Given $I\subseteq V$, the subgraph of $G$ induced by $I$ is denoted $G[I]$, and $G\setminus I$ means $G[V\setminus I]$. A \EM{strongly connected component} (\EM{strong component} for short) of a digraph $G$ is an induced subgraph which is \EM{strongly connected} (\EM{strong} for short) and maximal for this property. A strong component $G[I]$ is \EM{initial} if $G$ has no arc from $V\setminus I$ to $I$, and \EM{terminal} if $G$ has no arc from $I$ to $V\setminus I$. A \EM{source} is a vertex of in-degree zero. A digraph is \EM{trivial} if it has a unique vertex and no arc. 

\medskip
A \EM{signed digraph} $G$ is a couple $(V,E)$ where $E\subseteq V^2\times\{-1,1\}$. If $(j,i,s)\in E$ then $G$ has an arc from $j$ to $i$ of sign $s$; we also say that $j$ is an in-neighbor of $i$ of sign $s$ and that $i$ is an out-neighbor of $j$ of sign $s$. We say that $G$ is \EM{simple} if it has not both a positive arc and a negative arc from one vertex to another, and \EM{full-positive} if all its arcs are positive. A subgraph of $G$ is a signed digraph $(V',E')$ with $V'\subseteq V$ and $E'\subseteq E$. Cycles and paths of $G$ are regarded as simple subgraphs (and thus have no repeated vertices). The \EM{sign of a cycle or a path} of $G$ is the product of the signs of its arcs. The underlying (unsigned) digraph of $G$ has vertex set $V$ and an arc from $j$ to $i$ if $G$ has a positive or a negative arc from $j$ to $i$. Every graph concept made on $G$ that does not involved signs are tacitly made on its underlying digraph. For instance, $G$ is strongly connected if its underlying digraph is.

%%%%%%%%%%%%%%%%%%%%%%%%%%%%%%%%%%%%%%%%%%%%%%%%%%%%%%%%%%%%%%%%%%%%%%%%%%%%%%%%%%%%%%%%%%%%%%%%%%%%%%%
\subsection{Configurations and words}
%%%%%%%%%%%%%%%%%%%%%%%%%%%%%%%%%%%%%%%%%%%%%%%%%%%%%%%%%%%%%%%%%%%%%%%%%%%%%%%%%%%%%%%%%%%%%%%%%%%%%%%

A \EM{configuration} on a finite set $V$ is a map $x$ from $V$ to $\{0,1\}$, and for $i\in V$ we write $x_i$ instead of $x(i)$. The set of configurations on $V$ is denoted $\B^V$. Given $I\subseteq V$, we denote by $x_I$ the restriction of $x$ to $I$. We denote by $e_I$ the configuration on $V$ defined by $(e_I)_i=1$ for all $i\in I$ and $(e_I)_i=0$ for all $i\in V\setminus I$. If $i\in V$, we write $e_i$ instead of $e_{\{i\}}$. If $x,y$ are configurations on $V$ then $x+y$ is the configuration $z$ on $V$ such that $z_i=x_i+y_i$ for all $i\in V$, where the sum is modulo two. Hence $x+e_i$ is the configuration obtained from $x$ by flipping component $i$. We denote by $\ZERO$ and $\ONE$ the all-zero and all-one configurations. 

\medskip
A word $w$ over $V$ is a finite sequence of elements in $V$; its length is denoted $|w|$ and the set of letters that appear in $w$ is denoted $\{w\}$. The concatenation of two words $u$ and $v$ is denoted $u,v$ or $uv$. The empty word, of length $0$, is denoted $\epsilon$.  

%%%%%%%%%%%%%%%%%%%%%%%%%%%%%%%%%%%%%%%%%%%%%%%%%%%%%%%%%%%%%%%%%%%%%%%%%%%%%%%%%%%%%%%%%%%%%%%%%%%%%%%
\subsection{Boolean networks} 
%%%%%%%%%%%%%%%%%%%%%%%%%%%%%%%%%%%%%%%%%%%%%%%%%%%%%%%%%%%%%%%%%%%%%%%%%%%%%%%%%%%%%%%%%%%%%%%%%%%%%%%

A \EM{Boolean network} (BN) with component set $V$ is a function $f:\B^V\to\B^V$. For $i\in V$, we denote by $f_i$ the map from $\B^V$ to $\B$ defined by $f_i(x)=f(x)_i$ for all configurations $x$ on $V$. For $a\in\B$, we say that $f_i$ is the \EM{$a$-constant function} if $f_i(x)=a$ for all $x\in\B^V$. 

\medskip
The \EM{signed interaction digraph} of $f$ is the signed digraph $G$ with vertex set $V$ such that, for all $i,j\in V$, there is a positive (negative) arc from $j$ to $i$ is there exists a configuration $x$ on $V$ with $x_j=0$ such that $f_i(x+e_j)-f_i(x)$ is positive (negative). Given a signed digraph $G$, a BN \EM{on} $G$ is a BN whose signed interaction digraph is~$G$. 

\medskip
Given a vertex $i$ in $G$, an \EM{$i$-unstable} configuration in $G$ is a configuration $x\in\B^V$ such that $x_j\neq x_i$ for every positive in-neighbor $j$ of $i$ in $G$ and $x_j= x_i$ for every negative in-neighbor $j$ of $i$ in $G$. Almost all results that gives relationships between $G$ and the dynamical properties of the BNs on $G$ use, at some point, the following property; the proof is easy and omitted.  

\begin{proposition}\label{pro:flipping}
If $i$ is a vertex of $G$ of in-degree at least one, and $x$ is an $i$-unstable configuration in $G$, then $f_i(x)\neq x_i$ for every BN $f$ on $G$. 
\end{proposition}

Let $f$ be a BN with component set $V$ and $G$ its signed interaction digraph. For $i\in V$, we say that $f_i$ is a \EM{conjunction} if, for all configurations $x$ on $V$, we have $f_i(x)=1$ if and only if $x_j=1$ for all positive in-neighbors $j$ of $i$ in $G$ and $x_j=0$ for all the negative in-neighbors $j$ of $i$ in $G$. Similarly, we say that $f_i$ is a \EM{disjunction} if, for all configurations $x$ on $V$, we have $f_i(x)=0$ if and only if $x_j=0$ for all positive in-neighbors $j$ of $i$ in $G$ and $x_j=1$ for all the negative in-neighbors $j$ of $i$ in $G$. We say that $f$ is an \EM{and-or-net} if $f_i$ is a conjunction or a disjunction for all $i\in V$. We say that $f$ is an \EM{and-net} (\EM{or-net}) if $f_i$ is a conjunction (disjunction) for all $i\in V$. If $G$ is simple and has maximum in-degree at most $2$, then every BN $f$ on $G$ is an and-or-net.   

\medskip
For all $i\in V$, we denote by $f^i$ the BN with component set $V$ defined as follows: for all configurations $x$ on $V$, we have $f^i(x)_i=f(x)_i$ and $f^i(x)_j=x_j$ for all $j\in V\setminus\{i\}$. Given a word $w=i_1,\dots,i_\ell$ over $V$, we set 
$f^w=f^{i_\ell}\circ f^{i_{\ell-1}}\circ\dots\circ f^{i_1}$. For convenience, $f^\epsilon$, where $\epsilon$ is the empty word, is the identity on $\B^V$. Also, if $w$ is a word and $i\in\{w\}\setminus V$ then $f^i$ is the identity on $V$. In this way, $f^w$ is well defined for any word $w$. We say that $w$ is a \EM{synchronizing word} for $f$ if $f^w$ is a constant function, and $f$ is \EM{synchronizing} if it has at least one synchronizing word.

\medskip
The functions $f^i$ define a deterministic finite automaton, called \EM{asynchronous automaton} of~$f$, in a natural way: the alphabet is $V$, the set of states is $\B^V$ and the transition function $\delta:\B^V\times [n]\to\B^V$ is defined by $\delta(x,i)=f^i(x)$ for all $(x,i)\in\B^V\times [n]$. So $f$ is synchronizing if and only if its asynchronous automaton is synchronizing in the usual sense. The \EM{state diagram} of the asynchronous automaton of $f$ is then the labelled digraph $\Gamma(f)$ with vertex set $\B^V$ and with an arc from $x$ to $f^i(x)$ labelled by $i$ for all $(x,i)\in\B^V\times [n]$. It is a very classical model for the dynamics of gene networks \cite{T73,TA90}. In this context, $\Gamma(f)$ is called the asynchronous transition graph of $f$.

\begin{example}\label{ex1}
Let $G$ be the following signed digraph, where green arcs are positive and red arcs are negative (this convention is used throughout the paper):
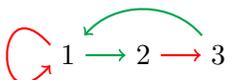
\[
\begin{tikzpicture}
\useasboundingbox (-0.9,-0.4) rectangle (2.2,0.7);
%\draw (-0.9,-0.4) rectangle (2.2,0.7);
\node (1) at (0,0){$1$};
\node (2) at (1,0){$2$};
\node (3) at (2,0){$3$};
\draw[thick,->,red] (1.150) .. controls (-1,+1) and (-1,-1) .. (1.210);
\path[thick,->]
(1) edge[Green] (2)
(2) edge[red] (3)
(3) edge[Green,bend right=50] (1)
;
\end{tikzpicture}
\]
It is simple and has exactly two cycles, both negative. Let $f$ be the and-net on $G$, that is: $f_1(x)=\overline{x_1}\land x_3$ and $f_2(x)=x_1$ and $f_3(x)=\overline{x_2}$ for all configurations $x$ on $[3]$. We have the following tables:
\[
\begin{array}{c|c}
x   & f(x)\\\hline
000 & 001\\ 
001 & 101\\
010 & 000\\
011 & 100\\
100 & 011\\
101 & 011\\
110 & 010\\
111 & 010\\
\end{array}
\quad
\begin{array}{c|c}
x   & f^1(x)\\\hline
000 & 000\\ 
001 & 101\\
010 & 010\\
011 & 111\\
100 & 000\\
101 & 001\\
110 & 010\\
111 & 011\\
\end{array}
\quad
\begin{array}{c|c}
x   & f^2(x)\\\hline
000 & 000\\ 
001 & 001\\
010 & 000\\
011 & 001\\
100 & 110\\
101 & 111\\
110 & 110\\
111 & 111\\
\end{array}
\quad
\begin{array}{c|c}
x   & f^3(x)\\\hline
000 & 001\\ 
001 & 001\\
010 & 010\\
011 & 010\\
100 & 101\\
101 & 101\\
110 & 110\\
111 & 110\\
\end{array}
\]
The state diagram of the asynchronous automaton of $f$ is the following:
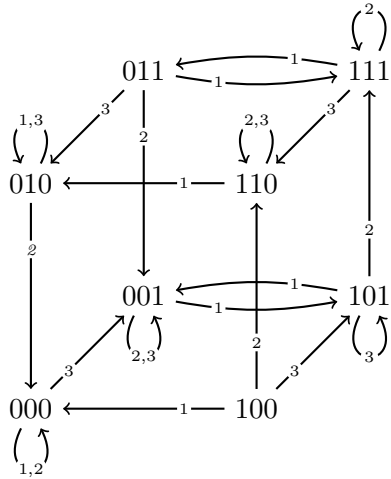
\[
\begin{array}{c}
%components
\begin{tikzpicture}
\node (000) at (0,0){{\em 000}};
\node (001) at (1.5,1.5){{\em 001}};
\node (010) at (0,3){{\em 010}};
\node (011) at (1.5,4.5){{\em 011}};
\node (100) at (3,0){{\em 100}};
\node (101) at (4.5,1.5){{\em 101}};
\node (110) at (3,3){{\em 110}};
\node (111) at (4.5,4.5){{\em 111}};
%loops
\draw[->,thick](000) .. node[outer sep=0,inner sep=0.5,fill=white,circle]{\tiny {\em 1,2}}  
controls ++(-0.5,-1) and ++(0.5,-1) .. (000);
\draw[->,thick](001) .. node[outer sep=0,inner sep=0.5,fill=white,circle]{\tiny {\em 2,3}} 
controls ++(-0.5,-1) and ++(0.5,-1) .. (001);
\draw[->,thick](010) .. node[outer sep=0,inner sep=0.5,fill=white,circle]{\tiny {\em 1,3}}
controls ++(0.5,1) and ++(-0.5,1) .. (010);
\draw[->,thick](101) .. node[outer sep=0,inner sep=0.5,fill=white,circle]{\tiny {\em 3}} 
controls ++(-0.5,-1) and ++(0.5,-1) .. (101);
\draw[->,thick](110) .. node[outer sep=0,inner sep=0.5,fill=white,circle]{\tiny {\em 2,3}} 
controls ++(0.5,1) and ++(-0.5,1) .. (110);
\draw[->,thick](111) .. node[outer sep=0,inner sep=0.5,fill=white,circle]{\tiny {\em 2}}
 controls ++(0.5,1) and ++(-0.5,1) .. (111);
%transitions
\path[thick,->,draw,black]
(000) edge node[near start,outer sep=0,inner sep=0.5,fill=white,circle]{\tiny {\em 3}} (001)
(001) edge[bend right=10] node[near start,outer sep=0,inner sep=0,fill=white,circle]{\tiny {\em 1}} (101)
(010) edge node[near start,outer sep=0,inner sep=0.5,fill=white,circle]{\tiny 2} (000)
(011) edge[bend right=10] node[near start,outer sep=0,inner sep=0,fill=white,circle]{\tiny {\em 1}} (111)
(011) edge node[near start,outer sep=0,inner sep=0.5,fill=white,circle]{\tiny {\em 2}} (001)
(011) edge node[near start,outer sep=0,inner sep=0.5,fill=white,circle]{\tiny {\em 3}} (010)
(100) edge node[near start,outer sep=0,inner sep=0,fill=white,circle]{\tiny {\em 1}} (000)
(100) edge node[near start,outer sep=0,inner sep=0.5,fill=white,circle]{\tiny {\em 2}} (110)
(100) edge node[near start,outer sep=0,inner sep=0.5,fill=white,circle]{\tiny {\em 3}} (101)
(101) edge[bend right=10] node[near start,outer sep=0,inner sep=0,fill=white,circle]{\tiny {\em 1}} (001)
(101) edge node[near start,outer sep=0,inner sep=0.5,fill=white,circle]{\tiny {\em 2}} (111)
(110) edge node[near start,outer sep=0,inner sep=0,fill=white,circle]{\tiny {\em 1}} (010)
(111) edge[bend right=10] node[near start,outer sep=0,inner sep=0,fill=white,circle]{\tiny {\em 1}} (011)
(111) edge node[near start,outer sep=0,inner sep=0.5,fill=white,circle]{\tiny {\em 3}} (110)
;
\end{tikzpicture}
\end{array}
\]
$f$ is synchronizing, since $w=231123$ is a synchronizing word for $f$. Indeed, $f^w(x)=001$ for all configurations $x$ on $[3]$, as shown below:
\[
\begin{array}{c}
000 \\ 
001 \\
010 \\
011 \\
100 \\
101 \\
110 \\
111 \\
\end{array}
\xrightarrow{f^2}
\begin{array}{c}
000 \\ 
010 \\
110 \\
111 \\
\end{array}
\xrightarrow{f^3}
\begin{array}{c}
001 \\ 
110 \\
\end{array}
\xrightarrow{f^1}
\begin{array}{c}
101 \\ 
010 \\
\end{array}
\xrightarrow{f^1}
\begin{array}{c}
001 \\ 
010 \\
\end{array}
\xrightarrow{f^2}
\begin{array}{c}
001 \\ 
000 \\
\end{array}
\xrightarrow{f^3}
\begin{array}{c}
001 \\ 
\end{array}
\] 
Let $h$ be the or-net on $G$, that is: $h_1(x)=\overline{x_1}\lor x_3$ and $h_2(x)=x_1$ and $h_3(x)=\overline{x_2}$ for all configurations $x$ on $[3]$. $h$ is synchronizing, since one can check that $w=231123$ is also a synchronizing word for $h$: we have $h^w(x)=110$ for all configurations $x$ on $[3]$. Since there are only two BNs on $G$ (namely $f$ and $h$),  $w=231123$ is a synchronizing word for every BN on $G$. 
\end{example}

\medskip
We now prove Proposition~\ref{pro:02} given in the introduction, that we restate.

\setcounter{proposition}{0}
\begin{proposition}
Suppose that $G$ is a simple signed digraph without vertices of in-degree $0$ or $2$. Then at least one BN on $G$ is not synchronizing.
\end{proposition} 
\setcounter{proposition}{2}

\begin{proof}
Let $V$ be the vertex set of $G$. Given an arc of $G$ from $j$ to $i$, we set $a_{ji}=0$ if this arc is positive and $a_{ji}=1$ otherwise. For all $i\in V$, let $d_i$ be the in-degree of $i$ in $G$, and suppose that $d_i\neq 0,2$. Let $f$ be the BN with component set $V$ defined as follows. Given a configuration $x$ on $V$, let $h_i(x)$ be the number of in-neighbors $j$ of $i$ with $x_j+a_{ji}=1$. For each vertex $i$, we fix an in-neighbor $\phi(i)$ of $i$ and we define $f_i$ as follows: for all configurations $x$ on $V$, 
\[
f_i(x)=
\left\{
\begin{array}{ll}
1&\textrm{ if $h_i(x)>d_i/2$,}\\[1mm]
x_{\phi(i)}+a_{\phi(i)i}&\textrm{ if $h_i(x)=d_i/2$},\\[1mm]
0&\textrm{ if $h_i(x)<d_i/2$.}
\end{array}
\right.
\]

\medskip
Let $H$ be the signed interaction digraph of $f$, and let us prove that $H=G$. One easily checks that $H$ is a subgraph of $G$. To prove the converse, suppose that $G$ has an arc from $j$ to $i$. If $d_i$ is odd there is $x$ such that $x_j+a_{ji}=0$ and $h_i(x)=(d_i-1)/2$. Thus $f_i(x)<f_i(x+e_j)$ so $H$ has an arc from $j$ to $i$ which is positive if $a_{ji}=0$ and negative otherwise. If $d_i$ is even then, since $d_i\geq 4$, there is $x$ such that $x_j+a_{ji}=x_{\phi(i)}+a_{\phi(i)i}=0$ and $h_i(x)=d_i/2$. We deduce that $f_i(x)<f_i(x+e_j)$ so $H$ has an arc from $j$ to $i$ which is positive if $a_{ji}=0$ and negative otherwise. Consequently, $H=G$, that is, $f$ is a BN on $G$. 

\medskip
We now prove that $f$ is not synchronizing. Let $x,y$ be opposite configurations, that is, $x_i\neq y_i$ for all $i\in V$. We have $h_i(y)=d_i-h_i(x)$, thus if $h_i(x)\neq d_i/2$ then $f_i(x)\neq f_i(y)$. Otherwise, $h_i(x)=h_i(y)=d_i/2$ thus $f_i(x)=x_{\phi(i)}+a_{\phi(i)i}\neq y_{\phi(i)}+a_{\phi(i)i}=f_i(y)$. So $f_i(x)\neq f_i(y)$ in every case. Consequently, $f^i(x)$ and $f^i(y)$ are opposite. We deduce that, for any word $w$, $f^w(x)$ and $f^w(y)$ are opposite, and thus $f$ is not synchronizing.
\end{proof}

%%%%%%%%%%%%%%%%%%%%%%%%%%%%%%%%%%%%%%%%%%%%%%%%%%%%%%%%%%%%%%%%%%%%%%%%%%%%%%%%%%%%%%%%%%%%%%%%%%%%%%%
\subsection{Switches}
%%%%%%%%%%%%%%%%%%%%%%%%%%%%%%%%%%%%%%%%%%%%%%%%%%%%%%%%%%%%%%%%%%%%%%%%%%%%%%%%%%%%%%%%%%%%%%%%%%%%%%%

Let $G=(V,E)$ be a signed digraph. Let $I\subseteq V$ and, for all $i\in V$, let $\sigma_I(i)=1$ if $i\in I$ and $\sigma_I(i)=-1$ if otherwise. The \EM{$I$-switch} of $G$ is the signed digraph $G^I=(V,E^I)$ with $E^I=\{(j,i,\sigma_I(j)\cdot s\cdot\sigma_I(i))\mid (j,i,s)\in E\}$; note that $G^I=G^{V\setminus I}$ and $(G^I)^I=G$. We say that $G$ is \EM{switch equivalent} to $H$ if $H=G^I$ for some $I\subseteq V$. Obviously, $G$ and $G^I$ have the same underlying digraph. Note also that $C$ is a cycle in $G$ if and only if $C^I$ is a cycle in $G^I$, and $C$ and $C^I$ have the same sign. Thus if $G$ has no positive (negative) cycles then every switch of $G$ has no positive (negative) cycles: this property is invariant by switch. The \EM{symmetric version} of $G$ is the signed digraph $G^s=(V,E^s)$ where $E^s=E\cup \{(i,j,s)\mid (j,i,s)\in E\}$. A basic result concerning switch is the following adaptation of Harary's theorem \cite{H53}.

\begin{proposition}\label{pro:harary}
A signed digraph $G$ is switch equivalent to a full-positive signed digraph if and only if $G^s$ has no negative cycles.
\end{proposition}

There is an analogue of the switch operation for BNs. Let $f$ be a BN with component set $V$ and $I\subseteq V$. The \EM{$I$-switch} of $f$ is the BN $h$ with component set $V$ defined by $h(x)=f(x+e_I)+e_I$ for all configurations $x$ on $V$; note that if $h$ is the $I$-switch of $f$ then $f$ is the $I$-switch of $h$. The analogy comes from first point of the following easy property. 

\begin{proposition}\label{pro:BN_switch}
If $h$ is the $I$-switch of $f$, then:
\begin{itemize}
\item the signed interaction digraph of $h$ is the $I$-switch of the signed interaction digraph of $f$;
\item if $f$ is an and-or-net, then $h$ is an and-or-net;
\item for any word $w$, $h^w$ is the $I$-switch of $f^w$.
\end{itemize}
\end{proposition}

\begin{example}\label{ex2}
Let $G$ and $H$ be the following signed digraphs:
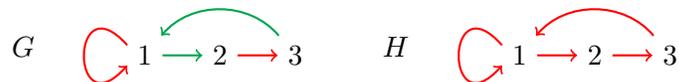
\[
G\quad
\begin{array}{c}
\begin{tikzpicture}
\useasboundingbox (-0.9,-0.4) rectangle (2.2,0.7);
%\draw (-0.9,-0.4) rectangle (2.2,0.7);
\node (1) at (0,0){$1$};
\node (2) at (1,0){$2$};
\node (3) at (2,0){$3$};
\draw[thick,->,red] (1.150) .. controls (-1,+1) and (-1,-1) .. (1.210);
\path[thick,->]
(1) edge[Green] (2)
(2) edge[red] (3)
(3) edge[Green,bend right=50] (1)
;
\end{tikzpicture}
\end{array}
\qquad
H\quad
\begin{array}{c}
\begin{tikzpicture}
\useasboundingbox (-0.9,-0.4) rectangle (2.2,0.7);
%\draw (-0.9,-0.4) rectangle (2.2,0.7);
\node (1) at (0,0){$1$};
\node (2) at (1,0){$2$};
\node (3) at (2,0){$3$};
\draw[thick,->,red] (1.150) .. controls (-1,+1) and (-1,-1) .. (1.210);
\path[thick,->]
(1) edge[red] (2)
(2) edge[red] (3)
(3) edge[red,bend right=50] (1)
;
\end{tikzpicture}
\end{array}
\]
Then $H$ is the $\{2,3\}$-switch of $G$. We deduce from the first point of Proposition~\ref{pro:BN_switch} that every BN on $H$ is the $I$-switch of a BN on $G$. As shown in Example~\ref{ex1}, $w=231123$ is a synchronizing word for every BN on $G$, and we deduce from the third point of Proposition~\ref{pro:BN_switch} that $w$ is also a synchronizing word for every BN on $H$. 
\end{example}

%%%%%%%%%%%%%%%%%%%%%%%%%%%%%%%%%%%%%%%%%%%%%%%%%%%%%%%%%%%%%%%%%%%%%%%%%%%%%%%%%%%%%%%%%%%%%%%%%%%%%%%
\section{Proof of Theorem~\ref{thm:main}}\label{sec:main}
%%%%%%%%%%%%%%%%%%%%%%%%%%%%%%%%%%%%%%%%%%%%%%%%%%%%%%%%%%%%%%%%%%%%%%%%%%%%%%%%%%%%%%%%%%%%%%%%%%%%%%%

In all this section, $G$ is a signed digraph with vertex set $V$, and $n=|V|$.

%%%%%%%%%%%%%%%%%%%%%%%%%%%%%%%%%%%%%%%%%%%%%%%%%%%%%%%%%%%%%%%%%%%%%%%%%%%%%%%%%%%%%%%%%%%%%%%%%%%%%%%
\subsection{Initial cycles}
%%%%%%%%%%%%%%%%%%%%%%%%%%%%%%%%%%%%%%%%%%%%%%%%%%%%%%%%%%%%%%%%%%%%%%%%%%%%%%%%%%%%%%%%%%%%%%%%%%%%%%%

An initial strong component of $G$ which is isomorphic to a cycle is an \EM{initial cycle} of $G$. The first case in Theorem~\ref{thm:main} is a consequence of the following easy property.

\begin{lemma}\label{lem:initial_cycle}
If $G$ has an initial cycle, then no BN on $G$ is synchronizing.
\end{lemma} 

\begin{proof}
Suppose that $G$ has an initial cycle with vertex set $I$, and let $f$ be a BN on $G$. Each vertex $i\in I$ has a unique in-neighbor  in $G$, say $\phi(i)$, and $\phi(i)\in I$. Then one can easily check that $f_i(x)=x_{\phi(i)}+a_i$ for all configurations $x$ on $V$, where $a_i=0$ if the arc from $\phi(i)$ to $i$ is positive, and $a_i=1$ otherwise. Let $x$ and $y$ be $I$-opposite configurations, that is, $x_i\neq y_i$ for all $i\in I$. Then, for all $i\in I$, we have $f_i(x)=x_{\phi(i)}+a_i\neq y_{\phi(i)}+a_i=f_i(y)$, thus $f^i(x)$ and $f^i(y)$ are $I$-opposite configurations. We deduce that $f^w(x)$ and $f^w(y)$ are $I$-opposite configurations for any word $w$, and thus $f$ is not synchronizing. 
\end{proof}

%%%%%%%%%%%%%%%%%%%%%%%%%%%%%%%%%%%%%%%%%%%%%%%%%%%%%%%%%%%%%%%%%%%%%%%%%%%%%%%%%%%%%%%%%%%%%%%%%%%%%%%
\subsection{Homogeneity}
%%%%%%%%%%%%%%%%%%%%%%%%%%%%%%%%%%%%%%%%%%%%%%%%%%%%%%%%%%%%%%%%%%%%%%%%%%%%%%%%%%%%%%%%%%%%%%%%%%%%%%%

Let us say that $G$ is \EM{and-or-synchronizing} if every and-or-net on $G$ is synchronizing. It remains to prove the second case in Theorem~\ref{thm:main}, and the difficult part is to show that if $G$ is strong, has no positive cycles, and is not a cycle, then $G$ is and-or-synchronizing. We proceed by induction on the number of vertices and, for inductive purpose, it is convenient to prove something stronger, relaxing the strong connectivity which is difficult to handle during inductive proofs.  

\medskip
If $G$ is no longer strongly connected, then Lemma~\ref{lem:initial_cycle} shows that the condition ``is not a cycle'' must be replaced by ``has no initial cycles''. So suppose that $G$ has no positive cycles and no initial cycles. Then $G$ is not necessarily and-or-synchronizing, as shown below, so additional assumptions are needed. 

\begin{example}
Let $G$ be the following signed digraph:
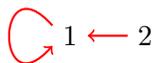
\[
\begin{tikzpicture}
\node (1) at (0,0){$1$};
\node (2) at (1,0){$2$};
\draw[thick,->,red] (1.150) .. controls (-1,+1) and (-1,-1) .. (1.210);
\path[thick,->]
(2) edge[red] (1)
;
\end{tikzpicture}
\]
Let $f$ be the and-net on $G$: $f_1(x)=\overline{x_1}\land \overline{x_2}$ and $f_2(x)=1$ for all configurations $x$ on~$[2]$. For $X=\{01,11\}$ we have $f^1(X)=X$ and $f^2(X)=X$. Hence, for any word $w$ we have $f^w(X)=X$ thus $f$ is not synchronizing.  
\end{example}

In the previous example, the source can be fixed to the state $1$, and the network then behaves as an ``isolated'' negative cycle of length one (since when $x_2=1$ we have $f_1(x)=\overline{x_1}$) which is not synchronizing (Lemma~\ref{lem:initial_cycle}). 

\medskip
To prevent this type of phenomenon, we can assume that $G$ has, in addition, no sources. But this is not enough, as shown by the next example. Actually, it shows that if $G$ has no positive cycles, no sources and no initial cycles, then $G$ is not necessarily and-or-synchronizing, even if each strong component of $G$ is and-or-synchronizing. Thus we have to impose some additional conditions about the connections between the strong components of $G$. 

\begin{example}
Let $H$ be the following signed digraph, which is strong, has no positive cycles and which is not a cycle:
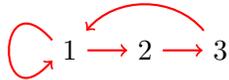
\[
\begin{tikzpicture}
\node (1) at (0,0){$1$};
\node (2) at (1,0){$2$};
\node (3) at (2,0){$3$};
\draw[thick,->,red] (1.150) .. controls (-1,+1) and (-1,-1) .. (1.210);
\draw[thick,->,white] (3.30) .. controls (3,+1) and (3,-1) .. (3.-30);
\path[thick,->]
(1) edge[red] (2)
(2) edge[red] (3)
(3) edge[red,bend right=50] (1)
;
\end{tikzpicture}
\]
By the second case of Theorem~\ref{thm:main}, $H$ is and-or-synchronizing (see Example~\ref{ex2}). Let $G$ be the following signed digraph, which has no positive cycles, no sources, and no initial cycles:
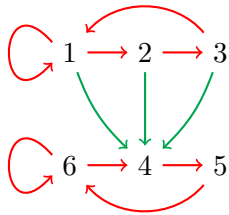
\[
\begin{tikzpicture}
\node (1) at (0,0){$6$};
\node (2) at (1,0){$4$};
\node (3) at (2,0){$5$};
\draw[thick,->,red] (1.150) .. controls (-1,+1) and (-1,-1) .. (1.210);
\draw[thick,->,white] (3.30) .. controls (3,+1) and (3,-1) .. (3.-30);
\path[thick,->]
(1) edge[red] (2)
(2) edge[red] (3)
(3) edge[red,bend left=50] (1)
;

\node (1') at (0,1.5){$1$};
\node (2') at (1,1.5){$2$};
\node (3') at (2,1.5){$3$};
\draw[thick,->,red] (1'.150) .. controls (-1,+2.5) and (-1,0.5) .. (1'.210);
\path[thick,->]
(1') edge[red] (2')
(2') edge[red] (3')
(3') edge[red,bend right=50] (1')
;

\path[thick,->]
(1') edge[Green, bend right=13] (2)
(2') edge[Green] (2)
(3') edge[Green,bend left=13] (2)
;
\end{tikzpicture}
\]
The two strong components of $G$ are and-or-synchronizing, since there are isomorphic to $H$. However, $G$ is not and-or-synchronizing. Indeed, let $f$ be the and-or-net such that all components of $f$ are conjunctions, excepted $f_6$ which is a disjunction. Let $X$ be the set of configurations $x$ on $[6]$ with $x_i=0$ for some $i\in [3]$, $x_4=0$ and $x_5=1$. One can easily check that $f^i(X)\subseteq  X$ for every $i\in [6]$. Furthermore, given $x\in X$, since $x_5=1$ we have $f_6(x)=\overline{x_5}\lor\overline{x_6}=\overline{x_6}$. Consequently, given $x,y\in X$ with $x_6\neq y_6$, we have $f^i(x),f^i(y)\in X$ and $f^i_6(x)\neq f^i_6(y)$ for every $i\in [6]$, and we deduce that $f^w_6(x)\neq f^w_6(y)$ for any word $w$. Hence $f$ is not synchronizing. 
\end{example}

To express the additional conditions between the strong components of $G$ we need some definitions. We say that a path $P$ of $G$ of is \EM{forward} if $P$ has at least one arc and the two vertices of the last arc belong to distinct strong components. Given a vertex $i\in V$, we say that $G$ is \EM{$i$-homogenous} if each strong component of $G$ contains a vertex $j$ such that all the forward paths from $j$ to $i$ have the same sign (obviously, if there is no forward paths from $j$ to $i$, then all the forward paths from $j$ to $i$ have the same sign). We say that $G$ is \EM{homogenous} if it is $i$-homogenous for every $i\in V$. 

\medskip
For instance, if $G$ is strong, then there are no forward paths, and $G$ is trivially homogenous, so homogeneity is indeed a relaxation of strong connectivity. Also, the $6$-vertex signed digraph of the previous example (which is not and-or-synchronizing) is not homogenous since, for $i=1,2,3$, it has a forward positive path and a forward negative path from $i$ to $4$. 

\medskip
It appears that homogeneity works and is well adapted for a proof by induction. So our aim is now to prove the following theorem, which generalizes the second case of Theorem~\ref{thm:main}. 

\begin{theorem}\label{thm:main2}
Suppose that $G$ is homogenous and has no positive cycles, no sources and no initial cycles. Then every and-or-net on $G$ has a synchronizing word of length at most $10(\sqrt{5}+1)^n$.
\end{theorem}

Thus Theorem~\ref{thm:main} is a consequence of Lemma \ref{lem:initial_cycle} (first case) and  Theorem~\ref{thm:main2} (second case).

\medskip
A useful observation is that the conditions of Theorem \ref{thm:main2} are invariant by switch. Indeed, we already mentioned that the absence of positive cycles is invariant by switch, and the following observation is straightforward to check. 

\begin{lemma}\label{lem:invariant_homogeneity}
If $G$ is homogenous, then every switch of $G$ is homogenous.
\end{lemma}

%%%%%%%%%%%%%%%%%%%%%%%%%%%%%%%%%%%%%%%%%%%%%%%%%%%%%%%%%%%%%%%%%%%%%%%%%%%%%%%%%%%%%%%%%%%%%%%%%%%%%%%
\subsection{Flipping a vertex}
%%%%%%%%%%%%%%%%%%%%%%%%%%%%%%%%%%%%%%%%%%%%%%%%%%%%%%%%%%%%%%%%%%%%%%%%%%%%%%%%%%%%%%%%%%%%%%%%%%%%%%%

The main tool is the following generalization of Lemma~\ref{lem:flipping_lemma}, obtained by replacing the strong connectivity of $G$ by the weaker assumption that $G$ is $i$-homogenous and without sources. 

\begin{lemma}\label{lem:no_fixed_vertex}
Suppose that $G$ is $i$-homogenous and has no positive cycles and no sources. Let $f$ be a BN on $G$. For every configuration $x$ on $V$, there is a word $w$ such that $f^w(x)_i\neq x_i$.
\end{lemma}

A shortest word $w$ with the above property is obviously of length $|w|\leq 2^{n-1}$ (the number of configurations $y$ with $y_i=x_i$), but we think that $|w|$ is sub-linear according to the number $2^n$ of configurations. We obtained such a sub-linear bound when $f$ is an and-or-net: $|w|< F(n+2)$ where $F(0),F(1),F(2)\dots$ is the Fibonacci sequence. 

\medskip
The word $w$ in Lemma~\ref{lem:no_fixed_vertex} is basically a concatenation of \EM{canalizing words}, defined as follows. Given a BN $f$ on $G$, a vertex $i\in V$ and $a\in\B$, a \EM{canalizing word from} $(i,a)$ is a word $w$ over $V\setminus i$ without repeated letters such that, for some configuration $b$ on $\{w\}$ we have:
\[
\forall x\in\B^V,\qquad x_i=a~\Rightarrow f^w(x)_{\{w\}}=b.
\]   
If $w$ is canalizing from $(i,a)$, the configuration $b$ with the above property is called the \EM{image} of $w$, and for each $j\in\{w\}$, we say that $j$ is \EM{canalized} to $b_j$ by $w$. By convention, the empty word is always a canalizing word.  

\begin{example}
Let $f$ be the and-net on $G$. Take $i\in V$ and let $I$ be the set of positive out-neighbors of $i$ distinct from $i$. For any configuration $x$ on $V$, if $x_i=0$ then $f_j(x)=0$ for all $j\in I$. We deduce that any permutation of $I$ is a canalizing word from $(i,0)$, which canalizes each member of $I$ to $0$. Similarly, if $J$ is the set of negative out-neighbors of $i$ distinct from $i$, then any permutation of $J$ is a canalizing word from $(i,1)$, which canalizes each member of $J$ to $0$. We deduce that if $i$ has at least one out-neighbor distinct from $i$, then there is a canalizing word from $(i,0)$ or $(i,1)$ which is of length at least one (this basic property is crucial for the bound of Lemma~\ref{lem:no_fixed_vertex_2} below). Suppose now that $G$ has a full-positive path $P$ of length $\ell\geq 1$ with vertices $i_0,i_1,\dots,i_\ell$ in order. Then $w=i_1i_2\dots i_\ell$ is a canalizing word from $(i_0,0)$ which canalizes $i_k$ to $0$ for $1\leq k\leq\ell$.
\end{example}

\medskip
The following bound on the length of the word $w$ in Lemma~\ref{lem:no_fixed_vertex}, which works for and-or-nets, is suited for a proof by induction; it implies the bound $|w|< F(n+2)$ announced above since $F(n+2-m)+m\leq F(n+2)$ for all $0\leq m\leq n$; this property will be used many times.

\begin{lemma}\label{lem:no_fixed_vertex_2}
Suppose that $G$ is $i$-homogenous and has no positive cycles and no sources. Let $f$ be an and-or-net on $G$. Let $x$ be a configuration on $V$ and let $v$ be a longest canalizing word from $(i,x_i)$. There is a word $w$ such that $f^w(x)_i\neq x_i$~and
\[
|w|< F(n-|v|+2)+|v|.
\]
\end{lemma}

We need the following lemma. It formalizes the intuitive fact that canalizations from $(i,a)$ operate through the paths of $G$ starting from $i$. 

\begin{lemma}\label{lem:diffusion}
Let $f$ be a BN on $G$, $i\in V$ and $a\in\B$. Let $w$ be a canalizing word from $(i,a)$ with image $b$, and suppose that no sources of $G$ are in $\{w\}$. For every $j\in\{w\}$, $G$ has a path from $i$ to $j$, whose internal vertices are in $\{w\}$, which is positive if $a=b_j$ and negative otherwise. 
\end{lemma}

\begin{proof}
We proceed by induction on $|w|$. If $|w|=0$ there is nothing to prove. Suppose that $|w|\geq 1$. Let $j$ be the last letter of $w$, and let $v$ be obtained from $w$ by removing the last letter (thus $v$ is the empty word if $|w|=1$). For convenience we set $b_i=a$ and $I=\{i\}\cup \{v\}$. We also write $\sigma(0)=-1$ and $\sigma(1)=1$. Hence, we have to prove that $G$ has a path from $i$ to $j$ of sign $\sigma(b_i)\cdot\sigma(b_j)$ whose internal vertices are in $\{w\}$. 

\medskip
Since $j$ is not a source, there is a configuration $x$ on $V$ such that $f_j(x)\neq b_j$. Let $z=x$ if $x_i=a$ and $z=x+e_i$ otherwise, so that $z_i=a$ in any case. Let $y=f^v(z)$. Since $z_i=a$ we have $f_j(y)=f^w(z)_j=b_j$ and $y_k=f^w(z)_k=b_k$ for $k\in I$. Thus $f_j(y)\neq f_j(x)$ and we deduce that $j$ has an in-neighbor $k$ in $G$ of sign $(y_k-x_k)\cdot (f_j(y)-f_j(x))=\sigma(y_k)\cdot \sigma(f_j(y))=\sigma(y_k)\cdot \sigma(b_j)$. Since $x$ and $y$ can only differ in components inside $I$, we have $k\in I$, and we deduce that the sign of the arc from $k$ to $j$ is $\sigma(b_k)\cdot \sigma(b_j)$. If $k=i$ we are done: the desired path is the arc from $i$ to~$j$. Otherwise, $k\in\{v\}$. Since $v$ is a canalizing word from $(i,a)$ with image $b_{\{v\}}$, by induction, $G$ has a path from $i$ to $k$ of sign $\sigma(b_i)\cdot \sigma(b_k)$ whose internal vertices are in $\{v\}$. By adding the arc from $k$ to $j$ to this path, we obtain a path of sign $\sigma(b_i)\cdot \sigma(b_k)\cdot \sigma(b_k)\cdot \sigma(b_j)=\sigma(b_i)\cdot \sigma(b_j)$  from $i$ to $j$ whose internal vertices are in $\{w\}$. This completes the induction step.   
\end{proof}

We are now ready to prove Lemmas~\ref{lem:no_fixed_vertex} and \ref{lem:no_fixed_vertex_2} together. The sketch is the following. Suppose that the initial state of $i$ is $0$, that is, $x_i=0$. So we want a word $w$ such that $f^w(x)_i=1$, that is, a word increasing the state of $i$. 

\medskip
We first suppose that $i$ is of out-degree zero. Since $G$ has no sources, it has is a non-trivial initial strong component $S$. Since $G$ is $i$-homogenous and $i$ is of out-degree zero, there is a vertex $j$ in $S$ such that all the paths from $j$ to $i$ have the same sign, say positive. Then $j$ can be seen as an activator of $i$, so in order to increase the state of $i$, it is preferable to put $j$ in state~$1$. This is possible by induction on $S$, since $S$ is $j$-homogenous. So we consider a word $u$ that puts $j$ in state~$1$. Then we consider a longest canalizing word $v$ from $(j,1)$. If $i$ is canalized by $v$, then it necessarily gets the state $1$, since all the paths from $j$ to $i$ are positive, and we are done: starting from $x$, the word $uv$ increases the state of $i$. Otherwise, we ``remove'' $j$ and the vertices canalized by $v$. We obtain a ``subnetwork'' whose interaction graph $H$ satisfies all the hypotheses. Using induction, we get, for the subnetwork, a word $w$ that increases the state of $i$ and we are done: starting from $x$, the word $uvw$ increases the state of $i$. We then prove that $|uvw|\leq F(n+2)$ when $f$ is an and-or-net (which is the bound of the statement since $i$ is of out-degree zero). The argument is roughly the following. The induction used to obtain $u$ is done on $S$, which has at most $n-1$ vertices, and thus $|u|< F(n+1)$. The induction used to obtained $w$ is made on $H$, which has $n-1-|v|$ vertices, and we easily obtain $|vw|<F(n+1)$. So if $|u|\leq 1$ then $|uvw|<F(n+2)$ as desired. If $|u|\geq 2$ then it means, and this is the key point, that there is a {\em non-empty} canalizing word from either $(j,0)$ or $(j,1)$ and, by carefully analyzing the inductive calls, we obtain $|u|\leq F(n)$ in the first case, and $|vw|\leq F(n)$ in the second. So $|uvw|<F(n)+F(n+1)=F(n+2)$ in any case.

\medskip
We then suppose that $i$ is of out-degree at least one. We consider a longest canalizing word $v$ from $(i,0)$. If the word $vi$ increases the state of $i$ then we are done since $|v|+1< F(n-|v|+2)+|v|$. Otherwise, we ``remove'' the vertices canalized by $v$ from the network, and we ``remove'' the out-going arcs of $i$, considering that each out-neighbor of $i$ behaves as if the state of $i$ was permanently $0$.  We obtain a ``subnetwork'' whose interaction graph $H$ satisfies all the hypotheses, and in which $i$ is of out-degree zero. By the first case, we get, for the subnetwork, a word $w$ that increases the state of $i$ and we are done: starting from $x$, the word $vw$ increases the state of $i$. Suppose now that $f$ is an and-or-net. Since $H$ has $n-|v|$ vertices we have $|w|< F(n-|v|+2)$ by the first case, so $|vw|< F(n-|v|+2)+|v|$ as desired. 

\medskip
We now proceed to the details.

\begin{proof}[\BF{Proof of Lemmas~\ref{lem:no_fixed_vertex} and \ref{lem:no_fixed_vertex_2}}]
Suppose that $G$ is $i$-homogenous and has no positive cycles and no sources. Let $f$ be a BN on $G$. Let $x$ be a configuration on $V$. We proceed by induction on the number $n$ of vertices in $G$. If $n=1$ the result is obvious: $i$ is the unique vertex and has a negative loop, so $f(x)\neq x$; hence we can take $w=i$ and since $1<F(3)=2$ we are done. So suppose that $n\geq 2$. We consider two cases. 

\bigskip
\BF{Case 1: $i$ is of out-degree zero}. We have to prove that $w$ exists, and $|w|< F(n+2)$ if $f$ is an and-or-net (since $i$ is of out-degree zero, the longest canalizing word from $(i,x_i)$ is the empty word). Since $i$ is of in-degree at least one, $G$ contains an initial component $S$ from which $i$ is reachable. Since $i$ is of out-degree zero, all the paths from $S$ to $i$ are forward. Since $G$ is $i$-homogenous, there is a vertex $j$ in $S$ such that all the paths from $j$ to $i$ have the same sign. Let $a=x_i+1$ if all the paths from $j$ to $i$ are positive, and $a=x_i$ otherwise. Let $m$ be the maximal size of a canalizing word from $(j,x_j)$ with only letters in $S$. We first prove, using the induction hypothesis, the following.

\bigskip
\noindent
(1) {\em There is a word $u$ such that $f^u(x)_j=a$, and if $f$ is an and-or-net then}
\[
|u|< F(n-m+1)+m\leq F(n+1).
\]

\medskip
\noindent
This is obvious if $x_j=a$ (we take $u$ the empty word) so suppose that $x_j\neq a$. Let $J$ be the vertex set of $S$ and $I=V\setminus J$. Let $h$ be the BN with component set $J$ defined by $h(y_J)=f(y)_J$ for all configurations $y$ on $V$ such that $y_I=x_I$. Since $G$ has no arc from $I$ to $J$, one can easily check that $S$ is the signed interaction digraph of $h$. Since $S$ is strong, $S$ is $j$-homogenous. Since $S$ is a subgraph of $G$, it has no positive cycles, and since $G$ has no sources, $S$ has no sources (that is, $S$ is non-trivial). Because $S$ has at most $n-1$ vertices, by induction, there is a word $u$ over $J$ such that $h^u(x_J)_j=a$ and we can then easily check that $f^u(x)_j=a$. If $f$ is an and-or-net then $h$ is also an and-or-net, and since $m$ is the maximal length of a canalizing word from $(j,x_j)$ in $h$, by induction we obtain $|u|< F(|J|-m+2)+m$. Since $|J|\leq n-1$ this proves (1).

\bigskip
Let $u$ be a shortest word as in (1). Let $v$ be a longest canalizing word from $(j,a)$ and $b$ its image. 

\bigskip
\noindent
(2) {\em If $i\in\{v\}$ then $f^{uv}(x)_i\neq x_i$ and, if $f$ is an and-or-net, then $|uv|< F(n+2)$.}

\medskip
\noindent
Suppose that $i\in\{v\}$. By Lemma~\ref{lem:diffusion}, $G$ has a path $P$ from $j$ to $i$ which is positive if and only if $a=b_i$. By the definition of $a$, we have $a\neq x_i$ if and only if $P$ is positive. Thus $x_i\neq b_i$  and since $f^u(x)_j=a$, we have $f^{uv}(x)_i=b_i\neq x_i$. If $f$ is an and-or-net, then $|u|< F(n+1)$ by (1) and since $|v|\leq n-1\leq F(n)$ we obtain $|uv|< F(n+2)$. This proves~(2). 

\bigskip
By (2) we can suppose that $i\not\in\{v\}$. Let $b'$ be the configuration on $I=\{j\}\cup\{v\}$ such that $b'_j=a$ and ${b'}_{\{v\}}=b$. Note that, for all configurations $y$ on $V$ with $y_j=a$ we have $f^v(y)_I=b'$. In particular, since $f^u(x)_j=a$, we have $f^{uv}(x)_I=b'$. 

\medskip
Let $h$ be the BN with component set $J=V\setminus I$ defined by $h(y_J)=f(y)_J$ for all configurations $y$ on $V$ such that $y_I=b'$. Let $H$ be the signed interaction digraph of $h$.

\bigskip
\noindent
(3) {\em $H$ has no positive cycles and no sources.}

\medskip
\noindent
One easily check that $H$ is a subgraph of $G$, and thus it has no positive cycles. Suppose, for a contradiction, that $H$ has a source $\ell$. Then it means that, for some $c$, we have $h_\ell(y_J)=f_\ell(y)=c$ for all configurations $y$ on $V$ with $y_I=b'$. But then, for all configurations $y$ on $V$ with $y_j=a$, we have $f^{v,\ell}(y)_I=f^v(y)_I=b'$ since $\ell\not\in I$ and thus $f^{v,\ell}(y)_\ell=f_\ell(f^v(y))=c$. Thus $v,\ell$ is a canalizing word from $(j,a)$ longer than $v$, a contradiction. This proves (3). 

\bigskip
\noindent
(4) {\em $H$ is $i$-homogenous.}

\medskip
\noindent
Suppose, for a contradiction, that $H$ is not $i$-homogenous. Then it has a strong component $F$ such that, for every vertex $\ell$ in $F$, $H$ has both a positive path and a negative path from $\ell$ to~$i$. Since $G$ is $i$-homogenous, $F$ is not a strong component of $G$ and thus $G$ has an arc from some vertex $k$ not in $F$ to some vertex $\ell$ in $F$. Let $y,z$ be two configurations on $V$ that only differ in component $k$ such that $f(y)_\ell\neq f(z)_\ell$. If $G$ has no arc from $I$ to $\ell$, then we can choose $y,z$ in such a way that $y_I=z_I=b'$. But then $h(y_J)_\ell\neq h(z_J)_\ell$ thus $H$ has an arc from $k$ to $\ell$, a contradiction. Hence $G$ has an arc from some vertex $k'\in I$ to $\ell$. By Lemma~\ref{lem:diffusion}, $G$ has a path from $j$ to $k'$ whose vertices are all in $I$, and thus it has a path $P$ from $j$ to $\ell$ whose vertices are in $I$ excepted $\ell$. Since $\ell$ is in $F$, $H$ has a positive path $P^+$ from $\ell$ to $i$ and a negative path $P^-$ from $\ell$ to $i$. Then $P\cup P^+$ and  $P\cup P^-$ are paths of $G$ from $j$ to $i$ of opposite signs, and this contradicts our choice of $j$. This proves (4). 

\bigskip
From (3) and (4) there is, by induction, a word $w$ over $J$ with $h^w(f^{uv}(x)_J)_i\neq x_i$. Since $f^{uv}(x)_I=b'$, and $w$ is a word over $J$, we have 
\[
f^{uvw}(x)_J=f^w(f^{uv}(x))_J=h^w(f^{uv}(x)_J)
\]
thus $f^{uvw}(x)_i\neq x_i$. 

\medskip
It remains to prove the upper bound on the length of $uvw$ when $f$ is an and-or-net. So suppose that $f$ is an and-or-net. Then $h$ is also an and-or-net, and since $H$ has $n-|v|-1$ vertices, we have $|w|<F(n-|v|+1)$ by induction, and thus 
\[\tag{5}
|vw|< F(n-|v|+1)+|v|\leq F(n+1). 
\] 
If $|u|\leq 1$ then $|uvw|< F(n+1)+1\leq F(n+2)$ as desired. So suppose that $|u|\geq 2$. Since $u$ is as short as possible, we have $x_j\neq a$ and $S$ has at least two vertices, so $j$ has an out-neighbor in $S$ distinct from $j$, say $k$. Since $f_k$ is a conjunction or a disjunction, the single letter $k$ is a canalizing word from $(j,c)$ for some $c\in\B$. If $c=x_j$ then we deduce that $m\geq 1$, so $|u|\leq F(n)$ by (1), and with (5) we obtain $|uvw|< F(n+2)$. If $c=a$ then we deduce that $|v|\geq 1$, so $|vw|\leq F(n)$ by (5) and with (1) we obtain $|uvw|< F(n+2)$. Hence the desired bound is always obtained. This proves the first case. 

\bigskip
\BF{Case 2: $i$ is of out-degree at least one.} 
Let $a=x_i$ and let $v$ be a longest canalizing word from $(i,a)$. If $f^{v,i}(x)_i\neq a$ then we are done since $|v,i|\leq n$ and $n< F(n-m+2)+m$ for all $0\leq m<n$. So we assume that $f^{v,i}(x)_i=a$. Let $b$ be the image  of $v$ and $I=\{v\}$.

\bigskip
\noindent
(1) {\em There is no $c$ such that $f_i(y)=c$ for all configurations $y$ on $V$ with $y_i=a$ and $y_I=b$.}

\medskip
\noindent
Suppose, for a contradiction, that $f_i(y)=c$ for all configurations $y$ on $V$ with $y_i=a$ and $y_I=b$. Let $y$ be a configuration on $V$ with $y_i=a$  such that, for all positive (negative) in-neighbors $j$ of $i$ with $j\not\in I$ we have $y_j\neq a$ ($y_j=a)$; if $j=i$ then $j$ is a negative in-neighbor of $i$ since $G$ has only negative cycles, and we have already imposed $y_j=a$. Suppose that $i$ has a positive in-neighbor $j\in I$. Since $G$ has only negative cycles, all the paths from $i$ to $j$ are negative. Hence, by Lemma~\ref{lem:diffusion}, we have $a\neq b_j$. Thus, for {\em every} positive in-neighbor $j$ of $i$ we have $f^v(y)_j\neq a$, because $f^v(y)_j=y_j\neq a$ if $j\not\in I$ and $f^v(y)_j=b_j\neq a$ otherwise. We prove similarly that, for every negative in-neighbor $j$ of $i$, we have $f^v(y)_j=a$. Hence $f^v(y)$ is an $i$-unstable configuration in $G$, and since $i$ is not a source, by Proposition~\ref{pro:flipping} we have $f_i(f^v(y))\neq a$. Since $f^v(y)_i=y_i=a$ and $f^v(y)_I=b$, we deduce that $c\neq a$. Now, since $f^v(x)_i=x_i=a$ and $f^v(x)_I=b$, we have $f_i(f^v(x))=c\neq a$ and thus $f^{v,i}(x)_i\neq a$, which contradicts our hypothesis. This proves (1). 

\bigskip
Given a configuration $y$ on $V$ we denote by $\tilde y$ the configuration on $V$ such that $\tilde y_i=a$ and $\tilde y_j=y_j$ for all $j\neq i$. Let $h$ be the BN with component set $J=V\setminus I$ such that $h(y_J)=f(\tilde y)_J$ for all configurations $y$ on $V$ with $y_I=b$. Let $H$ be the signed interaction digraph of $h$. Note that $H$ is a subgraph of $G$ in which $i$ is of out-degree zero.

\bigskip
\noindent
(2) {\em $H$ has no positive cycles and no sources.}

\medskip
\noindent
Since $H$ is a subgraph of $G$, it has no positive cycles. Suppose, for a contradiction, that $H$ has a source $k$. Then there is $c$ such that $h(y_J)_k=c$ for all configurations $y$ on $V$. So if $y_i=a$ and $y_I=b$ then $f(y)_J=h(y_J)$ and thus $f(y)_k=h(y_J)_k=c$. For all configurations $y$ on $V$ with $y_i=a$, we have $f^{v}(y)_i=a$ and $f^v(y)_I=b$ and we deduce that $f^{v,k}(y)_k=c$. By (1) we have $k\neq i$ and since $k\not\in I$ we deduce that $v,k$ is a canalizing word from $(i,a)$ longer than $v$, a contradiction. This proves (2). 

\bigskip
\noindent
(3) {\em $H$ is $i$-homogenous.}

\medskip
\noindent
Suppose, for a contradiction, that $H$ is not $i$-homogenous. Then it has a strong component $F$ such that, for any vertex $\ell$ in $F$, $H$ has both a positive path and a negative path from $\ell$ to $i$. Since $G$ is $i$-homogenous, $F$ is not a strong component of $G$ and thus $G$ has an arc from some vertex $k$ not in $F$ to some vertex $\ell$ in $F$. Let $y,z$ be two configurations on $V$ that only differ in component $k$ such that $f(y)_\ell\neq f(z)_\ell$. If $G$ has no arc from $I\cup\{i\}$ to $\ell$, then we can choose $y,z$ in such a way that $y_I=z_I=b$ and $y_i=z_i=a$. But then $h(y_J)_\ell=f(y)_\ell\neq f(z)_\ell= h(z_J)_\ell$ thus $H$ has an arc from $k$ to $\ell$, a contradiction. Hence $G$ has an arc from some vertex $k'\in I\cup\{i\}$ to $\ell$. By Lemma~\ref{lem:diffusion}, $G$ has a path from $i$ to $k'$ whose internal vertices are in $I$ (of length zero if $k=i$), and thus it has a path $P$ from $i$ to $\ell$ whose internal vertices are in $I$. Since $\ell$ is in $F$, $H$ has a positive path $P^+$ from $\ell$ to $i$ and a negative path $P^-$ from $\ell$ to $i$. Then $P\cup P^+$ and  $P\cup P^-$ are cycles of $G$ of opposite signs, so $G$ has a positive cycle, a contradiction. This proves~(3).

\bigskip
Since $i$ is of out-degree zero in $H$, we deduce from (2), (3) and the first case that there is a shortest word $w$ over $J$ such that $h^w(z_J)_i\neq a$ where $z=f^{v}(x)$. Since $w$ is of minimal length for this property, $i$ appears only once in $w$, in last position: $w=w',i$ for some word $w'$ not containing $i$. Since $z_I=b$, $z_i=a$ and $w'$ does not contain any letter in $I\cup\{i\}$, we deduce, setting $y=f^{w'}(z)$, that $y_J=h^{w'}(z_J)$. Furthermore, $y_I=b$ and $y_i=a$ so $f(y)_J=h(y_J)$. Consequently, 
\[
f^{vw}(x)_i=f(y)_i=h(y_J)_i=h(h^{w'}(z_J))_i=h^w(z_J)_i\neq a.
\] 

\medskip
It remains to prove the upper bound on the length of $vw$ when $f$ is an and-or-net. So suppose that $f$ is an and-or-net. Then $h$ is also an and-or-net, and since $H$ has $n-|v|$ vertices, we have $|w|< F(n-|v|+2)$ by the first case, and thus $|vw|< F(n-|v|+2)+|v|$ as desired. This proves the second case. 
\end{proof}

\begin{remark}\label{rem:att_size}
Suppose that $G$ is strong, non-trivial, and has no positive cycles. Let $f$ be any BN on $G$. It is proved in \cite{RC07} that the state diagram $\Gamma(f)$ of the asynchronous automaton of $f$ has a unique terminal strong component $A$. Furthermore, it is proved in \cite{A08} that $f$ has no fixed points, and thus $A$ contains at least $2$ configurations. But Lemma~\ref{lem:no_fixed_vertex} says something stronger: $A$ contains at least $n+1$ configurations. Indeed, for every configuration $x$ in $A$ and $i\in V$ we have $f^w(x)_i\neq x_i$ for some word $w$, and since $A$ is a terminal strong component of $\Gamma(f)$, we have $f^w(x)\in A$. Taking $w$ as short as possible, $i$ is the last letter of $w$, so $w=v,i$ for some word $v$, and $f^v(x)_i\neq f^w(x)_i$. Since $A$ is a terminal strong component of $\Gamma(f)$, the configurations $f^v(x)$ and $f^w(x)$ are both contained in $A$, and thay only differ in component $i$. This shows that, for every $i\in V$, there are two configurations in $A$ which differ only in component $i$ and this implies that $A$ contains at least $n+1$ configurations (see Appendix \ref{app:att_size} for details). We conjecture that the right lower bound is $2n$.
\end{remark}

%%%%%%%%%%%%%%%%%%%%%%%%%%%%%%%%%%%%%%%%%%%%%%%%%%%%%%%%%%%%%%%%%%%%%%%%%%%%%%%%%%%%%%%%%%%%%%%%%%%%%%%
\subsection{Synchronizing a vertex}
%%%%%%%%%%%%%%%%%%%%%%%%%%%%%%%%%%%%%%%%%%%%%%%%%%%%%%%%%%%%%%%%%%%%%%%%%%%%%%%%%%%%%%%%%%%%%%%%%%%%%%%

Suppose that $G$ is $i$-homogenous and has no positive cycles and no sources. Let $f$ be the and-net on $G$. In this subsection, we use Lemma~\ref{lem:no_fixed_vertex_2} to prove a ``local'' synchronization: given two configurations $x,y$, there is a word $w$ synchronizing $i$ at state $0$, that is, such that $f^w(x)_i=f^w(y)_i=0$. This will be the starting point for the ``global'' synchronization given in the next subsection. The argument is roughly the following. If $x_i=y_i=0$ there is nothing to prove: we can take $w=\epsilon$. So suppose, without loss of generality, that $x_i=1$. The key points can be explained by assuming that $i$ has exactly two in-neighbors, say $j$ and $\ell$, both positive and distinct from~$i$. An easy application of Lemma~\ref{lem:no_fixed_vertex_2} shows that there is a word $u$, not containing $i$, such that $f^u(x)$ is an $i$-unstable configuration in $G$. Hence, setting $x'=f^u(x)$ and $y'=f^u(y)$ we have the following situation:
\[
\begin{tikzpicture}
\node (j) at (0,1){$x'_j=0$};
\node (l) at (0,0){$x'_\ell=0$};
\node (i) at (2,0.5){$x'_i=1$};
\path[thick,->,Green]
(j) edge (i)
(l) edge (i)
;
\end{tikzpicture}
\qquad
\begin{tikzpicture}
\node (j) at (0,1){$y'_j=?$};
\node (l) at (0,0){$y'_\ell=?$};
\node (i) at (2,0.5){$y'_i=?$};
\path[thick,->,Green]
(j) edge (i)
(l) edge (i)
;
\end{tikzpicture}
\]
Hence $f_i(x')=0$ and if $f_i(y')=0$ then we are done with $w=v,i$. So suppose that $f_i(y')=1$. Since $f_i$ is a conjunction, this implies $y'_j=y'_\ell=1$, so we have the following situation: 
\[
\begin{tikzpicture}
\node (j) at (0,1){$x'_j=0$};
\node (l) at (0,0){$x'_\ell=0$};
\node (i) at (2,0.5){$x'_i=1$};
\path[thick,->,Green]
(j) edge (i)
(l) edge (i)
;
\end{tikzpicture}
\qquad
\begin{tikzpicture}
\node (j) at (0,1){$y'_j=1$};
\node (l) at (0,0){$y'_\ell=1$};
\node (i) at (2,0.5){$y'_i=?$};
\path[thick,->,Green]
(j) edge (i)
(l) edge (i)
;
\end{tikzpicture}
\]
Setting $x''=f^i(x')$ and $y''=f^i(y')$ we obtain a clear situation:
\[
\begin{tikzpicture}
\node (j) at (0,1){$x''_j=0$};
\node (l) at (0,0){$x''_\ell=0$};
\node (i) at (2,0.5){$x''_i=0$};
\path[thick,->,Green]
(j) edge (i)
(l) edge (i)
;
\end{tikzpicture}
\qquad
\begin{tikzpicture}
\node (j) at (0,1){$y''_j=1$};
\node (l) at (0,0){$y''_\ell=1$};
\node (i) at (2,0.5){$y''_i=1$};
\path[thick,->,Green]
(j) edge (i)
(l) edge (i)
;
\end{tikzpicture}
\]
By Lemma~\ref{lem:no_fixed_vertex_2}, there is a shortest word $v$ such that $f^v(y'')_i=0$. Since $y''_i=1$, $v$ has a shortest prefix $\tilde v$ whose last letter is $j$ or $\ell$. Suppose that this is $\ell$. Then $i$ and $j$ do not appear in $\tilde v$ and thus, setting $x'''=f^{\tilde v}(x'')$ and $y'''=f^{\tilde v}(y'')$, we obtain the following situation:
\[
\begin{tikzpicture}
\node (j) at (0,1){$x'''_j=0$};
\node (l) at (0,0){$x'''_\ell=?$};
\node (i) at (2,0.5){$x'''_i=0$};
\path[thick,->,Green]
(j) edge (i)
(l) edge (i)
;
\end{tikzpicture}
\qquad
\begin{tikzpicture}
\node (j) at (0,1){$y'''_j=1$};
\node (l) at (0,0){$y'''_\ell=0$};
\node (i) at (2,0.5){$y'''_i=1$};
\path[thick,->,Green]
(j) edge (i)
(l) edge (i)
;
\end{tikzpicture}
\]
Since $x'''_j=y'''_\ell=0$ we have $f_i(x''')=f_i(y''')=0$, thus $w=u,i,\tilde v,i$ has the desired properties. Furthermore, from the bound given in Lemma~\ref{lem:no_fixed_vertex_2}, we obtain $|u|,|v|<F(n+2)-1$ so that $|w|<2F(n+2)$. The case where $i$ has more than two in-neighbors, with possibly some negative, is very similar. To treat the case where $i$ has a loop or is of in-degree one, some additional easy arguments have to be given. 

\medskip
We now proceed to the details, starting with the full statement of the ``local'' synchronization. 

%In this subsection, we use Lemma~\ref{lem:no_fixed_vertex_2} to prove the following.

\begin{lemma}\label{lem:fixing_a_vertex}
Suppose that $G$ is $i$-homogenous and has no positive cycles and no sources. Suppose also that $G$ has a positive path $P$ from a vertex $j$ to $i$ such that $j$ is of in-degree at least two, and all the vertices of $P$ distinct from $j$ are of in-degree one (if $i$ is of in-degree at least two, this path always exists and is the trivial graph whose unique vertex is $i$). Let $f$ be the and-net on $G$ and let $x,y$ be configurations on $V$. There is a word $w$ such that $f^w(x)_i=f^w(y)_i=0$ and 
\[
|w|< 2F(n+2).
\] 
\end{lemma}

We need the following lemma, which shows that any configuration can be sent by some word onto an $i$-unstable configuration. 

\begin{lemma}\label{lem:coherent_local_configuration}
Suppose that $G$ is $i$-homogenous and has no positive cycles and no sources. Let $f$ be an and-or-net on $G$, let $x$ be a configuration on $V$, and let $r$ be the maximum length of a canalizing word from $(i,x_i)$. There is a word $w$, which does not contain $i$, such that $f^w(x)$ is an $i$-unstable configuration in $G$ and 
\[
|w|< F(n-r+1)+r-1.
\]
\end{lemma}

\begin{proof}
Suppose that $x_i=0$. Let $h$ be the and-or-net on $G$ such that $h_j=f_j$ for all $j\neq i$ and $h_i$ is a conjunction. By Lemma \ref{lem:no_fixed_vertex_2}, there is a shortest word $v$ such that $h^v(x)_i=1$ and $|v|<F(n-r+1)+r$ (because $r$ does not depend on $f_i$). Since $v$ is a shortest word, $i$ appears once in $w$, in last position. Let $w$ obtained from $v$ by removing the last letter~$i$. Then $h_i(h^w(x))=1$, and since $h_i$ is a conjunction, we deduce that $h^w(x)$ is an $i$-unstable configuration in $G$. Since $i\not\in\{w\}$ and $h_j=f_j$ for $j\neq i$, we have $h^w(x)=f^w(x)$. Thus $w$ has the desired properties. If $x_i=1$ the proof is similar, excepted that $h_i$ is a disjunction. 
\end{proof}

\begin{proof}[{\BF{Proof of Lemma~\ref{lem:fixing_a_vertex}}}]
If $x_i=y_i=0$ there is nothing to prove, so suppose, without loss of generality, that $x_i=1$. Let $r$ be the maximal length of a canalizing word from $(i,1)$. We prove that there is a word $w$ such that $f^w(x)_i=f^w(y)_i=0$ with 
\begin{enumerate}
\item $|w|\leq 2F(n-r+2)+2r-2$ if $i$ is of in-degree at least two,
\item $|w|\leq 2F(n+2)-1$ if $i$ is of in-degree one. 
\end{enumerate}

\medskip
Suppose first that $i$ is of in-degree at least two. By Lemma \ref{lem:coherent_local_configuration} there is a word $u$, which does not contain $i$, such that $f^u(x)$ is an $i$-unstable configuration in $G$ and 
\[
|u|\leq F(n-r+2)+r-2.
\]
Let $x'=f^u(x)$ and $y'=f^u(y)$. Since $x'$ is $i$-unstable, we have $f(x')_i=0$. Hence if $f(y')_i=0$ we are done by taking $w=u,i$. So suppose that $f(y')_i=1$. Let $x''=f^i(x')$ and $y''=f^i(y')$. Hence we have $x''_i=0$ and $y''_i=1$. We consider two cases. Suppose first that $i$ has a loop. Then it is negative thus $f(y'')_i=0$. Let $j$ be an in-neighbor of $i$ distinct from $i$. Then $x''_j=x'_j$ and since $x'$ is $i$-unstable, we have $x'_j=0$ if $j$ is a positive in-neighbor of $i$, and $x'_j=1$ otherwise, and we deduce that $f(x'')_i=0$ (because $f_i$ is a conjunction). So we are done with $w=u,i,i$. For the second case, suppose that $i$ has no loop. By Lemma~\ref{lem:no_fixed_vertex_2}, there is a shortest word $v$  such that $f^v(y'')_i=0$ and 
\[
|v|\leq F(n-r+2)+r-1.
\]
Since $y''_i=1$ and $v$ is as short as possible, $i$ is the last letter of $v$ and at least one in-neighbor of $i$ appears in $v$. Let $\ell$ be the first in-neighbor of $i$ that appears $v$, and let $\tilde v$ the shortest prefix of $v$ containing $\ell$. Let $x'''=f^{\tilde v}(x'')$ and $y'''=f^{\tilde v}(y'')$. We have $y'''_\ell\neq y''_\ell$ since $v$ is as short as possible and $\ell$ appears once in $\tilde v$. Furthermore, $y''_\ell=y'_\ell$ since $\ell\neq i$ (because $i$ has no loop). So $y'''_\ell\neq y'_\ell$, and since $f(y')_i=1$, we deduce that  $f(y''')_i=0$ (because $f_i$ is a conjunction). Let $j$ be an in-neighbor of $i$ distinct from $\ell$. Since $j$ does not appear in $\tilde v$ and $j\neq i$ (because $i$ has no loop), we have $x'''_j=x''_j=x'_j$. Since $x'$ is an $i$-unstable configuration in $G$, we have $x'_i=0$ if $j$ is a positive in-neighbor of $i$, and $x'_i=1$ otherwise. Consequently, $f(x''')_i=0$ (because $f_i$ is a conjunction). So we are done with $w=u,i,\tilde v,i$, observing that $|\tilde v|\leq |v|-1$.  

\medskip
Suppose now that $i$ is of in-degree one, and let $P$ be as in the statement. Let $s$ be the maximal length of a canalizing word from $(j,1)$. Since $G$ is $i$-homogenous, it is also $j$-homogenous and, by the first case, there is word $u$ such that $f^u(x)_j=f^u(y)_j=0$ and $|u|\leq 2F(n-s+2)+2s-2$. Let $v$ be an enumeration of the vertices of $P\setminus j$ in order. Since all the vertices of $P\setminus j$ are of in-degree one, $v$ is a canalizing word from $(j,1)$ and $(j,0)$. Thus $|v|\leq s$ and, since $P$ is positive and $f^u(x)_j=f^u(y)_j=0$, we obtain $f^{uv}(x)_i=f^{uv}(y)_i=0$ (by Lemma \ref{lem:diffusion}). Furthermore, 
\[
|uv|\leq 2F(n-s+2)+3s-2.
\]
So it is sufficient to prove that, for all $n\geq 1$ and $0\leq s<n$, 
\[
2F(n-s+2)+3s-2\leq 2F(n+2)-1.
\]
We proceed by induction on $n$. The case $n=1$ is obvious, so suppose that $n\geq 2$. If $s=n-1$ then the inequality to prove becomes $3n\leq 2F(n+2)$ which is easy to check. So suppose that $0\leq s\leq n-2$. Then, using the induction for the first inequality, we obtain
\[
\begin{array}{lll}
2F(n-s+2)+3s-2
&=&2F((n-1)-s+3)+3s-2\\
&=&2F((n-1)-s+1)+2F((n-1)-s+2)+3s-2\\
&\leq & 2F((n-1)-s+1)+2F((n-1)+2)-1\\
&\leq & 2F(n)+2F(n+1)-1\\
&= & 2F(n+2)-1,
\end{array}
\]
completing the induction. 
\end{proof}

%%%%%%%%%%%%%%%%%%%%%%%%%%%%%%%%%%%%%%%%%%%%%%%%%%%%%%%%%%%%%%%%%%%%%%%%%%%%%%%%%%%%%%%%%%%%%%%%%%%%%%%
\subsection{Synchronizing two configurations}
%%%%%%%%%%%%%%%%%%%%%%%%%%%%%%%%%%%%%%%%%%%%%%%%%%%%%%%%%%%%%%%%%%%%%%%%%%%%%%%%%%%%%%%%%%%%%%%%%%%%%%%

Suppose that $G$ is $i$-homogenous and has no positive cycles, no sources and no initial cycles. Let $f$ be the and-net on $G$. In this subsection, we use Lemmas~\ref{lem:no_fixed_vertex_2} and \ref{lem:fixing_a_vertex} to prove a more ``global'' synchronization property: given two configurations $x,y$, there is a word $w$ synchronizing $x$ and $y$, that is, such that $f^w(x)=f^w(y)$. As shown in the next subsection, this easily implies that every and-or-net on $G$ is synchronizing. The argument is roughly the following. Since $G$ has no positive cycles, it has a vertex $i_1$ with only negative out-neighbors. By Lemma \ref{lem:fixing_a_vertex}, there is a word $v^1$ of length at most $2F(n+2)$ synchronizing $i_1$ at state $0$. We then consider the ``subnetwork'' $f'$ obtained by fixing $i_1$ at state $0$ and removing $i_1$. Since $i_1$ has only negative out-neighbors, the signed interaction digraph of $f'$ is $G\setminus i_1$. If $G\setminus i_1$ still satisfies all the conditions, then we can synchronize a new vertex $i_2$, with only negative out-neighbors, with a word $v^2$ of length at most $2F(n+1)$. Repeating this argument, we possible synchronize all components, obtaining a word $w=v^1v^2\dots $ with $f^w(x)=f^w(y)$ and $|w|\leq 2F(n+2)+2F(n+1)+\cdots +2F(3)$. Note that $|w|\leq 2F(n+4)-2$ since the sum of the $n$ first Fibonacci numbers is $F(n+2)-1$. However, this synchronizing process cannot be always completed, since $G\setminus i_1$ does not necessarily satisfy the appropriate conditions: it can have initial cycles. But these initial cycles are ``controlled'' by $i_1$ and, with fastidious arguments involving Lemmas~\ref{lem:no_fixed_vertex_2}, we can construct a word $\tilde v^1$ of length $|\tilde v^1|\leq 3F(n+2)$ synchronizing both $i_1$ and the vertices that belongs to the initial cycles of $G\setminus i_1$. From this point, we can ``remove'' $i_1$ and continue the synchronizing process without trouble. We eventually get a word $w=\tilde v^1\tilde v^2\dots$ such that $f^w(x)=f^w(y)$ and $|w|\leq 3F(n+4)-3$. 

\begin{lemma}\label{lem:2-synchronization}
Suppose that $G$ is homogenous and has no positive cycles, no sources and no initial cycles. Let $f$ be the and-net on $G$. For every configurations $x,y$ on $V$, there is a word $w$ such that $f^w(x)=f^w(y)$ and 
\[
|w|\leq 3F(n+4)-3.
\]
\end{lemma}

\begin{proof}
For inductive purpose, we prove the following claim. For $0\leq m<n$, let 
\[
g(n,m)=\sum_{t=m+1}^n 3F(t+2).
\]

\bigskip
\noindent
{\bf Claim:} {\em Suppose that $G$ is homogenous and has no positive cycles and no sources. Let $I$ be the set of vertices of $G$ that belong to an initial cycle of $G$. Let $f$ be the and-net on $G$. For every configurations $x,y$ on $V$ with $x_I=y_I$, there is a word $w$ such that $f^w(x)=f^w(y)$ and}
\[
|w|\leq g(n,|I|).
\]

\bigskip
This implies the statement since if $G$ has no initial cycles then $I$ is empty, thus $x_I=y_I$ is always true, and we obtain 
\[
|w|\leq g(n,0)=\sum_{t=1}^{n} 3F(t+2)\leq 3\sum_{t=1}^{n+2} F(t)=3F(n+4)-3,
\]
where we use, for the last equality, the fact that $\sum_{t=1}^n F(t)=F(n+2)-1$ for all $n\geq 1$, which is easy to prove by induction on $n$. 

\medskip
The claim is proved by induction on $n-|I|$, that is, the number of vertices that do not belong to an initial cycle. If $n-|I|=0$ then $I=V$ so $x_I=y_I$ means $x=y$ and there is nothing to prove (we can take $w=\epsilon$). So suppose that $n-|I|\geq 1$. 
We first prove that, up to a switch, some some practical assumptions on the signs of $G$ can be made. 

\bigskip
\noindent
(1) {\em We can suppose that every $i\in V\setminus I$ with in-degree one in $G$ has a positive in-neighbor.}

\medskip
\noindent
Let $U$ be the set of vertices of $G$ of in-degree one that do not belong to $I$. Then $G[U]$ is a disjoint union of out-trees, say $T_1,\dots,T_k$. For $\ell\in [k]$, let $U_\ell$ be the vertex set of $T_\ell$. By Proposition~\ref{pro:harary}, there is $J^1_\ell\subseteq U_\ell$ such that the $J^1_\ell$-switch of $T_\ell$ is full-positive. Let $J^2_\ell=U_\ell\setminus J^1_\ell$. The $J^2_\ell$-switch of $T_\ell$ is also full-positive. Let $i_\ell$ be the source of $T_\ell$ and suppose, without loss, that $i_\ell\in J^1_\ell$. We set $J_\ell=J^1_\ell$ if the in-neighbor of $i_\ell$ in $G$ is negative, and $J_\ell=J^2_\ell$ otherwise. Hence, in the $J_\ell$-switch of $G$, the in-neighbor of each $i\in U_\ell$ is positive. Let $J=J_1\cup\dots\cup J_k$ and let $G'$ be the $J$-switch of $G$. Then, in $G'$, the in-neighbor of each $i\in U$ is positive. Let $f'$ be the $J$-switch of $f$. Since $J\subseteq U$, each vertex $i$ of in degree at least two in $G$ (or $G'$) is not in $J$, and thus $f'_i$ is a conjunction. So $f'$ is the and-net of $G'$. By Lemma~\ref{lem:invariant_homogeneity}, $G'$ satisfies the condition of the claim. From Proposition \ref{pro:BN_switch}, for any word $w$ we have $f^w(x)=f^w(y)$ if and only if $f'^w(x+e_J)=f'^w(y+e_J)$. Hence, without loss of generality, we can suppose that $G=G'$ and $f'=f$. This proves (1).

\bigskip
So in the following, we suppose every $i\in V\setminus I$ of in-degree one in $G$ has a positive in-neighbor. Since $n-|I|\geq 1$, $G$ has a terminal strong component $S$ which is not an initial cycle. Let $i$ be a vertex in $S$ with only negative out-neighbors; such a vertex exists since otherwise $G$ has a (full-)positive cycle. Let $H=G\setminus i$.

\bigskip
\noindent
(2) {\em $H$ has no positive cycles and no sources.}

\medskip
\noindent
Since $H$ is a subgraph of $G$ it has no positive cycles. If $j$ is a source of $H$, then $i$ is the unique in-neighbor of $j$ in $G$ and thus, by hypothesis, the arc from $i$ to $j$ is positive, which contradicts our choice of $i$. This proves (2).

\bigskip
\noindent
(3) {\em $H$ is homogenous.}

\medskip
\noindent
Let $j$ be a vertex in $H$ and suppose, for a contradiction, that $H$ has a strong component $F$ such that, for every vertex $\ell$ in $F$, $H$ has both a positive and negative forward path from $\ell$ to $j$. Since $G$ is $j$-homogenous, $S$ is not a strong component of $G$, so $G$ has an arc from $i$ to some vertex $\ell$ in $F$, and thus $H$ has a positive forward path $P^+$ from $\ell$ to $j$ and a negative forward path $P^-$ from $\ell$ to $j$. Since $i$ is in a terminal strong component of $G$ and $G$ has a path from $i$ to $j$, $G$ has also a path $P$ from $j$ to $i$. Since $P\setminus i$ is a path of $H$ starting from $j$ and $P^+,P^-$ are forward paths of $H$ ending in $j$, $P^+$ and $P^-$ intersect $P\setminus i$ in $j$ only. Thus $P^+\cup P$ and $P^-\cup P$ are paths of $G$ from $\ell$ to $i$ with opposite signs. Since $G$ has an arc from $i$ to $\ell$, we deduce that $G$ has a positive cycle, a contradiction. This proves (3).  

\bigskip
Let $K$ be the set of vertices of $H$ that belongs to an initial cycles of $H$. Since $i$ is not in an initial cycle of $G$, all the initial cycles of $G$ are initial cycles of $H$, so $I\subseteq K$. From the induction hypothesis, we get that two configurations that ``agree'' on $i$ and $K$ can be synchronized:

\bigskip
\noindent
(4) {\em If $x,y$ are configurations on $V$ with $x_i=y_i=0$ and $x_K=y_K$, then there is a word $w$ such that $f^w(x)=f^w(y)$ and $|w|\leq g(n-1,|K|)$.}

\medskip
\noindent
Let $x,y$ be configurations on $V$ with $x_i=y_i=0$ and $x_K=y_K$. Let $h$ be the and-net on $H$. From (2), (3) and the induction hypothesis, there is a word $w$ on $V'=V\setminus i$ such that $h^w(x_{V'})=h^w(y_{V'})$ of length at most $g(n-1,|K|)$. Since $x_i=y_i=0$ and $i$ has only negative out-neighbors and $i\not\in\{w\}$, we deduce that $f^w(x)_{V'}=h^w(x_{V'})$ and  $f^w(y)_{V'}=h^w(y_{V'})$ and thus $f^w(x)=f^w(y)$. This proves (4). 

\bigskip
So, by (4), to prove the claim it is sufficient to prove that, for every configurations $x,y$ on $V$ with $x_I=y_I$, there is a word $w$ such that $f^w(x)_i=f^w(y)_i=0$ and $f^w(x)_K=f^w(x)_K$. This is done in several steps. 

\bigskip
\noindent
(5) {\em If $x,y$ are configurations $x,y$ on $V$ with $x_I=y_I$, then $f^w(x)_I=f^w(y)_I$ for any word $w$.}

\medskip
\noindent
Let $x,y$ be two configurations on $V$ with $x_I=y_I$. It is sufficient to prove that, for any $j\in V$, we have $f^j(x)_I=f^j(y)_I$. This is obvious if $j\not\in I$. Otherwise, $j$ has a unique in-neighbor, say $k$, and since $k\in I$ we have $x_k=y_k$ and we deduce that $f_j(x)=f_j(y)$ and so $f^j(x)_I=f^j(y)_I$. This proves (5). 

\bigskip
\noindent
(6) {\em Each initial cycle of $H$ which is not an initial cycle of $G$ has a unique negative arc.}

\medskip
\noindent
Suppose, for a contradiction, that $H$ has an initial cycle $C$, which is not an initial cycle of $G$, such that $C$ does not contain a unique negative arc. Since $C$ is negative, it has at least three negative arcs. Let $(j_k,i_k)$, $k=1,2,3$, be three negative arcs in $C$, and let $P_k$ be the path from $i_k$ to $i_{k+1}$ contained in $C$, where $i_4$ means $i_1$. Hence $C=P_1\cup P_2\cup P_3$. The three negative arcs can be chosen consecutively in $C$, that is, in such a way that $C=P_1\cup P_2\cup P_3$ and, for $k=1,2$, $P_k$ has a unique negative arc, which is $(j_{k+1},i_{k+1})$. Thus $P_1,P_2$ are negative, and so is $P_3$. Also, there is a negative arc from $i$ to each of $i_1,i_2,i_3$. Indeed, if there is no arc from $i$ to $i_k$ then $i_k$ is of in-degree one in $G$, so $(j_k,i_k)$ is positive by the hypothesis resulting from (1), a contradiction. Thus $(i,i_k)$ exists and is negative by the choice of $i$. Let $Q$ be a shortest path of $G$ from $C$ to $i$ (it exists since $i$ is in a terminal strong component of $G$) and let $j$ be its initial vertex. Then $j$ belongs to $P_k\setminus i_{k+1}$ for some $1\leq k\leq 3$. Let $P'_k$ be the path from $i_k$ to $j$ contained in $P_k$. Since $P_{k-1}$ is negative ($P_0$ means $P_3$), $P'_k$ and $P_{k-1}\cup P'_k$ have distinct signs. Since $G$ has a negative arc from $i$ to $i_k$ and $i_{k-1}$ ($i_0$ means $i_3$) we deduce that $G$ has a positive path $P^+$ from $i$ to $j$ and a negative path $P^-$ from $i$ to $j$ such that $P^+,P^-$ are internally vertex-disjoint from $Q$. Hence, $P^+\cup Q$ and $P^-\cup Q$ are cycles with distinct signs, and thus $G$ has a positive cycle, a contradiction. This proves (6).

\bigskip
\noindent
(7) {\em Let $C$ be an initial cycle of $H$, which is not an initial cycle of $G$, and let $L$ its vertex~set. If $x,y$ are configurations on $V$ and $x_i=1$, then there is a word $w$ over $L$ such that $f^w(x)_L=f^w(y)_L=\ZERO$ and $|w|\leq 2|L|$.} 

\medskip
\noindent
By (6), $C$ has a unique negative arc, say $(j_1,i_1)$. By the hypothesis resulting from (1), $i_1$ is of in-degree at least two, so $G$ has an arc from $i$ to $i_1$, which is negative by the choice of $i$. Let $u=i_1i_2\dots i_\ell$ be an enumeration of the vertices of $C$ in the order. 

\medskip
Let $x$ be a configuration on $V$ with $x_i=1$. Since $G$ has a negative arc from $i$ to $i_1$, we have $f(x)_{i_1}=0$. Since there is a positive arc from $i_k$ to $i_{k+1}$ for all $1\leq k<\ell$, we deduce that $f^u(x)_L=\ZERO$. Furthermore, since $i\not\in\{u\}$ we have $f^u(x)_i=x_i=1$ and the same argument shows that $f^{uu}(x)_L=\ZERO$ by $(a)$. Consequently,
\[\tag{$a$}
x_i=1\quad\Rightarrow\quad f^u(x)_L=\ZERO\quad\textrm{and}\quad f^{uu}(x)_L=\ZERO.
\]

\medskip
Consider now a configuration $x$ on $V$ with $x_i=0$. Let $x(0)=x$ and $x(k)=f^{i_k}(x(k-1))$ for $1\leq k\leq \ell$, so that $x(\ell)=f^u(x)$. Note that $x(k)_i=0$ for $0\leq k\leq\ell$ since $i\not\in\{u\}$. We will prove, by induction on $k$ from $1$ to $\ell$, that $x(k)_{i_k}=a$, where $a=x_{i_\ell}+1$. Since $i_1$ has exactly two in-neighbors, $i$ and $i_\ell$, both negative, and $x(0)_i=0$, we have $f(x(0))_{i_1}=(x(0)_i+1)\land (x(0)_{i_\ell}+1)=1\land a=a$,  and thus $x(1)_{i_1}=a$. Let $1<k\leq\ell$. By induction, $x(k-1)_{i_{k-1}}=a$. Since there is a positive arc from $i_{k-1}$ to $i_k$, and since $i_k$ has at most one possible other in-neighbor, which is $i$ and then negative, we have either $f(x(k-1))_{i_k}=x(k-1)_{i_{k-1}}=a$ if there is no arc from $i$ to $i_k$ or $f(x(k-1))_{i_k}=x(k-1)_{i_{k-1}}\land (x(k-1)_i+1)=a\land 1=a$. Thus $x(k)_{i_k}=a$ in both cases. This completes the induction step. Consequently, if $x_{i_\ell}=1$ then $f^u(x)_L=\ZERO$, and if $x_{i_\ell}=0$ then $f^u(x)_L=\ONE$. So if $x_{i_\ell}=0$ then $f^u(x)_{i_\ell}=1$ and, since $f^{u}(x)_i=x(\ell)_i=0$, the same argument shows that $f^{uu}(x)_L=\ZERO$. Consequently, 
\[\tag{$b$}
x_i=0\quad\Rightarrow\quad f^u(x)_L=\ZERO\quad\textrm{or}\quad f^{uu}(x)_L=\ZERO.
\]

\medskip
Let $x,y$ be configurations on $V$ with $x_i=1$. By $(a)$ and $(b)$ we have $f^u(x)_L=f^u(y)_L=\ZERO$ or  $f^{uu}(x)_L=f^{uu}(y)_L=\ZERO$. This proves (7).

\bigskip
Let $J=K\setminus I$. So $J$ is the set of vertices if $H$ which belongs to an initial cycle of $H$ which is not an initial cycle of $G$. 

\bigskip
\noindent
(8) {\em If $x,y$ are configurations on $V$ and $x_i=1$, then there is a word $w$ over $J$ such that $f^w(x)_J=f^w(y)_J=\ZERO$ and $|w|\leq 2|J|$.} 

\medskip
\noindent
If $J=\emptyset$ there is nothing to prove. Otherwise, $H$ has $p\geq 1$ initial cycles $C_1,\dots,C_p$ which are not initial cycles of $G$. For $1\leq k\leq p$, let $L_k$ be the vertex set of $C_k$, so that $J=L_1\cup\dots \cup L_p$. By (7), for each $1\leq k\leq p$, there is a word $u^k$ over $L_k$ with $f^{u^k}(x)_{L_k}=f^{u^k}(y)_{L_k}=\ZERO$ and $|u^k|\leq 2|L_k|$. Since, for every distinct $r,s\in [p]$, $L_r\cap L_s=\emptyset$ and $G$ has no arc between $L_r$ and $L_s$, we deduce that $f^w(x)_J=f^w(y)_J=\ZERO$ with $w=u^1,\dots,u^p$. This proves (8).

\bigskip
We are now in position to prove that, for any pair of configurations $x,y$ on $V$ with $x_I=y_I$, there is a word $w$ such that $f^w(x)_i=f^w(y)_i=0$ and $f^w(x)_K=f^w(x)_K$. We consider two cases, giving (9) and (10) below. The first uses only Lemma~\ref{lem:no_fixed_vertex_2} (giving a bound on $|w|$ of order $F(n+2)$), while the second uses  Lemmas~\ref{lem:no_fixed_vertex_2} and \ref{lem:fixing_a_vertex} (giving a bound on $|w|$ of order $3F(n+2)$). 

\bigskip
\noindent
(9) {\em Suppose that $G$ has an arc from $K$ to $i$. If $x,y$ are configurations on $V$ with $x_I=y_I$, then there is a word $w$ such that $f^w(x)_i=f^w(y)_i=0$ and $f^w(x)_K=f^w(y)_K$ with}
\[
|w|\leq F(n+2)+2|J|.
\]

\medskip
\noindent
By Lemma~\ref{lem:no_fixed_vertex_2}, there is a word $u$ over $V$ such that $f^u(x)_i=1$ and $|u|< F(n+2)$. So by (8) there is a word $v$ over $J$ such that $f^{uv}(x)_J=f^{uv}(y)_J=\ZERO$ with $|v|\leq 2|J|$. We deduce from (5) that $f^{uv}(x)_K=f^{uv}(y)_K$. Suppose that $G$ has an arc from some $j\in K$ to $i$. Then $j$ belongs to an initial cycle $C$ of $H$ which is not an initial cycle of $G$, so $j\in J$. By (6), $C$ has a unique negative arc, say $(k,\ell)$. By the hypothesis resulting from (1), $\ell$ is of in-degree at least two and thus $G$ has an arc from $i$ to $\ell$, which is negative by the choice of $i$. Hence, the path from $\ell$ to $j$ contained in $C$ is full-positive and it forms, with the arc from $i$ to $\ell$ and the arc from $j$ to $i$, a cycle, which is negative by hypothesis. We deduce that the arc from $j$ to $i$ is positive. Since $f^{uv}(x)_j=f^{uv}(y)_j=0$ we have $f_i(f^{uv}(x))=f_i(f^{uv}(y))_i=0$. Since $i\not\in K$, we deduce that $w=uv,i$ has the desired properties. This proves (9). 

\bigskip
\noindent
(10) {\em Suppose that $G$ has no arc from $K$ to $i$. If $x,y$ are configurations $x,y$ on $V$ with $x_I=y_I$, then there is a word $w$ such that $f^w(x)_i=f^w(y)_i=0$ and $f^w(x)_K=f^w(y)_K$ with}
\[
|w|\leq 3F(n+2)+2|J|.
\]

\medskip
\noindent
By Lemma~\ref{lem:no_fixed_vertex_2}, there is a word $u^1$ over $V$ such that $f^{u^1}(x)_i=1$ and $|u^1|< F(n+2)$. Let $P$ be a path of $G$ from some vertex $j$ of in-degree at least two to $i$ such that all the vertices of $P$ distinct from $j$ are of in-degree one (if $i$ is of in-degree at least two, then $P$ is the path of length zero containing~$i$); this path exists since $G$ has no sources and $i$ does not belong to an initial cycle of $G$. By the hypothesis resulting from (1), $P$ is full-positive. So, by Lemma~\ref{lem:fixing_a_vertex}, there is a shortest word $u^2$ such that $f^{u^1u^2}(x)_i=f^{u^1u^2}(y)_i=0$ with $|u^2|\leq 2F(n+2)$. Let $u$ be obtained from $u^1u^2$ by removing the last letter. Since $f^{u^1}(x)_i=1$, $u^2$ is not empty, and since $u^2$ is as short as possible we deduce that the last letter is $i$. So $u,i=u^1u^2$ and since $u^2$ is as short as possible, $f^u(x)_i=1$ or $f^u(y)_i=1$. Suppose that $f^u(x)_i=1$, the other case being symmetric. By (8) there is a word $v$ over $J$ such that $f^{uv}(x)_J=f^{uv}(y)_J=\ZERO$ and $|v|\leq 2|J|$. We deduce from (5) that $f^{uv}(x)_K=f^{uv}(y)_K$. Since $G$ has no arc from $K$ to $i$, we have $f^{uv,i}(x)_i=f^{u,i}(x)_i=0$ and $f^{uv,i}(y)_i=f^{u,i}(y)_i=0$. Since $i\not\in K$, we deduce that $w=uv,i$ has the desired property. This proves (10).

\bigskip
We can now prove the claim, by combining (4) with (9) and (10). Let $x,y$ be configurations on $V$ with $x_I=y_I$. By (9) and (10) there is a word $u$ such that $f^u(x)_i=f^u(y)_i=0$ and $f^u(x)_K=f^u(y)_K$ with 
\[
|u|\leq 3F(n+2)+2|J|.
\]
We deduce from (4) that there is a word $v$ such that $f^{uv}(x)=f^{uv}(y)$ with 
\[
|v|\leq g(n-1,|I|+|J|).
\]
If $|J|=0$ then 
\[
|uv|\leq 3F(n+2)+g(n-1,|I|)=g(n,|I|).
\]
If $|J|\geq 1$, then using the fact that $2m\leq 3F(m+2)$ for all $m\geq 1$ to get the second inequality, we obtain 
\[
\begin{array}{rcl}
|uv|&\leq& 	3F(n+2)+2|J|+g(n-1,|I|+|J|) \\[1mm]
&=& 		g(n,|I|+|J|)+2|J| \\[1mm]
&\leq &		g(n,|I|+|J|)+3F(|I|+|J|+2)\\[1mm]
&=&			g(n,|I|+|J|-1)\\[1mm]
&\leq &		g(n,|I|).
\end{array}
\] 
This proves the claim. 
\end{proof}

%%%%%%%%%%%%%%%%%%%%%%%%%%%%%%%%%%%%%%%%%%%%%%%%%%%%%%%%%%%%%%%%%%%%%%%%%%%%%%%%%%%%%%%%%%%%%%%%%%%%%%%
\subsection{Global synchronization}
%%%%%%%%%%%%%%%%%%%%%%%%%%%%%%%%%%%%%%%%%%%%%%%%%%%%%%%%%%%%%%%%%%%%%%%%%%%%%%%%%%%%%%%%%%%%%%%%%%%%%%%

A classical result of \v{C}ern\'y \cite{C64} is that, if any two states of a deterministic finite automaton can be sent to the same state by some word, then this automaton is synchronizing. Here is an adaptation of this observation to our context. 

\begin{lemma}\label{lem:2-synch}
Let $f$ be a BN with component set $V$ and suppose that, for every configurations $x,y$ on $V$, there is a word $u$ such that $f^u(x)=f^u(y)$ and $|u|\leq k$. For every non-empty subset $X\subseteq \B^V$, there is a word $w$ such that $|f^w(X)|=1$ and $|w|\leq k(|X|-1)$. Taking $X=\B^V$, we deduce that $f$ has a synchronizing word of length at most $k(2^n-1)$.
\end{lemma}

\begin{proof}
We proceed by induction on $|X|$. If $|X|=1$ the result is obvious, since we can take $w=\epsilon$. So suppose that $|X|\geq 2$. Let $x,y\in X$, distinct, and let $u$ be a word such that $f^u(x)=f^u(y)$ and $|u|\leq k$. Then $X'=f^u(X)$ is of size at most $|X|-1$ and thus, by induction, there is a word $v$ such that $|f^v(X')|=1$ and $|v|\leq k(|X'|-1)$. Setting $w=uv$ we obtain $|f^w(X)|=|f^v(X')|=1$ and $|w|\leq k|X'|\leq k(|X|-1)$, completing the induction. 
\end{proof}

From this observation and Lemma~\ref{lem:2-synchronization} we deduce the following. 

\begin{lemma}\label{lem:synchro_and_net}
Suppose that $G$ is homogenous and has no positive cycles, no sources and no initial cycles. The and-net on $G$ has a synchronizing word of length at most $10(\sqrt{5}+1)^n$.
\end{lemma}

\begin{proof}
Let $f$ be the and-net on $G$. From Lemmas~\ref{lem:2-synch} and \ref{lem:2-synchronization}, we deduce that $f$ has a synchronizing word of length at most $3(F(n+4)-1)(2^n-1)$. By the well known Binet's formula, we have $F(n+4)=(\varphi^{n+4}-\psi^{n+4})/\sqrt{5}$, where $\varphi=(\sqrt{5}+1)/2$ is the golden number and $\psi=1-\varphi$. We deduce that $F(n+4)\leq (\varphi^{n+4}/\sqrt{5})+1$ and thus 
\[
3(F(n+4)-1)(2^n-1)\leq \frac{3}{\sqrt{5}}\varphi^{n+4}(2^n-1)\leq \frac{3}{\sqrt{5}}\varphi^{n+4}2^n=\frac{3\varphi^4}{\sqrt{5}}(\sqrt{5}+1)^n.
\]
Since $(3\varphi^4/\sqrt{5})\sim 9.19... < 10$ this proves the lemma.
\end{proof}

Using the switch operation, we obtain Theorem \ref{thm:main2}, that we restate.

\setcounter{theorem}{3}

\begin{theorem}
Suppose that $G$ is homogenous and has no positive cycles, no sources and no initial cycles. Then every and-or-net on $G$ has a synchronizing word of length at most $10(\sqrt{5}+1)^n$.
\end{theorem}

\begin{proof}
Let $f$ be an and-or-net on $G$, and let $I$ be the set of $i\in V$ such that $f_i$ is a disjunction. Let $G'$ be the $I$-switch of $G$ and $f'$ the $I$-switch of $f$. By Proposition \ref{pro:BN_switch}, $f'$ is a BN on $G'$, and one easily check that $f'$ is actually the and-net on $G'$. We deduce from Lemma \ref{lem:invariant_homogeneity} that $G'$ is homogenous and has no positive cycles, no sources and no initial cycles. Thus, by Lemma \ref{lem:synchro_and_net}, $f'$ has a synchronizing word $w$ of length at most $10(\sqrt{5}+1)^n$ and, by Proposition \ref{pro:BN_switch}, $w$ is a synchronizing word for~$f$.
\end{proof}

%%%%%%%%%%%%%%%%%%%%%%%%%%%%%%%%%%%%%%%%%%%%%%%%%%%%%%%%%%%%%%%%%%%%%%%%%%%%%%%%%%%%%%%%%%%%%%%%%%%%%%%
\section{Proof of Theorem \ref{thm:fast}}\label{sec:fast}
%%%%%%%%%%%%%%%%%%%%%%%%%%%%%%%%%%%%%%%%%%%%%%%%%%%%%%%%%%%%%%%%%%%%%%%%%%%%%%%%%%%%%%%%%%%%%%%%%%%%%%%

Let $G$ be a signed digraph with vertex set $V$, and $n=|V|$. Suppose that $G$ satisfies the conditions of Theorem \ref{thm:fast}, that is: $G$ is strong, has no positive cycles, and is not a cycle. These conditions are invariant by switch and so, to prove the theorem, it is sufficient to prove that the and-net on some switch of $G$ has a synchronizing word of length at most $5n(\sqrt{2})^n$. This version is more convenient, because of the regularity of and-nets. The idea is to consider a switch $G'$ containing a spanning full-positive out-tree, which exists since $G'$ is strongly connected. The out-tree $T$ ``speeds up'' the synchronization of the and-net $f$ on $G'$: for instance, denoting $r$ the root of $T$ and $u$ is a topological sort of $T\setminus r$, one can easily check that $f^u(x)=f^u(y)=\ZERO$ whenever $x_r=y_r=0$. So any configurations $x,y$ can be synchronized quickly if there is a short word $v$ synchronizing the root $r$ at state $0$, that is, such that $f^v(x)_r=f^v(y)_r=0$. It appears that, with some tricky arguments, the root can by synchronized at state $0$ with a word of length at most $4n$. Consequently, any two configurations can be synchronized with a word of length at most $5n$. By \v{C}ern\'y's observation (Lemma~\ref{lem:2-synch}), we deduce that $f$ has a synchronizing word of length at most $5n(2^n-1)$. To replace $2$ by $\sqrt{2}$ in the bound, we need additional arguments. If $G$ has a loop, a direct argument shows that $f$ has a synchronizing word of length at most $7n\leq 5n(\sqrt{2})^n$. Otherwise $G'$ has no loop and no positive cycles, and so, by a recent and difficult result of Millani, Steiner and Wiederrecht \cite{MSW19}, there is a subset $I\subseteq V$ of size at most $n/2$ such that $G'\setminus I$ is acyclic. By classical arguments, the topological sort $v$ of $G'\setminus I$ sends any two configurations $x,y$ with $x_I=y_I$ on the same configuration. It follows that $v$ sends every configuration inside a set $X$ of at most $2^{n/2}$ configurations. Using \v{C}ern\'y's observation, we deduce that there is a word $w$ of length at most $5n(|X|-1)\leq 5n(\sqrt{2})^n-5n$ with $|f^w(X)|=1$. Then $v,w$ is a synchronizing word for $f$ of length at most $5n(\sqrt{2})^n-4n$ and we are done. 

\medskip
We now proceed to the details. We first show that the presence of a spanning full-positive out-tree allows a fast synchronization of two configurations. The precise statement is the following (note that it allows the presence of some positive cycles).

\begin{lemma}\label{lem:fast-2}
Suppose that $G$ is simple and has a spanning full-positive out-tree $T$, where the root $r$ has only negative in-neighbors and is of in-degree at least two. Suppose also that $G$ has no full-positive cycle, and that each positive cycle containing $r$ has at least four negative arcs. Let $f$ be the and-net on $G$. For every configurations $x,y$ on $V$ there is a word $w$ such that $f^w(x)=f^w(y)$ and $|w|\leq 5n$. Furthermore, if $r$ has a loop, then $f$ has a synchronizing word of length at most $7n$. 
\end{lemma}

\begin{proof}
Let $u$ be a topological sort of $T\setminus\{r\}$. Since $T$ is full-positive and $f$ is an and-net, for every configuration $x$ on $V$ with $x_{r}=0$ we have $f^u(x)=\ZERO$. So for every configurations $x,y$ on $V$ with $x_{r}=y_{r}=0$ the lemma holds by taking $w=u$. For the other cases, we need additional arguments. Let $x$ be a configuration on $V$. We denote by $\ZERO(x)$ the set of $i\in V$ with $x_i=0$, and $\ONE(x)=V\setminus \ZERO(x)$. We say that a configuration $x$ on $V$ is \EM{regular} if, for all arcs of $T$ from $j$ to $i$,  we have $x_j\geq x_i$. 

\bigskip
\noindent
(1) {\em For every configuration $x$ on $V$, there is a word $v(x)$ over $\ONE(x)$, without repeated letters and not containing $r$, such that $f^{v(x)}(x)$ is regular.}

\medskip
\noindent
Let $T'$ be the subgraph of $T$ obtained by removing $r$ and each vertex in $\ZERO(x)$. Let $v=v(x)$ be a topological sort of $T'$. If $v=\epsilon$ then $x=\ZERO$ or $x=e_{r}$, thus $f^v(x)=x$ is regular. So suppose that $v$ is not empty and suppose, for a contradiction, that $y=f^v(x)$ is not regular. Then $T$ has an arc from $j$ to $i$ with $y_j<y_i$. If $i$ is not in $v$, then $y_i=x_i=0$, a contradiction. So $i$ is in $v$, hence $v=v_1,i,v_2$ for some words $v_1,v_2$. Let $z=f^{v_1}(x)$. Since $i$ is not in $v_2$, we have $f_i(z)=y_i=1$. Since there is a positive arc from $j$ to $i$, we deduce that $z_j=1$. If $j$ is not in $v_2$, then $y_j=z_j=1$, a contradiction. So $j$ is in $v_2$ and this contradicts the fact that $v$ is a topological sort of $T'$. This proves (1).

\bigskip
Given a path $P$ in $G$ and a configuration $x$ on $V$, we write $x_P=\ZERO$ ($x_P=\ONE$) to means that $x_k=0$ ($x_k=1$) for all vertices $k$ in $P$. The key observation is the following. 

\bigskip
\noindent
(2) {\em Let $x$ be a regular configuration on $V$ with $x_{r}=1$. Let $j$ be an in-neighbor of~$r$. There is a word $w(x,j)$ over $\ZERO(x)$, without repeated letters, such that $f^{w(x,j)}(x)_j=1$ and such that the subgraph of $G$ induced by $\{w(x,j)\}$ has a full-positive path from each vertex to $j$.}

\medskip
\noindent
Let $I(x)$ be the set of vertices $i\in V$ such that $G$ has a full-positive path $P$ from $i$ to $j$ with $x_P=\ZERO$; we have $I(x)\subseteq \ZERO(x)$. We proceed by induction on $|I(x)|$. If $|I(x)|=0$ then $x_j=1$ and $w(x,j)$ is the empty word. So suppose that $|I(x)|\geq 1$. 

\medskip
We first prove that $f_i(x)=1$ for some $i\in I(x)$. Suppose, for a contradiction, that $f_i(x)=0$ for all $i\in I(x)$. Let $i\in I(x)$. Since $f_i(x)=0$, $i$ has an in-neighbor $k$ in $G$ such that either $x_k=0$ and $k$ is a positive in-neighbor of $i$ or $x_k=1$ and $k$ is a negative in-neighbor of $i$. Suppose first that $k$ is a negative in-neighbor of $i$. Let $Q$ be the path from $r$ to $k$ contained in $T$, which is full-positive. Since $x$ is regular and $x_k=1$, we have $x_Q=\ONE$. Since $i\in I(x)$, there is a full-positive path $P$ from $i$ to $j$ with $x_P=\ZERO$. Thus $Q$ and $P$ are disjoint full-positive paths, from $r$ to $k$ and from $i$ to $j$, respectively. Since $(k,i)$ and $(j,r)$ are negative arcs of $G$, we deduce that $G$ has a cycle with exactly two negative arcs and containing $r$, which contradicts the hypothesis of the statement. So $x_k=0$ and $k$ is a positive in-neighbor of $i$. If $k$ is in $P$ then $G$ has a full-positive cycle, a contradiction. So $k$ is not in $P$ and, by adding to $P$ the positive arc from $k$ to $i$, we obtain a full-positive path $P'$ from $k$ to $j$ with $x_{P'}=\ZERO$. Hence, $k\in I(x)$. This proves that each vertex in $I(x)$ has a positive in-neighbor in $I(x)$. But then $G[I(x)]$ has a full-positive cycle, a contradiction.

\medskip
Therefore, there is at least one $i\in I(x)$ such that $f_i(x)=1$. Let $k$ be the in-neighbor of $i$ in $T$ (it exists since $x_i<x_{r}$ and thus $i\neq r$). Since $f_i(x)=1$ we have $x_k=1$ (since $T$ is full-positive). Let $x'=f^i(x)$. Then $x'$ is regular (because $x'_i=x'_k=x_k=1$). Furthermore, $I(x')\subseteq I(x)$, since $x'\geq x$, and $i\not\in I(x')$ since $x'_i=1$. Thus, by induction, there is a word $w(x',j)$ with the properties of the statement (with $x'$ instead of $x$). Then $w(x,j)=i,w(x',j)$ has the desired properties. This proves (2). 

\bigskip
For every configuration $x$ on $V$ and every negative in-neighbor $j$ of $r$, we define $w(x,j)$ to be a word as in (2) if $x_{r}=1$, and we set $w(x,j)=\epsilon$ otherwise. 

\medskip
Let $x,y$ be configurations on $V$ with $x_{r}=1$ or $y_{r}=1$. We will prove that there is a word $w(x,y)$ of length at most $5n$ that sends both $x$ and $y$ on $\ZERO$, that is, $f^{w(x,y)}(x)=f^{w(x,y)}(y)=\ZERO$ (for the case $x_{r}=y_{r}=0$ it is sufficient to take $w(x,y)=u$ as already mentioned at the beginning of the proof). Let $j$ and $k$ be distinct in-neighbors of $r$, which exist and are negative by hypothesis. Since $G$ has no full-positive cycle, either $G$ has no full-positive paths from $j$ to $k$, or $G$ has no full-positive paths from $k$ to $j$. Suppose, without loss of generality, that $G$ has no full-positive paths from $j$ to $k$.

\medskip
Suppose that $x_{r}=y_{r}=1$. Let 
\[
\begin{array}{l}
x^1=f^{v(x)}(x),\\
y^1=f^{v(x)}(y),
\end{array}
\quad
\begin{array}{l}
x^2=f^{w(x^1,j)}(x^1),\\
y^2=f^{w(x^1,j)}(y^1),
\end{array}
\quad
\begin{array}{l}
x^3=f^{v(y^2)}(x^2),\\
y^3=f^{v(y^2)}(y^2),
\end{array}
\quad
\begin{array}{l}
x^4=f^{w(y^3,k)}(x^3),\\
y^4=f^{w(y^3,k)}(y^3).
\end{array}
\]
By (1), $x^1$ is regular, and since $v(x)$ does not contain $r$, we have $x^1_{r}=x_{r}=1$. So, by (2), $x^2_j=1$. If $y^2_j=1$ then $f_{r}(x^2)=f_{r}(y^2)=0$ so $f^u(f^{r}(x^2))=f^u(f^{r}(y^2))=\ZERO$. Thus
\[
v(x),w(x^1,j),r,u
\]
sends both $x$ and $y$ on $\ZERO$ and is of length at most $3n-2$. So suppose that $y^2_j=0$. By (1), $y^3$ is regular, and $j$ does not appears in $v(y^2)$ since $y^2_j=0$. Hence, $x^3_j=x^2_j=1$. Furthermore, $r$ does not appear in $v(x),w(x^1,j),v(y^2)$, thus $y^3_{r}=y_{r}=1$. Consequently, by (2), we have $y^4_k=1$. Since $G$ has no full-positive path from $j$ to $k$, $j$ is not contained in $w(y^3,k)$, and thus $x^4_j=x^3_j=1$. Hence, $x^4_j=y^4_k=1$. So $f_{r}(x^4)=f_{r}(y^4)=0$. Consequently, $f^u(f^{r}(x^4))=f^u(f^{r}(y^4))=\ZERO$ and we deduce that 
\[
v(x),w(x^1,j),v(y^2),w(y^3,k),r,u. 
\]
sends both $x$ and $y$ on $\ZERO$ and is of length at most $5n-4$. 

\medskip
Finally, suppose that $x_{r}<y_{r}$ (the case $x_{r}>y_{r}$ is similar by symmetry). Let 
\[
\begin{array}{l}
x^1=f^u(x)\\
y^1=f^u(y)
\end{array}
\quad
\begin{array}{l}
x^2=f^{v(y^1)}(x^1)\\
y^2=f^{v(y^1)}(y^1)
\end{array}
\quad
\begin{array}{l}
x^3=f^{w(y^2,j)}(x^2)\\
y^3=f^{w(y^2,j)}(y^2)
\end{array}
\quad
\begin{array}{l}
x^4=f^{r}(x^3)\\
y^4=f^{r}(y^3)
\end{array}
\quad
\begin{array}{l}
x^5=f^{w(x^4,k)}(x^4)\\
y^5=f^{w(x^4,k)}(y^4)
\end{array}
\]
We have $x^1=\ZERO$. Since $r$ is not in $v(y^1)$ we deduce that $x^2=\ZERO$, and since $r$ is furthermore not in $u$, we have $y^2_{r}=y_{r}=1$. By (1), $y^2$ is regular and so by (2) we have $y^3_j=1$, and since $r$ is not in $w(y^2,j)$ we have $x^3=\ZERO$ because $x^2=\ZERO$. Since $r$ has only negative in-neighbors, we have $x^4=e_{r}$, which is regular, and thus by (2) we have $x^5_k=1$. Since $G$ has no full-positive paths from $j$ to $k$, we have $j\neq r$ and since $j$ is not contained in $w(x^2,k)$, we obtain $1=y^3_j=y^4_j=y^5_j$. So $x^5_k=y^5_j=1$, and we deduce that $f_{r}(x^5)=f_{r}(y^5)=0$. Consequently, $f^u(f^{r}(x^5))=f^u(f^{r}(y^5))=\ZERO$ thus 
\[
u,v(y^1),w(y^2,j),r,w(x^4,k),r,u.
\]
sends both $x$ and $y$ on $\ZERO$ and is of length at most $5n-3$. 

\medskip
Thus, for every configurations $x,y$ on $V$, there is a word $w$ of length at most $5n-3$ such that $f^{w}(x)=f^{w}(y)$. 

\medskip
Suppose now that $r$ has a loop. Let $z=f^u(e_{r})$, let $w$ be a word of length at most $5n-3$ such that $f^{w}(\ZERO)=f^{w}(z)=\ZERO$, and let us prove that $u,r,u,w$, which is of length at most $7n-2$, sends all the configurations on $\ZERO$, and thus synchronizes $f$. Let $x$ be any configuration on $V$, and let
\[
x^1=f^u(x),\qquad x^2=f^{r}(x^1),\quad x^3=f^u(x^2),\quad x^4=f^w(x^3). 
\]
We have to prove that $x^4=\ZERO$ and by the choice of $w$ it is sufficient to prove that $x^3$ is $\ZERO$ or $z$. Suppose first that $x_{r}=1$. Then $x^1_{r}=x_{r}=1$ since $r$ is not in $u$. Since $r$ has a negative loop we have $x^2_{r}=0$  and thus $x^3=\ZERO$. Suppose now that $x_{r}=0$. Then $x^1=\ZERO$ and thus $x^2=e_{r}$ since $r$ has only negative in-neighbors. So $x^3=z$ as desired.
\end{proof}

Observing that if $G$ satisfies the conditions of Theorem~\ref{thm:fast} then, up to a switch, it satisfies the conditions of the previous lemma, we deduce the following. 

\begin{lemma}\label{lem:fast-3}
Suppose that $G$ is strong, has no positive cycles, and is not a cycle. There is an and-or-net $f$ on $G$ such that the following holds. For every configurations $x,y$ on $V$, there is a word $w$ such that $f^w(x)=f^w(y)$ and $|w|\leq 5n$. Furthermore, if $G$ has a loop, then $f$ has a synchronizing word of length at most $7n$. 
\end{lemma}

\begin{proof}
Since $G$ is strong and is not a cycle, it has a vertex of in-degree at least two and any vertex with a loop is of in-degree at least two. Let $r$ be a vertex with a loop, if it exists, and let $r$ be any vertex of in-degree at least two otherwise.  Since $G$ strong, it has a spanning out-tree $T$ rooted in $r$. By Proposition \ref{pro:harary}, there is $I\subseteq V$ such that $T^I$, the $I$-switch of $T$, is full-positive. Hence $T^I$ is a spanning full-positive out-tree of $G^I$, the $I$-switch of $G$, rooted in $r$. Since $G$ has no positive cycles, $G^I$ has no positive cycles, thus all the in-neighbors of $r$ are negative. Hence $G^I$ satisfies the conditions of Lemma~\ref{lem:fast-2}. So, denoting $h$ the and-net on $G^I$, we deduce that for every configurations $x,y$ on $V$, there is a word $w$ such that $h^w(x)=h^w(y)$ and $|w|\leq 5n$. Furthermore, if $G$ has a loop then $r$ has a loop and we deduce that $h$ has a synchronizing word of length at most $7n$. Let $f$ be the $I$-switch of $h$. By Proposition~\ref{pro:BN_switch}, $f$ is an and-or-net on $G$ with the desired properties. 
\end{proof}

From the previous lemma and Lemma \ref{lem:2-synch}, we deduce that if $G$ satisfies the conditions of Theorem~\ref{thm:fast}, then there is an and-or-net on $G$ with a synchronizing word of length at most $5n(2^n-1)$. To obtain the bound $5n(\sqrt{2})^n$, we need additional arguments. We first need the following easy observation, which implies a basic result of \cite{AGRS20}: if $G$ is acyclic, then there is a word $w$ of length $n$ which synchronizes every BN on $G$. 

\begin{lemma}\label{lem:acyclic_image}
Let $I\subseteq V$ such that $G\setminus I$ is acyclic, and let $w$ be a topological sort of $G\setminus I$. Let $f$ be a BN on $G$. For every configurations $x,y$ on $V$ with $x_I=y_I$ we have $f^w(x)=f^w(y)$. 
\end{lemma}

\begin{proof}
Let $x,y$ be configurations on $V$ with $x_I=y_I$. Let us write $w=i_1,\dots,i_\ell$, where $\ell=|V\setminus I|$. Let $I_0=I$ and, for $1\leq k\leq \ell$, let $I_k=I_{k-1}\cup\{i_k\}$. We will prove, by induction on $k$, that 
\[\tag{1}
f^w(x)_{I_k}=f^w(y)_{I_k}.
\]
Since $I\cap\{w\}=\emptyset$, we have $f^w(x)_I=x_I=y_I=f^w(y)_I$, thus (1) holds for $k=0$. Suppose that $1\leq k\leq \ell$. Let $u,v$ such that $w=u,i_k,v$, and let 
\[
\begin{array}{l}
x^1=f^u(x),\\
y^1=f^u(y),
\end{array}
\quad
\begin{array}{l}
x^2=f^{i_k}(x^1),\\
y^2=f^{i_k}(y^1),
\end{array}
\quad
\begin{array}{l}
x^3=f^v(x^2),\\
y^3=f^v(y^2).
\end{array}
\]
By induction, $x^3_{I_{k-1}}=x^3_{I_{k-1}}$. Since $I_{k-1}\cap \{i_k,v\}=\emptyset$, this implies $x^1_{I_{k-1}}=y^1_{I_{k-1}}$. Since all the in-neighbors of $i_k$ are in $I_{k-1}$, we deduce that $x^2_{i_k}=x^2_{i_k}$ and since $i_k\not\in\{v\}$ this implies $x^3_{i_k}=x^3_{i_k}$. Thus $x^3_{I_k}=y^3_{I_k}$, completing the induction step. So (1) holds for all $0\leq k\leq\ell$. Since $I_\ell=V$, we obtain $f^w(x)=f^w(y)$. 
\end{proof}

Next, we need the following recent result, showing that signed digraphs without positive cycles are rather well structured. 

\begin{theorem}[Millani, Steiner and Wiederrecht \cite{MSW19}]\label{thm:chromatic}
Suppose that $G$ has no positive cycles and no loops. There are disjoint subsets $I,J\subseteq V$ with $I\cup J=V$ such that $G[I]$ and $G[J]$ are acyclic.  
\end{theorem}

We are now in position to prove Theorem~\ref{thm:fast}, that we restate.

\setcounter{theorem}{1}
\begin{theorem}
Let $G$ be a strongly connected signed digraph on $[n]$ without positive cycles, which is not a cycle. At least one and-or-net on $G$ has a synchronizing word of length at most $5n(\sqrt{2})^n$. \end{theorem}
\setcounter{theorem}{5}

\begin{proof}
If $G$ has a loop then, by Lemma~\ref{lem:fast-3}, there is an and-or-net on $G$ with a synchronizing word of length at most $7n\leq 5n(\sqrt{2})^n$. So suppose that $G$ has no loops. By Theorem~\ref{thm:chromatic}, there is a subset $I\subseteq V$ of size at most $n/2$ such that $G\setminus I$ is acyclic. By Lemma~\ref{lem:fast-3}, there is an and-or-net $f$ on $G$ such that, for every configurations $x,y$ on $V$, there is a word $v$ such that $f^v(x)=f^v(y)$ and $|v|\leq 5n$. Let $u$ be a topological sort of $G\setminus I$, and let $X=f^u(\B^V)$. If $|X|>2^{|I|}$ then there are $x,y\in X$ with $x\neq y$ and $x_I=y_I$. Since $x,y\in X$, there are configurations $x',y'$ on $V$ with $x=f^u(x')$ and $y=f^u(y')$. Since $I\cap\{u\}=\emptyset$, we have $x'_I=x_I=y_I=y'_I$. But then, by Lemma~\ref{lem:acyclic_image}, we have $f^u(x')=f^u(y')$ and thus $x=y$, a contradiction. Hence $|X|\leq 2^{|I|}\leq 2^{n/2}$. By Lemma~\ref{lem:2-synch}, there is a word $w$ such that $|f^w(X)|=1$ and $|w|\leq 5n(|X|-1)$. Thus $u,w$ is a synchronizing word for $f$ of length at most 
\[
n-|I|+5n(|X|-1)\leq n+5n(2^{n/2}-1)= 5n(\sqrt{2})^n-4n.
\] 
\end{proof}

%%%%%%%%%%%%%%%%%%%%%%%%%%%%%%%%%%%%%%%%%%%%%%%%%%%%%%%%%%%%%%%%%%%%%%%%%%%%%%%%%%%%%%%%%%%%%%%%%%%%%%%
\section{Proof of Theorem~\ref{thm:complexity}}\label{sec:complexity}
%%%%%%%%%%%%%%%%%%%%%%%%%%%%%%%%%%%%%%%%%%%%%%%%%%%%%%%%%%%%%%%%%%%%%%%%%%%%%%%%%%%%%%%%%%%%%%%%%%%%%%%

Recall that if $G$ is a simple signed digraph with maximum in-degree at most two, then every BN on $G$ is an and-or-net. Hence, Theorem~\ref{thm:complexity} exactly says that the following two decisions problems are coNP-hard. Note that, since Robertson, Seymour and Thomas \cite{RST99} and McCuaig \cite{McC04} proved independently that we can decide in polynomial time if a signed digraph has a positive cycle, the conditions on the inputs of these two problems can be checked in polynomial time.

\bigskip
\begin{minipage}{13.5cm}
{\sc Strong-Synchronizing-Problem}
\begin{itemize}
\item
{\sc Input:} A simple signed digraph $G$, {\em strongly connected}, with maximum in-degree at most two, and containing a vertex meeting every cycle.
\item
{\sc Question:} Is every BN on $G$ synchronizing? 
\end{itemize}
\end{minipage}

\bigskip
\begin{minipage}{13.5cm}
{\sc Negative-Synchronizing-Problem}
\begin{itemize}
\item
{\sc Input:} A simple signed digraph $G$, {\em without positive cycles}, with maximum in-degree at most two, and containing a vertex meeting every cycle.
\item
{\sc Question:} Is every BN on $G$ synchronizing?
\end{itemize}
\end{minipage}

\bigskip
We prove that these two problems are coNP-hard with reductions from 3-SAT. So consider a 3-CNF formula $\psi$ over a set of $n\geq 2$ variables $\LAMBDA=\{\lambda_1,\dots,\lambda_n\}$ consisting of $m\geq 1$ clauses $\MU=\{\mu_1,\dots,\mu_m\}$. To each variable $\lambda_r$ is associated a positive literal $\lambda^+_r$ and a negative literal $\lambda^-_r$. The resulting sets of positive and negative literals are denoted $\LAMBDA^+$ and $\LAMBDA^-$. Each clause $\mu_s$ is a subset of $\LAMBDA^+\cup\LAMBDA^-$ of size three, and we write $\mu_r=\{\mu_{r,1},\mu_{r,2},\mu_{r,3}\}$. An \EM{assignment} for $\psi$ is regarded as a configuration $z$ on $[n]$. A positive literal $\lambda^+_r$ is satisfied by $z$ if $z_r=1$, and a negative literal $\lambda^-_r$ is satisfied by $z$ if $z_r=0$. A clause is satisfied by $z$ if at least one of its literals is satisfied by $z$. The formula $\psi$ is \EM{satisfied} by $z$ (or $z$ is a \EM{satisfying assignment} for $\psi$) if every clause in $\MU$ is satisfied by $z$. We say that $\psi$ is \EM{satisfiable} if it has at least one satisfying assignment. 

\medskip
Our reductions from 3-SAT are based on the simple signed digraph $H_\psi$ defined as follows; see Figure~\ref{fig:H_psi} for an illustration:
\begin{itemize}
\item
The vertex set is $\LAMBDA^+\cup\LAMBDA^-\cup\ELL\cup \MU'\cup\MU\cup \C$, 
where $\ELL=\{\ell_0,\dots,\ell_n\}$, $\MU'=\{\mu'_1,\dots,\mu'_m\}$ and $\C=\{c_1,\dots,c_m\}$; there are thus $3n+3m+1$ vertices. 
\item The arcs are,	 for all $r\in [n]$ and $s\in [m]$, 
\begin{itemize}
\item
$(\ell_{r-1},\lambda^+_r),(\ell_{r-1},\lambda^-_r)$, 
\item
$(\lambda^+_r,\ell_r),(\lambda^-_r,\ell_r)$,
\item
$(c_1,\ell_0)$,
\item
$(\mu_{s,1},\mu'_s),(\mu_{s,2},\mu'_s)$,$(\mu_{s,3},\mu_s),(\mu'_s,\mu_s)$,
\item
$(\mu_s,c_s),(c_{s+1},c_s)$, where $c_{m+1}$ means $\ell_n$.
\end{itemize} 
\item
For all $s\in [m]$, the arc $(\mu_s,c_s)$ is negative, and all the other arcs are positive. 
\end{itemize}

\begin{figure}
\[
\begin{tikzpicture}
\node (Lp1) at (1.5,1){$\lambda^+_1$};
\node (Lm1) at (1.5,-0.5){$\lambda^-_1$};
\node (Lp2) at (4.5,1){$\lambda^+_2$};
\node (Lm2) at (4.5,-0.5){$\lambda^-_2$};
\node (Lp3) at (7.5,1){$\lambda^+_3$};
\node (Lm3) at (7.5,-0.5){$\lambda^-_3$};
\node (1) at (0,0.25){$\ell_0$};
\node (2) at (3,0.25){$\ell_1$};
\node (3) at (6,0.25){$\ell_2$};
\node (4) at (9,0.25){$\ell_3$};
\node (mu1') at (3,-3){$\mu'_1$};
\node (mu2') at (6,-3){$\mu'_2$};
\node (mu1) at (3,-4.5){$\mu_1$};
\node (mu2) at (6,-4.5){$\mu_2$};
\node (c1) at (3,-6){$c_1$};
\node (c2) at (6,-6){$c_2$};
\path[Green,->,very thick]
(1) edge (Lp1)
(1) edge (Lm1)
(Lp1) edge (2)
(Lm1) edge (2)
(2) edge (Lp2)
(2) edge (Lm2)
(Lp2) edge (3)
(Lm2) edge (3)
(3) edge (Lp3)
(3) edge (Lm3)
(Lp3) edge (4)
(Lm3) edge (4)
(4) edge[bend left=30] (c2)
(c2) edge (c1)
(c1) edge[bend left=30] (1)
;
\path[Green,->,very thick]
(Lp1) edge[densely dashed,bend left=10] (mu1')
(Lm2) edge[densely dashed] (mu1')
(Lm3) edge[densely dashed,bend right=10] (mu1)
(mu1') edge (mu1)
(mu1) edge[red,-|] (c1);
\path[Green,->,very thick]
(Lm1) edge[densely dashed] (mu2')
(Lp2) edge[densely dashed,bend left=10] (mu2')
(Lm3) edge[densely dashed,bend left=10] (mu2)
(mu2') edge (mu2)
(mu2) edge[red,-|] (c2)
;
\end{tikzpicture}
\]
\caption{\label{fig:H_psi}
The signed digraph $H_\psi$ for the 3-CNF formula $\psi=(\lambda_1\lor \neg\lambda_2\lor\neg\lambda_3)\land(\neg\lambda_1\lor \lambda_2\lor \neg\lambda_3)$. Using our notations, the set of variables is $\LAMBDA=\{\lambda_1,\lambda_2,\lambda_3\}$ and the set of clauses is $\MU=\{\mu_1,\mu_2\}$ with $\mu_1=\{\lambda^+_1,\lambda^-_2,\lambda^-_3\}$ and $\mu_2=\{\lambda^-_1,\lambda^+_2,\lambda^-_3\}$. Clauses are encoded through dashed arrows. 
}
\end{figure}
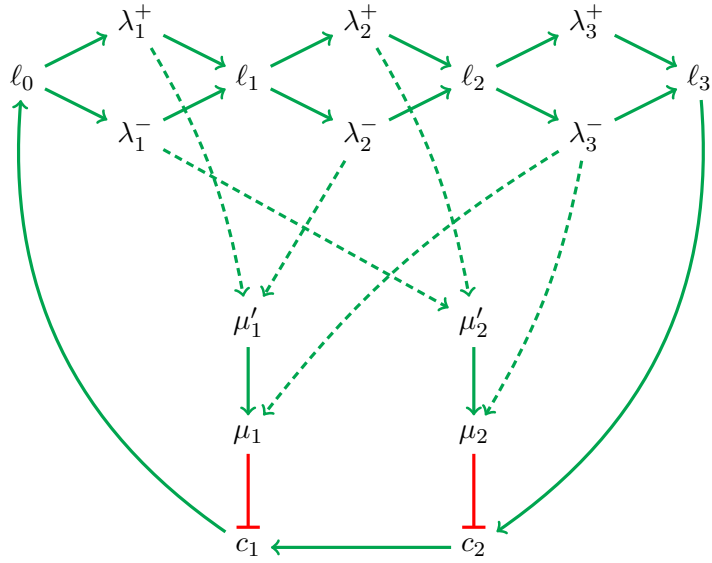

\begin{figure}
\[
\begin{tikzpicture}
\node (t) at (4.5,4){$t$};
\node (Lp1) at (1.5,1){$\lambda^+_1$};
\node (Lm1) at (1.5,-0.5){$\lambda^-_1$};
\node (Lp2) at (4.5,1){$\lambda^+_2$};
\node (Lm2) at (4.5,-0.5){$\lambda^-_2$};
\node (Lp3) at (7.5,1){$\lambda^+_3$};
\node (Lm3) at (7.5,-0.5){$\lambda^-_3$};
\node (1) at (0,0.25){$\ell_0$};
\node (2) at (3,0.25){$\ell_1$};
\node (3) at (6,0.25){$\ell_2$};
\node (4) at (9,0.25){$\ell_3$};
\node (mu1') at (3,-3){$\mu'_1$};
\node (mu2') at (6,-3){$\mu'_2$};
\node (mu1) at (3,-4.5){$\mu_1$};
\node (mu2) at (6,-4.5){$\mu_2$};
\node (c1) at (3,-6){$c_1$};
\node (c2) at (6,-6){$c_2$};
\node (q) at (1.65,2.7){$q$};
\path[Green,->,very thick]
(1) edge (Lp1)
(1) edge (Lm1)
(Lp1) edge (2)
(Lm1) edge (2)
(2) edge (Lp2)
(2) edge (Lm2)
(Lp2) edge (3)
(Lm2) edge (3)
(3) edge (Lp3)
(3) edge (Lm3)
(Lp3) edge (4)
(Lm3) edge (4)
(4) edge[bend left=30] (c2)
(c2) edge (c1)
(c1) edge[bend left=30] (1)
;
\path[Green,->,very thick]
(Lp1) edge[densely dashed,bend left=10] (mu1')
(Lm2) edge[densely dashed] (mu1')
(Lm3) edge[densely dashed,bend right=10] (mu1)
(mu1') edge (mu1)
(mu1) edge[red,-|] (c1);
\path[Green,->,very thick]
(Lm1) edge[densely dashed] (mu2')
(Lp2) edge[densely dashed,bend left=10] (mu2')
(Lm3) edge[densely dashed,bend left=10] (mu2)
(mu2') edge (mu2)
(mu2) edge[red,-|] (c2)
;
\path[Green,->,thick]
(Lm1) edge[densely dashed] (mu2')
(Lp2) edge[densely dashed,bend left=10] (mu2')
(Lm3) edge[densely dashed,bend left=10] (mu2)
(mu2') edge (mu2)
(mu2) edge[red,-|] (c2)
;
\path[Green,->,thick]
(t) edge (Lp1)
(t) edge (Lm1)
(t) edge (Lp2)
(t) edge[bend left=25] (Lm2)
(t) edge (Lp3)
(t) edge (Lm3)
(1) edge[bend left=55] (t)
(1) edge[red,bend left=10,-|] (q)
(q) edge[bend left=10] (t)
;
\end{tikzpicture}
\]
\caption{\label{fig:G_psi}
The signed digraph $G_\psi$ for the 3-CNF formula of Figure~\ref{fig:H_psi}. 
}
\end{figure}
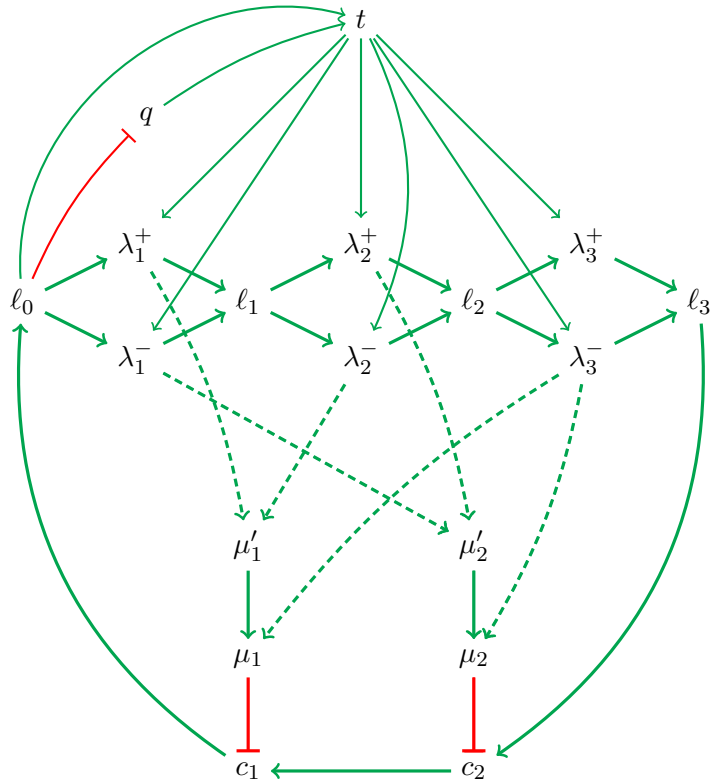

This signed digraph $H_\psi$ has been recently introduced in \cite{BDPR22}. Actually, in that paper, the authors consider the signed digraph $H'_\psi$ obtained from $H_\psi$ by adding, for all $r\in [n]$, a new vertex $i_r$ with a positive arc $(i_r,\lambda^+_r)$ and a negative arc $(i_r,\lambda^-_r)$, and then prove the following: $\psi$ is satisfiable if and only if there is a BN on $H'_\psi$ with at least two fixed points.

\medskip
Here, we present two adaptations of this construction to the context of synchronization. For the first, let $G_\psi$ be the signed digraph obtained from $H_\psi$ by adding: two new vertices, $q$ and $t$; a negative arc $(\ell_0,q)$; two positive arcs $(\ell_0,t),(q,t)$; and a positive arc $(t,i)$ for all $i\in\LAMBDA^+\cup\LAMBDA^-$. An illustration is given in Figure \ref{fig:G_psi}. Note that $G_\psi$ is strongly connected, simple and has maximum in-degree two (thus every BN on $G_\psi$ is an and-or-network). Also, $\ell_0$ meets every cycle. Thus $G_\psi$ is a valide instance of the {\sc Strong-Synchronizing-Problem}. 

\medskip
The following result shows that $G_\psi$ as the same property that $H'_\psi$: $\psi$ is satisfiable if and only if there is a BN on $H'_\psi$ with at least two fixed points. But it gives a stronger conclusion when $\psi$ is not satisfiable: every BN on $G$ is synchronizing. Thus $\psi$ is not satisfiable if and only if every BN on $G_\psi$ is synchronizing. Since the SAT problem is NP-complete and $G_\psi$ is a valid instance of the {\sc Strong-Synchronizing-Problem}, we deduce that this problem is coNP-hard.   

\begin{theorem}\label{thm:fixed_point}
The following conditions are equivalent:
\begin{enumerate}
\item $\psi$ is not satisfiable. 
\item Every BN on $G_\psi$ has at most one fixed point.
\item There is a word of length $12n+12m+11$ which synchronizes every BN on~$G_\psi$.
\item Every BN on $G_\psi$ is synchronizing.
\end{enumerate}
\end{theorem}

\begin{proof}
\BF{($1\Rightarrow 2$)} Suppose that $\psi$ is not satisfiable and let $f$ be a BN on $G_\psi$. We will prove that $f$ has at most one fixed point. Suppose, for a contradiction, that $f$ has two distinct fixed points $x$ and $y$. As argued in \cite{A08}, $G_\psi$ has a positive cycle $C$ such that $x_i\neq y_i$ for every vertex $i$ in $C$. Since  $x_q=f_q(x)=\neg x_{\ell_0}$ and  $y_q=f_q(y)=\neg y_{\ell_0}$, we have $f_t(x)=f_t(y)=0$ if $f_t$ is a conjunction and $f_t(x)=f_t(y)=1$ if $f_t$ is a disjunction, and thus $x_t=y_t$. So $t$ is not in $C$. We deduce that
\begin{itemize}
\item
$C$ is full-positive,
\item
$C$ contains all the vertices in $\ELL\cup\C$,
\item
$C$ contains exactly one of $\lambda^+_r,\lambda^-_r$ for each $r\in [n]$.
\end{itemize}
Suppose, without loss of generality, that $x_{\ell_0}<y_{\ell_0}$. We prove, by induction on $s$ from $1$ to $m$ that $x_{c_s}<y_{c_s}$. Since $c_1$ is the unique in-neighbor of $\ell_0$ we have $x_{c_1}<y_{c_1}$. Let $1<s\leq m$. By induction, $x_{c_{s-1}}<y_{c_{s-1}}$. Suppose, for a contradiction, that $x_{c_s}>y_{c_s}$. If $f_{c_{s-1}}$ is a disjunction then, since $x_{c_s}=1$, we have $f_{c_{s-1}}(x)=1\neq x_{c_{s-1}}$, a contradiction, and if $f_{c_{s-1}}$ is a conjunction then, since $y_{c_s}=0$, we have $f_{c_{s-1}}(y)=0\neq y_{c_{s-1}}$, a contradiction. We deduce that $x_{c_s}\leq y_{c_s}$ and since the equality is not possible we obtain $x_{c_s}<y_{c_s}$. This completes the induction, and thus $x_{\C}=\ZERO$ and $y_{\C}=\ONE$. Since $x_{\ell_0}\leq y_{\ell_0}$ and $x_t=y_t$, we have $x_i\leq y_i$ for $i\in\{\lambda^+_1,\lambda^-_1\}$, and thus $x_{\ell_1}\leq y_{\ell_1}$. Repeating this argument we get $x_i\leq y_i$ for every $i\in \ELL\cup\LAMBDA^+\cup\LAMBDA^-$. Consequently, $x_i\leq y_i$ for all $i\in \MU'$ and thus $x_i\leq y_i$ for all $i\in \MU$. We deduce that $x_{\mu_s}=y_{\mu_s}$ for all $s\in [m]$. Indeed, otherwise $x_{\mu_s}<y_{\mu_s}$ so $f_{c_s}(x)=1\neq x_{c_s}$ if $f_i$ is a disjunction and $f_{c_s}(y)=0\neq y_{c_s}$ if $f_i$ is a conjunction. Thus $x_{\mu_s}=y_{\mu_s}$. Let us prove that $\mu_s$ contains a literal $i$ with $x_i=y_i$. Indeed, otherwise we have $x_i<y_i$ for the three literals $i$ contained in $\mu_s$, which implies $x_{\mu'_s}<y_{\mu'_s}$ and thus $x_{\mu_s}<y_{\mu_s}$, a contradiction. Thus each clause $\mu_s$ contains a literal which is not in $C$ (since $x_i\neq y_i$ for all vertex $i$ in $C$). Let $z\in\B^n$ defined by $z_r=1$ if $\lambda^+_r$ is not in $C$ and $z_r=0$ if $\lambda^-_r$ is not in $C$; there is no ambiguity since exactly one of $\lambda^+_r,\lambda^-_r$ is in $C$. All the literals not in $C$ are satisfied by $z$, and thus each clause contains a literal satisfied by $z$. So $\psi$ is satisfiable, a contradiction. Thus $f$ has indeed at most one fixed point. 

\bigskip
\noindent
\BF{($2\Rightarrow 3$)} Suppose that every BN on $G_\psi$ has at most one fixed point. Let $u$ be a topological sort of $G_\psi\setminus\{\ell_0,q,t\}$; so $|u|=3n+3m$. The first two letters of $u$ are $\lambda^+_1,\lambda^-_1$ and we denote by $u'$ the suffix of $u$ obtained by removing these two letters; so $|u'|=3n+3m-2$. Let $V$ be the vertex set of $G_\psi$; so $|V|=3n+3m+3$. We will prove that the word $\W$ defined by  
\[
\W=v,\ell_0,v,w\quad\textrm{where}\quad v=q,t,u\quad\textrm{and}\quad  w=\ell_0,u,t,\ell_0,\lambda^+_1,\lambda^-_1,t,u',\ell_0,q
\]
is a synchronizing word for every BN $f$ on $G_\psi$. Since $|v|=|V|-1$ and $|w|=2|V|$, we have $|\W|=4|V|-1=12n+12m+11$ as desired.

\medskip
Let $f$ be a BN on $G_\psi$. For $a=0,1$, let $f^a$ be the BN with component set $V$ such that $f^a_{\ell_0}$ is the $a$-constant function and $f^a_i=f_i$ for every $i\neq\ell_0$. Note that, for every word $w$ not containing $\ell_0$ and any configuration $x$ on $V$, we have $f^w(x)=(f^a)^w(x)$. Let $H$ be obtained from $G_\psi$ by deleting $(c_1,\ell_0)$.  Then $H$ is the signed interaction digraph of $f^a$, and it is acyclic. Hence $f^a$ has a unique fixed point; let $x$ be the fixed point of $f^0$ and $y$ the fixed point of $f^1$. Obviously, $x_{\ell_0}<y_{\ell_0}$. It is clear that if $f(x)_{\ell_0}<f(y)_{\ell_0}$ then $x$ and $y$ are fixed points of $f$, a contradiction. Thus $f(x)_{\ell_0}\geq f(y)_{\ell_0}$. 

\bigskip
\noindent
(1) {\em We have $x_q>y_q$ and,  for every configuration $z$ on $V$,
\begin{itemize}
\item $f^v(z)=x$ if $z_{\ell_0}=0$,
\item $f^v(z)=y$ if $z_{\ell_0}=1$,
\item $f^u(y+e_{\ell_0})=x+e_q$,
\item $f^u(x+e_{\ell_0})=y+e_q$.
\end{itemize}}

\medskip
\noindent
We have $x_q=f_q(x)=\neg x_{\ell_0}=1$ and $y_q=f_q(y)=\neg y_{\ell_0}=0$. We deduce that $f_t(x)=f_t(y)=1$ if $f_t$ is a disjunction and $f_t(x)=f_t(y)=0$ if $f_t$ is a conjunction. So $x_t=y_t$ in all cases. Let $z$ be a configuration on $V$. Since $\ell_0,v$ is a topological sort of $H$ and $x$ is a fixed point of $f^0$, we deduce from Lemma~\ref{lem:acyclic_image} that $(f^0)^{\ell^0,v}(z)=x$. So if $z_{\ell_0}=0$ then $f^v(z)=(f^0)^v(z)=(f^0)^{\ell_0,v}(z)=x$. This proves the first item. Let $z=y+e_{\ell_0}$. Since $z_{\ell_0}=0$ we have $f^{q,t,u}(z)=f^v(z)=x$ by the first item and since $z_t=y_t=x_t$ we deduce that $f^{q,u}(z)=x$. Since $G_\psi$ has no arc from $q$ to a vertex in $\{u\}$ we have $f^u(z)_i=f^{q,u}(z)_i=x_i$ for all $i\neq q$. Since $f^u(z)_q=z_q=y_q\neq x_q$ we deduce that $f^u(z)=x+e_q$. This proves the third item. We prove similarly the second and fourth items, and (1) follows.

\bigskip
\noindent
(2) {\em If $f^{\ell_0,v,w}(x)=f^{\ell_0,v,w}(y)$ then $\W$ is a synchronizing word for $f$.}

\medskip
\noindent
Let $z,z'$ be any configurations on $V$. If $z_{\ell_0}=z'_{\ell_0}$ then $f^v(z)=f^v(z')$ by the first two items of (1) and thus $f^\W(z)=f^\W(z')$. Suppose, without loss of generality, that $z_{\ell_0}< z'_{\ell_0}$. Then, by the first two items of (1) we have $f^v(z)=x$ and $f^v(z')=y$, and thus if $f^{\ell_0,v,w}(x)=f^{\ell_0,v,w}(y)$ we have $f^\W(z)=f^\W(z')$. This proves (2).

\bigskip
\noindent
(3) {\em If $f(x)_{\ell_0}=f(y)_{\ell_0}$ then $\W$ is a synchronizing word for $f$.}

\medskip
\noindent
Suppose that $f(x)_{\ell_0}=f(y)_{\ell_0}$, that is, $f^{\ell_0}(x)_{\ell_0}=f^{\ell_0}(y)_{\ell_0}$. Then, by the first two items of (1) we have $f^{\ell_0,v}(x)=f^{\ell_0,v}(y)$, thus $f^{\ell_0,v,w}(x)=f^{\ell_0,v,w}(y)$ and by (2) $\W$ is a synchronizing word for $f$. This proves (3). 

\bigskip
By (3) we can suppose that $f(x)_{\ell_0}\neq f(y)_{\ell_0}$ and thus $f(x)_{\ell_0}> f(y)_{\ell_0}$. Consequently, by (1) we have $f^{\ell_0,v}(x)=y$ and $f^{\ell_0,v}(y)=x$. Thus, by (2), it is sufficient to prove that $f^w(x)=f^w(y)$. Since $f(x)_{\ell_0}> f(y)_{\ell_0}$ we have $f^{\ell_0}(x)=x+e_{\ell_0}$ and $f^{\ell_0}(y)=y+e_{\ell_0}$ and from the last two items of (1) we have $f^u(x+e_{\ell_0})=y+e_q$ and $f^u(y+e_{\ell_0})=x+e_q$. Since $(y+e_q)_{\ell_0}=(y+e_q)_q=1$ we have $f^t(y+e_q)_t=1$ and since $(x+e_q)_{\ell_0}=(x+e_q)_q=0$ we have $f^t(y+e_q)_t=0$. Since $G_\psi$ has no arc from $\{q,t\}$ to $\ell_0$ we have $f^{t,\ell_0}(y+e_q)_{\ell_0}=f(y)_{\ell_0}=0$ and $f^{t,\ell_0}(x+e_q)_{\ell_0}=f(x)_{\ell_0}=1$. Thus, summing up, we have 
\begin{itemize}
\item 
$f^{\ell_0,u}(x)=y+e_q$ and $f^{t,\ell_0}(y+e_q)_{\ell_0}< f^{t,\ell_0}(y+e_q)_t$,
\item
$f^{\ell_0,u}(y)=x+e_q$ and $f^{t,\ell_0}(x+e_q)_{\ell_0}> f^{t,\ell_0}(x+e_q)_t$.
\end{itemize} 
Setting $y'=f^{t,\ell_0,\lambda^+_1,\lambda^-_1}(y+e_q)$ and $x'=f^{t,\ell_0,\lambda^+_1,\lambda^-_1}(x+e_q)$ we deduce that 
\begin{itemize}
\item
$y'_{\lambda^+_1}=x'_{\lambda^+_1}=0$ if $f_{\lambda^+_1}$ is a conjunction,
\item
$y'_{\lambda^+_1}=x'_{\lambda^+_1}=1$ if $f_{\lambda^+_1}$ is a disjunction,
\item
$y'_{\lambda^-_1}=x'_{\lambda^-_1}=0$ if $f_{\lambda^-_1}$ is a conjunction,
\item
$y'_{\lambda^-_1}=x'_{\lambda^-_1}=1$ if $f_{\lambda^-_1}$ is a disjunction.
\end{itemize} 
Furthermore, since 
\begin{itemize}
\item
$y'_{\ell_0}=f^{t,\ell_0}(y+e_q)_{\ell_0}=0$ and $y'_q=(y+e_q)_q=1$,
\item
$x'_{\ell_0}=f^{t,\ell_0}(x+e_q)_{\ell_0}=1$ and $x'_q=(x+e_q)_q=0$, 
\end{itemize}
we have 
\begin{itemize}
\item
$f^t(y')_t=f^t(x')_t=0$ if $f_t$ is a conjunction,
\item
$f^t(y')_t=f^t(x')_t=1$ if $f_t$ is a disjunction.
\end{itemize}
Hence, setting $I=\{t,\lambda^+_1,\lambda^-_1\}$, we have $f^t(y')_I=f^t(x')_I$. Since $G_\psi\setminus I$ is acyclic and $u',\ell_0,q$ is a corresponding topological sort, by Lemma~\ref{lem:acyclic_image} we have $f^{t,u',\ell_0,q}(y')=f^{t,u',\ell_0,q}(x')$, and thus $f^w(x)=f^w(y)$ as desired. 

\bigskip
\noindent
\BF{($3 \Rightarrow 4$)} This is obvious.

\bigskip
\noindent
\BF{($4 \Rightarrow 1$)} Suppose that $z$ is a satisfying assignment of $\psi$, and let us prove that some BN on $G_\psi$ is not synchronizing. Let $\LAMBDA^1$ be the literals satisfied by $z$ and $\LAMBDA^0$ the literals not satisfied by $z$, hence $(\LAMBDA^0,\LAMBDA^1)$ is a balanced partition of $\LAMBDA^+\cup\LAMBDA^-$. Let $f$ be the and-or-net on $G_\psi$ such that:
\begin{itemize}
\item
for all $i\in \LAMBDA^1\cup\MU\cup\MU'\cup\C\cup\{t\}$, $f_i$ is a disjunction;
\item
for all $i\in\LAMBDA^0\cup\ELL$, $f_i$ is a conjunction. 
\end{itemize}

\medskip
Let $\MU'^1$ be the vertices in $\MU'$ with at least one in-neighbor in $\LAMBDA^1$ and $\MU'^0=\MU'\setminus\MU'^1$. Let $I=\LAMBDA^1\cup\MU'^1\cup \MU$ and $J=\LAMBDA^0\cup\ELL\cup\C\cup\MU'^0$. Let us prove that the configurations $x,y$ on $V$ defined as follows are fixed points of $f$:
\begin{itemize}
\item
$x_q=1$, $x_t=1$, $x_I=\ONE$, $x_J=\ZERO$,
\item
$y_q=0$, $y_t=1$, $y_I=\ONE$, $y_J=\ONE$.
\end{itemize}

\bigskip
\noindent
(4) {\em $f(x)_i=x_i$ and $f(y)_i=y_i$ for all $i\in \{q,t\}\cup I$.}

\medskip
\noindent
Since $x_q\neq x_{\ell_0}$ and $y_q\neq y_{\ell_0}$, we have $f(x)_q=\neg x_{\ell_0}=x_q$ and $f(y)_q=\neg y_{\ell_0}=y_q$, and since $f_t$ is a disjunction, we also deduce that $f_t(x)=f_t(y)=1=x_t=y_t$. We now prove that $f(x)_I=f(y)_I=\ONE$. If $i\in\LAMBDA^1\cup\MU'^1$, then $f_i$ is a disjunction and $i$ has an in-neighbor  $j\in\LAMBDA^1\cup\{t\}$, so $x_j=y_j=1$ and we deduce that $f_i(x)=f_i(y)=1$. Next, consider a clause $\mu_s$ in $\MU$. Then $f_{\mu_s}$ is a disjunction and, since $z$ is a satisfying assignment, $\mu_s$ contains a literal $i\in\LAMBDA^1$ and thus $x_i=y_i=1$. If $i$ is an in-neighbor of $\mu_s$ then $f_{\mu_s}(x)=f_{\mu_s}(y)=1$. Otherwise, $i$ is an in-neighbor of $\mu'_s$, so $\mu'_s\in\MU'^1$. Thus $x_{\mu'_s}=y_{\mu'_s}=1$ and we deduce that $f_{\mu_s}(x)=f_{\mu_s}(y)=1$. This proves~(4).

\bigskip
\noindent
(5) {\em $f(x)_J=\ZERO$.}

\medskip
\noindent
If $i\in\MU'^0$ then for every in-neighbor $j$ of $i$ we have $j\in\LAMBDA^0$, thus $x_j=0$, and we deduce that $f_i(x)=0$. If $i\in\LAMBDA^0\cup \ELL$, then $f_i$ is a conjunction and has an in-neighbor $j\in\LAMBDA^0\cup\ELL\cup\C$. Thus $x_j=0$ we deduce that $f_i(x)=0$. Setting $c_{m+1}=\ell_n$, for all $s\in [m]$  we have $f_{c_s}(x)=\neg x_{\mu_s}\lor x_{c_{s+1}}=\neg 1\lor 0=0$. This proves (5). 

\bigskip
\noindent
(6) {\em $f(y)_J=\ONE$.}

\medskip
\noindent
If $i\in\LAMBDA^0\cup \ELL\cup \MU'^0$ then for every in-neighbor $j$ of $i$ we have $j\neq q$, thus $y_j=1$, and we deduce that $f_i(y)=1$. If $i\in \C$ then $i$ has a positive in-neighbor $j\in\C\cup\{\ell_n\}$, thus $x_j=1$, and since $f_i$ is a disjunction, we have $f_i(x)=1$. This proves (6). 

\bigskip
By (4), (5) and (6), $x$ and $y$ are fixed points of $f$ and thus $f$ is not synchronizing. 
\end{proof}

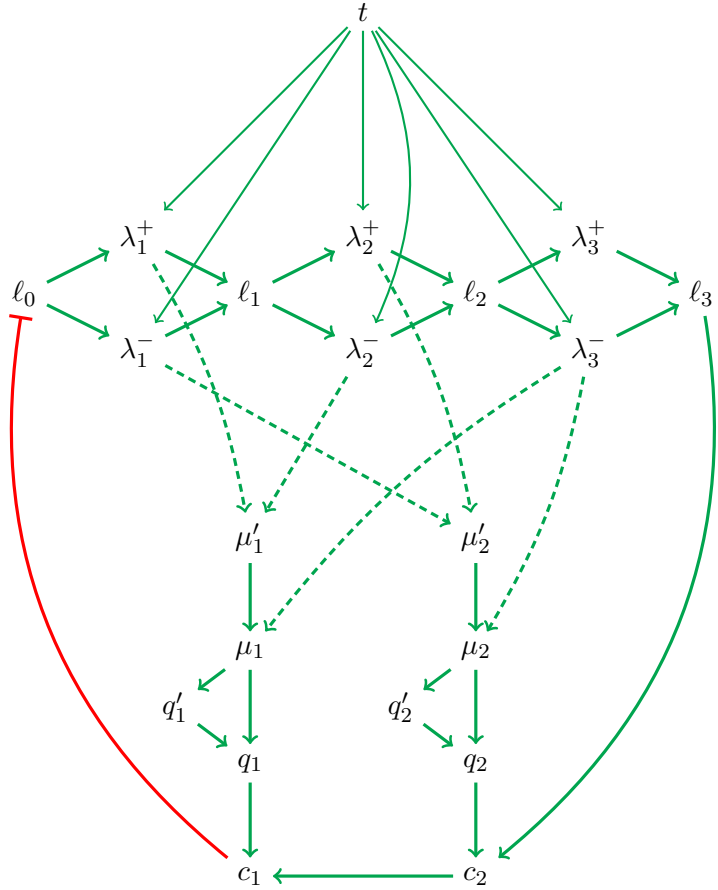
\begin{figure}
\[
\begin{tikzpicture}
\node (t) at (4.5,4){$t$};
\node (Lp1) at (1.5,1){$\lambda^+_1$};
\node (Lm1) at (1.5,-0.5){$\lambda^-_1$};
\node (Lp2) at (4.5,1){$\lambda^+_2$};
\node (Lm2) at (4.5,-0.5){$\lambda^-_2$};
\node (Lp3) at (7.5,1){$\lambda^+_3$};
\node (Lm3) at (7.5,-0.5){$\lambda^-_3$};
\node (1) at (0,0.25){$\ell_0$};
\node (2) at (3,0.25){$\ell_1$};
\node (3) at (6,0.25){$\ell_2$};
\node (4) at (9,0.25){$\ell_3$};
\node (mu1') at (3,-3){$\mu'_1$};
\node (mu2') at (6,-3){$\mu'_2$};
\node (mu1) at (3,-4.5){$\mu_1$};
\node (mu2) at (6,-4.5){$\mu_2$};
\node (c1) at (3,-7.5){$c_1$};
\node (c2) at (6,-7.5){$c_2$};
\node (q1) at (3,-6){$q_1$};
\node (q2) at (6,-6){$q_2$};
\node (q'1) at (2,-5.25){$q'_1$};
\node (q'2) at (5,-5.25){$q'_2$};

\path[Green,->,very thick]
(1) edge (Lp1)
(1) edge (Lm1)
(Lp1) edge (2)
(Lm1) edge (2)
(2) edge (Lp2)
(2) edge (Lm2)
(Lp2) edge (3)
(Lm2) edge (3)
(3) edge (Lp3)
(3) edge (Lm3)
(Lp3) edge (4)
(Lm3) edge (4)
(4) edge[bend left=30] (c2)
(c2) edge (c1)
(c1) edge[red,-|,bend left=30] (1)
;
\path[Green,->,very thick]
(Lp1) edge[densely dashed,bend left=10] (mu1')
(Lm2) edge[densely dashed] (mu1')
(Lm3) edge[densely dashed,bend right=10] (mu1)
(mu1') edge (mu1)
(mu1) edge (q1)
(mu1) edge (q'1)
(q'1) edge (q1)
(q1) edge (c1)
;
\path[Green,->,very thick]
(Lm1) edge[densely dashed] (mu2')
(Lp2) edge[densely dashed,bend left=10] (mu2')
(Lm3) edge[densely dashed,bend left=10] (mu2)
(mu2') edge (mu2)
(mu2) edge (q2)
(mu2) edge (q'2)
(q'2) edge (q2)
(q2) edge (c2)
;
\path[Green,->,thick]
(t) edge (Lp1)
(t) edge (Lm1)
(t) edge (Lp2)
(t) edge[bend left=25] (Lm2)
(t) edge (Lp3)
(t) edge (Lm3)
;
\end{tikzpicture}
\]
\caption{\label{fig:F_psi}
The signed digraph $F_\psi$ for the 3-CNF formula of Figure~\ref{fig:H_psi}. 
}
\end{figure}

We now present our second construction. Let $F_\psi$ be the signed digraph obtained from $H_\psi$~by: adding a new vertex $t$ and a positive arc $(t,i)$ for all $i\in\LAMBDA^+\cup\LAMBDA^-$; making negative the arc $(c_1,\ell_0)$; deleting the arc $(\mu_s,c_s)$ for $s\in [m]$ and adding, for all $s\in [m]$, two new vertices $q_s,q'_s$ and the positive arcs $(\mu_s,q'_s),(\mu_s,q_s),(q'_s,q_s)$ and $(q_s,c_s)$. An illustration is given in Figure \ref{fig:F_psi}. Thus $F_\psi$ has a unique negative arc, $(c_1,\ell_0)$, and since all the cycles of $F_\psi$ contain this arc, all the cycles of $F_\psi$ are negative (and $\ell_0$ meets every cycle). Also $F_\psi$ is simple and has maximum in-degree two, thus every BN on $F_\psi$ is an and-or-net. 

\medskip
We prove below that $\psi$ is not satisfiable if and only if every BN on $F_\psi$ is synchronizing. Since the SAT problem is NP-complete and $F_\psi$ is a valid instance of the {\sc Negative-Synchronizing-Problem}, we deduce that this problem is coNP-hard, thus the proof of Theorem~\ref{thm:complexity} is completed. 

\begin{theorem}\label{thm:every_synch}
The following conditions are equivalent:
\begin{enumerate}
\item $\psi$ is not satisfiable.
\item Every BN on $F_\psi$ is synchronizing.
\end{enumerate}
\end{theorem}

\begin{proof}
\BF{($1 \Rightarrow 2$)} Suppose that $\psi$ is not satisfiable and suppose, for a contradiction, that some BN $f$ on $F_\psi$ is not synchronizing. Since $t$ is a source, there is $a\in\B$ such that $f_t$ is the $a$-constant function. Let $u$ be a longest canalizing word from $(t,a)$. Let $b$ be the configuration on $I=\{t\}\cup\{u\}$ such that $b_t=a$ and $b_{\{u\}}$ is the image of $u$. Hence for all configurations $x$ on $V$ we have $f^{t,u}(x)_I=b$. Let $V$ be the vertex set of $F_\psi$. If $I=V$ then $t,u$ is a synchronizing word for $f$, a contradiction. Thus $I$ is a strict subset of $V$. Let $J=V\setminus I$ and let $h$ be the BN with component set $J$ defined by $h(x_J)=f(x)_J$ for all configurations $x$ on $V$ with $x_I=b$. Let $H$ be the signed interaction digraph of $h$, which is a subgraph of $F_\psi[J]$. 

\bigskip
\noindent
(1) {\em $h$ is not synchronizing.}

\medskip
\noindent
Suppose that $h$ has a synchronizing word $v$ and let $w=t,u,v$. For all configurations $x,y$ on $V$ we have $f^{t,u}(x)_I=f^{t,u}(y)_I=b$, and since $v$ is a word over $J$ we deduce that $f^v(f^{t,u}(x))_J=h^v(f^{t,u}(x)_J)=h^v(f^{t,u}(y)_J)=f^v(f^{t,u}(y))_J$ and thus $f^w(x)=f^w(y)$. Hence $w$ is a synchronizing word for $f$, a contradiction. This proves (1). 

\bigskip
\noindent
(2) {\em $H$ has no sources.}

\medskip
\noindent
Suppose, for a contradiction, that $H$ has a source $i$. Then there is $c$ such that $h(x_J)_i=c$ for all configurations $x$ on $V$. So if $x_I=b$ then $f(x)_J=h(y_J)$ and thus $f(x)_i=h(x_J)_i=c$. For all configurations $x$ on $V$ with $x_t=a$, we have $f^u(x)_I=b$ and we deduce that $f^{u,i}(x)_i=c$. Since $i\not\in I$, $u,i$ is a canalizing word from $(t,a)$ longer than $u$, a contradiction. This proves (2). 

\bigskip
\noindent
(3) {\em $H$ has a unique initial strong component, say $S$, and is homogenous.}

\medskip
\noindent
Let $S$ be an initial strong component of $H$. Since $H$ has no source, $S$ contains a cycle, and since all the cycles of $F_\psi$ contains the arc $(c_1,\ell_0)$, this arc is in $S$ and $S$ is the unique initial strong component of $H$. Consequently, there is a path from $\ell_0$ to each vertex in $H$, and since all the paths of $F_\psi$ starting from $\ell_0$ are full-positive, we deduce that $H$ is homogenous. This proves (3).

\bigskip
\noindent
(4) {\em $S$ is a cycle.}

\medskip
\noindent
If not then $H$, by (2) and (3), $H$ has no sources and no initial cycles, and since $H$ has only negative cycles, we deduce from Theorem~\ref{thm:main2} that $h$ is synchronizing, and this contradicts (1). This proves (4). 

\bigskip
\noindent
(5) {\em If $F_\psi[J]$ has an arc from $j$ to $i$ and $F_\psi$ has no arc from $I$ to $i$, then $H$ has an arc from $j$~to~$i$.}

\medskip
\noindent
Suppose that $F_\psi[J]$ has an arc from $j$ to $i$ and $F_\psi$ has no arc from $I$ to $i$. There is a configuration $x$ on $V$ such that $f_i(x)\neq f_i(x+e_j)$, and since $F_\psi$ has no arc from $I$ to $i$ we can choose $x$ in such a way that $x_I=b$. Then $h_i(x_J)=f_i(x)$ and since $(x+e_j)_I=b$ we have $h_i(x_J+e_j)=h_i((x+e_j)_J)=f_i(x+e_j)$. Hence $h_i(x_J)\neq h_i(x_J+e_j)$ thus $H$ has an arc from $j$ to $i$. This proves (5). 

\bigskip
\noindent
(6) {\em For all $s\in [m]$, if $\mu_s,c_s\in J$ then $(\mu_s,q_s),(q'_s,q_s),(q_s,c_s)$ are arcs of $H$.}

\medskip
\noindent
Suppose that $\mu_s,c_s\in J$. We deduce that $q'_s\in J$. Indeed, the unique in-neighbor of $q'_s$ is $\mu_s$, which is in $J$, and so from Lemma~\ref{lem:diffusion} we have $q'_s\in J$. So the two in-neighbors of $q_s$  are in $J$ and from the same lemma we deduce that $q_s\in J$. So $\mu_s,q'_s,q_s\in J$ hence $F_\psi$ has no arc from $I$ to $q_s$ and by (5) we deduce that $(\mu_s,q_s)$ and $(q'_s,q_s)$ are arcs of $H$. If $(q_s,c_s)$ is not in $H$, we deduce from (5) that $c_{s+1}\in I$, where $c_{m+1}$ means $\ell_n$, but then $c_s$ is a source of $H$ and this contradicts (2). Thus $(q_s,c_s)$ is in $H$. This proves (6).

\bigskip
\noindent
(7) {\em $\MU\subseteq I$.}

\medskip
\noindent
Suppose, for a contradiction, that some clause $\mu_s$ is in $H$, and let $s$ be minimal for that property. If $H$ has a cycle which does not contains $c_s$, then this cycle has necessarily a vertex in $\{\mu_1,\dots,\mu_{s-1}\}$, and this contradicts the choice of $s$. Thus every cycle of $H$ contains $c_s$. Thus $c_s$ is in $S$, and we deduce from (6) that $(\mu_s,q_s),(q'_s,q_s),(q_s,c_s)$ are arcs of $H$. Hence $q_s$ is in $S$  and has two in-neighbors in $S$. Thus $S$ is not a cycle, and this contradicts (4). This proves (7).

\bigskip
By (4) and (7), $S$ is a cycle disjoint from $\MU$. We deduce that $S$ contains each vertex in $\ELL\cup\C$ and exactly one of $\lambda^+_r,\lambda^-_r$ for each $r\in [n]$. Let $z\in\B^n$ defined by $z_r=1$ if $\lambda^+_r$ is not in $S$ and $z_r=0$ if $\lambda^-_r$ is not in $S$; there is no ambiguity since exactly one of $\lambda^+_r,\lambda^-_r$ is in $S$. By (7) and Lemma~\ref{lem:diffusion}, for each $s\in [m]$, $G[I]$ has a path from $t$ to $\mu_s$, and thus each clause $\mu_s$ contains a literals in $I$. This literal is not in $C$, thus it is satisfied by $z$. So $\psi$ is satisfied by $z$, a contradiction. We deduce that if $\psi$ is not satisfiable, then every BN on $F_\psi$ is synchronizing. 

\bigskip
\noindent
\BF{($2 \Rightarrow 1$)} Suppose that $z$ is a satisfying assignment of $\psi$, and let us prove that some BN on $F_\psi$ is not synchronizing. Let $\LAMBDA^1$ be the literals satisfied by $z$ and $\LAMBDA^0$ the literals not satisfied by $z$, hence $(\LAMBDA^0,\LAMBDA^1)$ is a balanced partition of $\LAMBDA^+\cup\LAMBDA^-$. Let $f$ be any and-or-net on $F_\psi$ such that:
\begin{itemize}
\item
$f_t$ is the $1$-constant function,
\item
for all $i\in\LAMBDA^1\cup\MU\cup\MU'$, $f_i$ is a disjunction,
\item
for all $i\in \ELL\cup\LAMBDA^0\cup \C$, $f_i$ is a conjunction.
\end{itemize}

\medskip
Let $\MU'^1$ be the vertices in $\MU'$ with at least one in-neighbors in $\LAMBDA^1$. Let $\q=\{q_1,\dots,q_m\}$ and $\q'=\{q'_1,\dots,q'_m\}$. Let $I=\{t\}\cup \LAMBDA^1\cup\MU'^1\cup \MU\cup\q\cup\q'$, and let $X$ be the set of configurations $x$ on the vertex set of $F_\psi$ with $x_I=\ONE$.

\bigskip
\noindent
(1) {\em $f(X)\subseteq X$.}

\medskip
\noindent
Let $x\in X$. We have $f_t(x)=1$. For every $i\in\LAMBDA^1$, $f_i$ is a disjunction and since $x_t=1$ we have $f_i(x)=1$. For every $i\in\MU'^1$, $f_i$ is a disjunction and $i$ has an in-neighbor $j\in\LAMBDA^1$, and since $x_j=1$ we deduce that $f_i(x)=1$. Next, consider a clause $\mu_s$ in $\MU$. Then $f_{\mu_s}$ is a disjunction and, since $z$ is a satisfying assignment, $\mu_s$ contains a literal $i\in\LAMBDA^1$ and thus $x_i=1$. If $i$ is an in-neighbor of $\mu_s$ then $f_{\mu_s}(x)=1$. Otherwise, $i$ is an in-neighbor of $\mu'_s$, so $\mu'_s\in\MU'^1$. Thus $x_{\mu'_s}=1$ and we deduce that $f_{\mu_s}(x)=1$. Finally every $i\in\q\cup \q'$ has only in-neighbors in $\q'\cup\MU\subseteq I$, and we deduce that $f_i(x)=1$. Thus $f(x)_I=\ONE$, that is, $f(x)\in X$. This proves (1). 

\bigskip
Let $J=\ELL\cup\LAMBDA^0\cup\C$. So $(I,J)$ is a partition of the vertices of $F_\psi$. Note $F_\psi[J]$ is a negative cycle since exactly one of $\lambda^+_r,\lambda^-_r$ is in $\LAMBDA^0$ for each $r\in [n]$.   

\bigskip
\noindent
(2) {\em Let $i\in J$ and let $j$ its in-neighbor in $F_\psi[J]$. For all $x,y\in X$, if $x_j\neq y_j$ then $f_i(x)\neq f_i(y)$.}

\medskip
\noindent
Let $x,y\in X$ with $x_j\neq y_j$. If $i=\ell_0$ then $j=c_1$ and $f_i(x)=\neg x_j\neq \neg y_j=f_i(y)$. Otherwise, $j$ is a positive in-neighbor of $i$. Furthermore, $i$ has exactly one in-neighbor $k\neq j$, which is positive and belongs to $\{t\}\cup\LAMBDA^1\cup\q\subseteq I$. Hence $x_k=y_k=1$, and since $f_i$ is a conjunction, we obtain $f_i(x)=x_j\land x_k=x_j\land 1=x_j$ and $f_i(y)=y_j\land y_k=y_j\land 1=y_j$. Thus $f_i(x)\neq f_i(y)$. This proves~(2).

\bigskip
Let $x$ and $y$ be $J$-opposite configurations in $X$, that is, $x,y\in X$ and $x_i\neq y_i$ for all $i\in J$. Let $i\in I$. By (1) we have $f(x),f(y)\in X$ and we deduce that $f^i(x)=x$ and $f^i(y)=y$, thus $f^i(x)$ and $f^i(y)$ are $J$-opposite configurations in $X$. Let $i\in J$. Obviously,  $f^i(x),f^i(y)\in X$. Furthermore, by (2) we have $f_i(x)\neq f_i(y)$, thus $f^i(x)$ and $f^i(y)$ are $J$-opposite configurations in $X$. So for any vertex $i$ in $F_\psi$, $f^i(x)$ and $f^i(y)$ are $J$-opposite configurations in $X$. This implies that $f^w(x)$ and $f^w(y)$ are $J$-opposite configurations in $X$ for any word $w$, and thus $f$ is not synchronizing. 
\end{proof}

%%%%%%%%%%%%%%%%%%%%%%%%%%%%%%%%%%%%%%%%%%%%%%%%%%%%%%%%%%%%%%%%%%%%%%%%%%%%%%%%%%%%%%%%%%%%%%%%%%%%%%%
\section{Concluding remarks}\label{sec:conclu}
%%%%%%%%%%%%%%%%%%%%%%%%%%%%%%%%%%%%%%%%%%%%%%%%%%%%%%%%%%%%%%%%%%%%%%%%%%%%%%%%%%%%%%%%%%%%%%%%%%%%%%%

As said in the introduction (and proved in Appendix \ref{app:one_fixed_point}), if $G$ has no negative cycles, every synchronizing BN on $G$ has a synchronizing word of length $n$: this solves (in a strong form) \v{C}ern\'y conjecture when $G$ has no negative cycles. It would be nice to prove the conjecture when $G$ has no positive cycles. In this paper, we were only able to do that for and-or-nets in the strong case: if $G$ has no positive cycles and is strongly connected, every synchronizing and-or-net on $G$ has a synchronizing word of length at most $10(\sqrt{5}+1)^n$. Another perspective is to prove the \v{C}ern\'y conjecture for and-or-nets (independently of the interaction digraph).  

\medskip
We feel that the upper bound $10(\sqrt{5}+1)^n$ can be widely improved, and we were not able to find any non-trivial lower bound. Besides, it seems possible to have better expressions, using some features of $G$ instead of the number $n$ of vertices only. For example, let $\tau(G)$ be the minimum size of a feedback vertex set of $G$. Using Lemma~\ref{lem:acyclic_image}, one can easily replace the bound $10(\sqrt{5}+1)^n$ by the better bound $\tau(G)+3F(n+4)2^{\tau(G)}$. However, for the particular case $\tau(G)\leq 1$, the following gives something much better (see Appendix \ref{app:tau1} for the proof). 

\begin{proposition}\label{pro:tau1}
Let $G$ be a strongly connected signed digraph on $[n]$ without positive cycles, which is not a cycle. If $\tau(G)\leq 1$ then there is a word $w$ of length $3n-1$ which synchronizes every and-or-net on $G$.
\end{proposition}

Following the approach initiated in \cite{AGRS20}, it would also be interesting to study the length of words that synchronize families of synchronizing BNs. The previous proposition gives an example, where the family is the set of and-or-nets on $G$. Let us discuss another bigger family, namely, the families $\mathcal{F}^+(n)$ of synchronizing BNs with component set $[n]$ and whose signed interaction digraphs have no negative cycles. Let $\ell^+(n)$ be the minimum length of a word synchronizing all the members of $\mathcal{F}^+(n)$. Since each member of $\mathcal{F}^+(n)$ has a synchronizing word of length $n$, considering the concatenation of these words, we obtain $\ell^+(n)\leq n|\mathcal{F}^+(n)|$, which is doubly exponential in $n$, since $\mathcal{F}^+(n)$ is. But this is far from the true value of $\ell^+(n)$, which is only quadratic in $n$ (see Appendix \ref{app:one_fixed_point} for a proof). 

\begin{proposition}\label{pro:l+}
$\ell^+(n)=n^2-o(n^2)$. 
\end{proposition}

By symmetry, it would be interesting to study the minimum length $\ell^-(n)$ of a word that synchronizes the families of synchronizing BNs with component set $[n]$ whose signed interaction digraphs have no positive cycles. 

%%%%%%%%%%%%%%%%%%%%%%%%%%%%%%%%%%%%%%%%%%%%%%%%%%%%%%%%%%%%%%%%%%%%%%%%%%%%%%%%
\paragraph{Acknowledgments}
%%%%%%%%%%%%%%%%%%%%%%%%%%%%%%%%%%%%%%%%%%%%%%%%%%%%%%%%%%%%%%%%%%%%%%%%%%%%%%%%

Julio Aracena was supported by ANID–Chile through Centro de Modelamiento Matemático (CMM), ACE210010 and FB210005, BASAL funds for center of excellence from ANID-Chile. Julio Aracena, Lilian Salinas and Adrien Richard were supported by ANID–Chile through ECOS C19E02. Adrien Richard was supported by 
the \emph{Young Researcher} project ANR-18-CE40-0002-01 ``FANs''.

\appendix

%%%%%%%%%%%%%%%%%%%%%%%%%%%%%%%%%%%%%%%%%%%%%%%%%%%%%%%%%%%%%%%%%%%%%%%%%%%%%%%%%%%%%%%%%%%%%%%%%%%%%%%
\section{Synchronizing BNs with a unique fixed point}\label{app:one_fixed_point}
%%%%%%%%%%%%%%%%%%%%%%%%%%%%%%%%%%%%%%%%%%%%%%%%%%%%%%%%%%%%%%%%%%%%%%%%%%%%%%%%%%%%%%%%%%%%%%%%%%%%%%%

In this appendix, we prove the two results given in the introduction before the statement of Theorem~\ref{thm:main}, which are Propositions \ref{pro:positive_G} and \ref{pro:negative_G} below. We then prove Proposition \ref{pro:l+} of Section~\ref{sec:conclu}. 

\medskip
The first tool is a classical result of Aracena, already mentioned. 

\begin{theorem}[Aracena \cite{A08}]\label{thm:A08}
Let $G$ be a signed digraph, and let $f$ be a BN on $G$.
\begin{enumerate}
\item If $G$ has no negative cycles, then $f$ has at least one fixed point. 
\item If $G$ has no positive cycles, then $f$ has at most one fixed point. 
\item If $G$ has no negative cycles and no sources, then $f$ has at least two fixed points. 
\item If $G$ has no positive cycles and no sources, then $f$ has no fixed points. 
\end{enumerate}
\end{theorem}

The second tool is the following. Given two words $u,v$, we say that $u$ \EM{contains} $v$ if $v$ can be obtained by deleting some letters in $u$. For instance, $aba$ contains $aa$. 

\begin{proposition}\label{pro:same_sign}
Let $G$ be a signed digraph with vertex set $V$ such that all the cycles of $G$ have the same sign. Let $f$ be a BN on $G$ with a unique fixed point. There is a permutation $\pi$ of $V$ such that any word containing $\pi$ synchronizes $f$. In particular, $\pi$ synchronizes $f$.
\end{proposition}

\begin{proof}
We proceed by induction on $|V|$. For $|V|=1$, since $f$ has a unique fixed point, $f$ is a constant function and the result is obvious. Suppose that $|V|\geq 2$. Since all the cycles of $G$ have the same sign and $f$ has a unique fixed point, we deduce from Theorem~\ref{thm:A08} that $G$ has a source, say $i$. Hence, $f_i$ is the $a$-constant function for some $a\in\B$. Let $I=V\setminus\{i\}$, where $V$ is the vertex set of $G$. Let $h$ be the BN with component set $I$ defined by $h(x_I)=f(x)_I$ for all configurations $x$ on $V$ with $x_i=a$. Since $f$ has a unique fixed point, $h$ has a unique fixed point, and since $f_i$ is the $a$-constant function, for every configuration $x$ on $V$ and every word $v$ we have:
\[\tag{1}
x_i=a~\Rightarrow~ f^v(x)_I=h^v(x_I)\textrm{ and }f^v(x)_i=a.
\]
Let $H$ be the signed interaction digraph of $G$. Since $H$ is a subgraph of $G$, all the cycles of $H$ have the same sign. Hence, by induction hypothesis, there is a permutation $\sigma$ of $I$ such that any word containing $\sigma$ synchronizes $h$. Let $\pi=i,\sigma$, and let $w$ be a word containing $\pi$. Since $w$ contains $i$, there is $u,v$ such that $w=u,i,v$ where $u$ does not contains $i$. So $v$ contains $\sigma$, thus $v$ synchronizes $h$. Let $x,y$ be any two configurations on $V$, and let $x'=f^{u,i}(x)$ and $y'=f^{u,i}(y)$. Since $f_i$ is the $a$-constant function, we have $x'_i=y'_i=a$, and since $v$ synchronizes $h$ we deduce from (1) that $f^v(x')_I=h^v(x'_I)=h^v(y'_I)=f^v(y')_I$ and $f^v(x')_i=f^v(y')_i=a$. Thus $f^v(x')=f^v(y')$ and this proves that $w$ synchronizes $f$. 
\end{proof}

We deduce the following two propositions. 

\begin{proposition}\label{pro:positive_G}
Let $G$ be an $n$-vertex signed digraph without negative cycles, and let $f$ be a BN on $G$. The following conditions are equivalent:
\begin{enumerate}
\item
$f$ has a unique fixed point.
\item
$f$ has a synchronizing word of length $n$. 
\item
$f$ is synchronizing.
\end{enumerate}
\end{proposition}

\begin{proof}
\BF{($1 \Rightarrow 2$)} is given by Proposition~\ref{pro:same_sign} and \BF{($2 \Rightarrow 3$)} is obvious. Suppose that $f$ is synchronizing. Then it has at most one fixed point. By Theorem~\ref{thm:A08}, $f$ has at least one fixed point. Thus $f$ has a unique fixed point. This proves \BF{($3 \Rightarrow 1$)}.
\end{proof}

\begin{proposition}\label{pro:negative_G}
Let $G$ be an $n$-vertex signed digraph without positive cycles, and let $f$ be a BN on $G$ with a fixed point. Then $f$ has a synchronizing word of length $n$. 
\end{proposition}

\begin{proof}
Since $f$ has a fixed point, by Theorem~\ref{thm:A08} it has a unique fixed point. We then deduce from Proposition~\ref{pro:same_sign} that $f$ has a synchronizing word of length $n$. 
\end{proof}

We now prove Proposition \ref{pro:l+} of Section~\ref{sec:conclu}.  Let us say that a word synchronizes a set of BNs if it synchronizes every BN in this set. Let $\mathcal{A}(n)$ be the set of BNs with component set $[n]$ and with an acyclic interaction digraph. A word $w$ is \EM{$n$-complete} if it contains all the permutations of $[n]$. Given a word $w$ and a BN $f$ with component set $V$, we say that $w$ \EM{fixes} $f$ if $f^w(x)$ is a fixed point of $f$ for all configurations $x$ on $V$. This notion was introduced in \cite{AGRS20}, where it is proved that a word fixes every BN in $\mathcal{A}(n)$ if and only if it is $n$-complete. It is clear that if $f$ has a unique fixed point, then $w$ fixes $f$ if and only if $w$ synchronizes $f$. Also, by Theorem~\ref{thm:A08}, every BN in $\mathcal{A}(n)$ has a unique fixed point. Consequently, we can restate the result of \cite{AGRS20} mentioned above as follows. 

\begin{theorem}[\cite{AGRS20}]\label{thm:n-complete}
A word synchronizes $\mathcal{A}(n)$ if and only if it is $n$-complete. 
\end{theorem}

Let $\mathcal{F}^+(n)$ be the set of synchronizing BNs with component set $[n]$ and whose signed interaction digraphs have no negative cycles. Let $\mathcal{F}(n)$ be the set of $f$ BN with component set $[n]$, with a unique fixed point, and such that all the cycles of the signed interaction digraph of $f$ have the same sign. Note that $\mathcal{A}(n)\subseteq \mathcal{F}^+(n)\subseteq \mathcal{F}(n)$; the first inclusion is obvious and the second follows from Proposition~\ref{pro:positive_G}. 

\begin{proposition}
The following conditions are equivalent:
\begin{enumerate}
\item $w$ is $n$-complete.
\item $w$ synchronizes $\mathcal{F}(n)$.
\item $w$ synchronizes $\mathcal{F}^+(n)$.
\item $w$ synchronizes $\mathcal{A}(n)$.
\end{enumerate}
\end{proposition}

\begin{proof}
\BF{($1 \Rightarrow 2$)} is given by Proposition~\ref{pro:same_sign}. \BF{($2 \Rightarrow 3$)} holds since $\mathcal{F}^+(n)\subseteq \mathcal{F}(n)$. \BF{($3 \Rightarrow 4$)} holds since $\mathcal{A}(n)\subseteq \mathcal{F}^+(n)$. \BF{($4 \Rightarrow 1$)} is given by Theorem~\ref{thm:n-complete}.
\end{proof} 

Let $\lambda(n)$ be the shortest length of an $n$-complete word, and let $\ell^+(n)$ the minimum length of a word that synchronizes $\mathcal{F}^+(n)$. By the previous proposition, $\ell^+(n)=\lambda(n)$, and it is proved in \cite{KK76} that $\lambda(n)=n^2-o(n^2)$. We thus obtain Proposition~\ref{pro:l+}. 

%%%%%%%%%%%%%%%%%%%%%%%%%%%%%%%%%%%%%%%%%%%%%%%%%%%%%%%%%%%%%%%%%%%%%%%%%%%%%%%%%%%%%%%%%%%%%%%%%%%%%%%
\section{Proof of Proposition~\ref{pro:tau1}}\label{app:tau1}
%%%%%%%%%%%%%%%%%%%%%%%%%%%%%%%%%%%%%%%%%%%%%%%%%%%%%%%%%%%%%%%%%%%%%%%%%%%%%%%%%%%%%%%%%%%%%%%%%%%%%%%

We need the following two lemmas. 

\begin{lemma}\label{lem:tau11}
Let $G$ be a strong signed digraph without positive cycles. Suppose that $G$ has a vertex $i$ that meets every cycle of $G$. Then $G$ is switch equivalent to a signed digraph $H$ in which all the in-coming arcs of $i$ are negative and all the other arcs are positive. 
\end{lemma}

\begin{proof}
Let $V$ be the vertex set of $G$. Let $G'$ be obtained by changing the sign of every in-coming arc of $i$. Then all the cycles of $G'$ are positive. By Proposition \ref{pro:harary}, the $I$-switch $H'$ of $G'$ is full positive for some $I\subseteq V$. Let $H$ by obtained from $H'$ by making negative the in-coming arcs of~$i$. One easily check that $H$ is the $I$-switch of $G$ and thus $H$ has the desired properties.
\end{proof}

\begin{lemma}\label{lem:tau12}
Let $G$ be a $n$-vertex strong signed digraph, which is not a cycle. Suppose that $G$ has a vertex $i$ such that: $i$ meets every cycle; all the in-coming arcs of $G$ are negative and all the other arcs are positive; $i$ is of in-degree at least two. There is a word $w$ of length $3n-1$ which synchronizes every BN $f$ on $G$ such that $f_i$ is a conjunction or a disjunction. 
\end{lemma}

\begin{proof}
Let $V$ be the vertex set of $G$, and let $u$ be a topological sort of $G\setminus i$.

\bigskip
\noindent
(1) {\em For every configuration $x$ on $V$ and prefix $v$ of $u$, we have $f^v(x)_j=x_i$ for all $j\in\{v\}\cup\{i\}$.}

\medskip
\noindent
We proceed by induction on $|v|$. This is obvious for $|v|=0$. Suppose that $|v|\geq 1$. Let $v'$ be the prefix of $u$ of length $|v|-1$, thus $v=v',k$ for some $k\in V$. Let $y=f^{v'}(x)$. For all $j\in\{v'\}\cup\{i\}$, we have $f^v(x)_j=y_j=x_i$, where the induction is used for the second equality. Since all the in-neighbors of $k$ are in $\{v'\}\cup\{i\}$ and positive, we deduce that $f_k(y)=x_i$, and since $f^v(x)_k=f_k(y)=x_i$ this completes the induction. This proves~(1).   

\bigskip
Let $v=\epsilon$ if $i$ has a loop, and let $v$ be the shortest prefix of $u$ containing an in-neighbor of $i$ otherwise. Note that $|v|\leq n-1$ since $|u|=n-1$. Let 
\[
w=u,i,v,i,u,
\]
of length $2n+|v|\leq 3n-1$. Let $f$ be a BN on $G$ such that $f_i$ is a conjunction or a disjunction. We will prove that $w$ synchronizes $f$. Let $x$ be any configuration on $V$. Let $y$ the configuration on $V$ such that $y_j=x_i$ for all $j\in V$. Let $I=\{v\}\cup\{i\}$, which contains an in-neighbor of $i$, and $z=y+e_I$. By (1) we have $f^u(x)=y$. Since $i$ has only negative in-neighbors, $f^i(y)=y+e_i$. We then deduce from (1) that $f^v(y+e_i)=y+e_I=z$. Let $j$ be an in-neighbor of $i$ in $I$. Since $i$ if of in-degree at least two, it has an in-neighbor $k\not\in I$. Then $z_j=y_j+1=x_i+1$ and $z_k=y_k=x_i$. Thus $z_j\neq z_k$ and we deduce that $f(z)_i=0$ if $f_i$ is a conjunction, and $f(z)_i=1$ if $f_i$ is a disjunction. We deduce from (1) that $f^{i,u}(z)=\ZERO$ if $f_i$ is a conjunction, and $f^{i,u}(z)=\ONE$ if $f_i$ is a disjunction. Since $f^{u,i,v}(x)=z$ we have $f^w(x)=f^{i,u}(z)$. Hence $f^w(x)=\ZERO$ if $f_i$ is a conjunction and $f^w(x)=\ONE$ if $f_i$ is a disjunction. Thus $w$ synchronizes~$f$. 
\end{proof}

We are now ready to prove Proposition~\ref{pro:tau1}, that we restate. 

\setcounter{proposition}{4}
\begin{proposition}
Let $G$ be a strongly connected signed digraph on $[n]$ without positive cycles, which is not a cycle. If $\tau(G)\leq 1$ then there is a word $w$ of length $3n-1$ which synchronizes every and-or-net on $G$.
\end{proposition}
\setcounter{proposition}{10}

\begin{proof}
If $\tau(G)=0$, then, by Lemma \ref{lem:acyclic_image}, there is a word of length $n$ which synchronizes every BN on $G$ and we are done. Suppose that $\tau(G)=1$. Since $G$ is not a cycle, it has a vertex $i$ of in-degree at least two that meets every cycle of $G$. By Lemma~\ref{lem:tau11}, $G$ is switch equivalent to a signed digraph $H$ which satisfies the conditions of Lemma~\ref{lem:tau12}. By this lemma, there is a word $w$ of length at most $3n-1$ which synchronizes every and-or-net on $H$. We then deduce from Proposition~\ref{pro:BN_switch} that $w$ synchronizes every and-or-net on $G$.   
\end{proof}

%%%%%%%%%%%%%%%%%%%%%%%%%%%%%%%%%%%%%%%%%%%%%%%%%%%%%%%%%%%%%%%%%%%%%%%%%%%%%%%%%%%%%%%%%%%%%%%%%%%%%%%
\section{A result mentioned in Remark \ref{rem:att_size}}\label{app:att_size}
%%%%%%%%%%%%%%%%%%%%%%%%%%%%%%%%%%%%%%%%%%%%%%%%%%%%%%%%%%%%%%%%%%%%%%%%%%%%%%%%%%%%%%%%%%%%%%%%%%%%%%%

\begin{proposition}
Let $X\subseteq \B^n$ such that, for all $i\in [n]$, there is $x,y\in X$ that only differ in component $i$. We have $|X|\geq n+1$.
\end{proposition}

\begin{proof}
We proceed by induction on $n$. The case $n=1$ is obvious. Suppose that $n\geq 2$ and let $I=[n-1]$. For $a=0,1$, let $X^a=\{x_I\mid x\in X, x_n=a\}$. For every $i\in I$, there is $x,y\in X$ that only differ in $i$, so $x_I$ and $y_I$ only differ in $i$, and $x_I,y_I\in X^0\cup X^1$. Since $X^0\cup X^1\subseteq \B^{n-1}$, by induction hypothesis, 
\[
n\leq |X^0\cup X^1|=|X^0|+|X^1|-|X^0\cap X^1|=|X|-|X^0\cap X^1|.
\]
Let $x,y\in X$ that only differ in $n$, say $x_n<y_n$. Then $x_I\in X^0$ and $y_I\in X^1$, and since $x_I=y_I$ we deduce that $|X^0\cap X^1|\geq 1$ and thus $|X|\geq n+1$. 
\end{proof}

%%%%%%%%%%%%%%%%%%%%%%%%%
\bibliographystyle{plain}
\bibliography{BIB}
%%%%%%%%%%%%%%%%%%%%%%%%%

\end{document}